\documentclass[11pt]{article}
\usepackage{amsfonts,amssymb,amsopn,amsmath,mathrsfs,theorem}
\usepackage{hyperref}
\usepackage{graphicx}
\usepackage{verbatim}
\usepackage{authblk}
\usepackage{color}
\usepackage{enumitem}
\usepackage[normalem]{ulem}
\usepackage{bm}

\newcommand{\detail}[1]{\emph{To be filled in.}}

\usepackage{tikz}
\usetikzlibrary{calc,arrows}

\newcounter{thm_counter}
\numberwithin{thm_counter}{section}

\newtheorem{lemma}[thm_counter]{Lemma}
\newtheorem{prop}[thm_counter]{Proposition}
\newtheorem{theorem}[thm_counter]{Theorem}

\newtheorem{cor}[thm_counter]{Corollary}

\theorembodyfont{\upshape}
\newtheorem{remark}[thm_counter]{Remark}

\numberwithin{equation}{section} 

\newcommand{\dis}{\displaystyle}

\newcommand{\noi}{\noindent}
\newcommand{\halmos}{\rule{1ex}{1.4ex}}
\newcommand{\QED}{\nopagebreak{\hspace*{\fill}$\halmos$\medskip}}

\newcommand{\med}{\medskip}
\newcommand{\quand}{\quad\mbox{and}\quad}

\newcommand{\bt}{\begin{theorem}}
\newcommand{\et}{\end{theorem}}
\newcommand{\bl}{\begin{lemma}}
\newcommand{\el}{\end{lemma}}
\newcommand{\bp}{\begin{prop}}
\newcommand{\ep}{\end{prop}}
\newcommand{\bcor}{\begin{cor}}
\newcommand{\ecor}{\end{cor}}
\newcommand{\br}{\begin{remark}}
\newcommand{\er}{\end{remark}}
\newcommand{\bcon}{\begin{conjecture}}
\newcommand{\econ}{\end{conjecture}}

\newenvironment{Proof}[1][]{\noi\textbf{Proof #1}}{\QED}
\newcommand{\bpro}{\begin{Proof}}
\newcommand{\epro}{\end{Proof}}

\newcommand{\be}{\begin{equation}}
\newcommand{\ee}{\end{equation}}
\newcommand{\ba}{\begin{array}}
\newcommand{\ea}{\end{array}}
\newcommand{\bc}{\be\begin{array}{r@{\,}c@{\,}l}}
\newcommand{\bac}{\begin{array}{r@{\,}c@{\,}l}}
\newcommand{\ec}{\end{array}\ee}

\newcommand{\ga}{\gamma}
\newcommand{\Ga}{\Gamma}
\newcommand{\de}{\delta}

\newcommand{\eps}{\varepsilon}
\newcommand{\la}{\lambda}
\newcommand{\La}{\Lambda}
\newcommand{\sig}{\sigma}

\newcommand{\Ai}{{\cal A}}
\newcommand{\Bi}{{\cal B}}
\newcommand{\Ci}{{\cal C}}
\newcommand{\Di}{{\cal D}}

\newcommand{\Fi}{{\cal F}}
\newcommand{\Gi}{{\cal G}}
\newcommand{\Hi}{{\cal H}}
\newcommand{\Ii}{{\cal I}}

\newcommand{\Ki}{{\cal K}}
\newcommand{\Li}{{\cal L}}

\newcommand{\Oi}{{\cal O}}

\newcommand{\Ri}{{\cal R}}

\newcommand{\Xc}{{\cal X}}
\newcommand{\Yi}{{\cal Y}}

\newcommand{\li}{\langle}
\newcommand{\re}{\rangle}
\newcommand{\lli}{\langle\!\langle}
\newcommand{\rre}{\rangle\!\rangle}

\newcommand{\desd}{\ensuremath{\Leftrightarrow}}
\newcommand{\volgt}{\ensuremath{\Rightarrow}}
\newcommand{\up}{\uparrow}
\newcommand{\down}{\downarrow}

\newcommand{\sub}{\subset}
\newcommand{\beh}{\backslash}
\newcommand{\asto}[1]{\underset{{#1}\to\infty}{\longrightarrow}}

\newcommand{\ti}{\tilde}

\newcommand{\ov}{\overline}
\newcommand{\un}{\underline}
\newcommand{\subb}[2]{_{\ba{c}\scriptstyle{#1}\\[-.15cm]\scriptstyle{#2}\ea}}

\newcommand{\pa}{\partial}
\newcommand{\ffrac}[2]{{\textstyle\frac{{#1}}{{#2}}}}

\newcommand{\cn}{\colon}

\newcommand{\half}{{[0,\infty)}}

\def\to{\rightarrow} 

\def\mb{\mathbb} 

\def\R{\mb{R}} 

\def\N{\mb{N}}
\def\Q{\mb{Q}}
\def\Z{\mb{Z}}

\def\~{\sim}
\def\li{\langle}

\def\qed{$\blacksquare$}
\def\1{\mathbbm{1}}

\fontsize{11pt}{16.0pt}
\selectfont

\def\epsilon{\varepsilon} 

\newcommand{\Ob}{{\mathbf O}}

\newcommand{\sg}{{\mathfrak s}}
\newcommand{\binf}{\bm{\infty}}
\newcommand{\pre}{\preceq}
\newcommand{\vhi}{\varphi}

\newcommand{\lo}{(\hspace{-3pt}(}
\newcommand{\ro}{)\hspace{-3pt})}
\newcommand{\lc}{[\![}
\newcommand{\rc}{]\!]}

\begin{document}

\hyphenation{equip-ped}
\hyphenation{re-de-rive}

\makeatletter\@addtoreset{equation}{section}
\makeatother\def\theequation{\thesection.\arabic{equation}} 

\renewcommand{\labelenumi}{{\rm (\roman{enumi})}}
\renewcommand{\theenumi}{\roman{enumi}}

\allowdisplaybreaks

\author[1]{Nic Freeman\thanks{School of Mathematics and Statistics, Hicks Building, University Of Sheffield, S3 7RH, United Kingdom, n.p.freeman@sheffield.ac.uk}}
\author[2]{Jan M.\ Swart\thanks{Institute of Information Theory and Automation, Pod vod\'arenskou v\v{e}\v{z}\'i 4, 18200 Praha 8, Czech republic, swart@utia.cas.cz}}
\affil[1]{School of Mathematics and Statistics, University Of Sheffield}
\affil[2]{The Czech Academy of Sciences, Institute of Information Theory and Automation}

\title{Skorokhod's topologies on path space}
\date{\today}
\maketitle

\begin{abstract}
Skorokhod’s J1 and M1 topologies are standard tools in proving limit theorems for stochastic processes. Motivated by applications, we extend these topologies so that they are capable of describing the convergence of a sequence of functions that are not all defined on the same domain. Traditionally, the J1 and M1 topologies are defined using time changes. Instead, we base our definitions on the point of view that the graph of a cadlag function can naturally be viewed as a compact set that is equipped with a total order. The distance between two graphs is then measured by matching points on one graph with points on the other graph in a way that respects the total order. We treat the J1 and M1 topologies in a unified framework and simplify the existing theory. We introduce a space of paths, elements of which are cadlag functions defined on an arbitrary closed subset of the real line. We show that this space is Polish and derive compactness criteria. Specialised to functions that are all defined on the same domain, this yields new proofs of known results.
\end{abstract}

\vspace{.5cm}

\noi
{\it Acknowledgments.} Work sponsored by GA\v{C}R grant 22-12790S.

\newpage

\tableofcontents

\newpage

\section{Introduction and main results}\label{S:intro}

\subsection{Introduction}

A function is \emph{cadlag} (from the French ``continue \`a droit, limite \`a gauche'') if it is right continuous with left limits. In his classical paper \cite{Sko56}, Skorokhod introduced four topologies on the space of real valued cadlag functions on a compact interval, which he called J1, J2, M1, and M2. His definitions have later been generalised to functions defined on unbounded intervals and, for the J1 and J2 topologies, taking values in general metrisable spaces. When restricted to continuous functions, all four topologies reduce to the topology of locally uniform convergence. Of Skorokhod's four topologies, his J1 topology has proved to be the most natural in many situations, in particular, when discussing convergence of Markov processes \cite{EK86}. For this reason, Skorokhod's J1 topology is now normally known as ``the Skorokhod topology''. Classical textbook discussions of the J1 topology can be found in \cite[Section~3.5]{EK86} and \cite[Section~12]{Bil99}. All four topologies introduced by Skorokhod are discussed in \cite[Section~11.5]{Whi02}.

In some applications, such as excursion theory \cite{Rog89} or the theory of the Brownian web \cite{FINR04,SSS17}, one needs a space whose elements are functions that are not all defined on the same domain. In particular, a central object in the theory of the Brownian web is a space of \emph{continuous upward paths} that are continuous functions that are defined on an interval of the form $[s,\infty)$ with $s\in\ov\R:=[-\infty,\infty]$ and that take values in $\ov\R$. The authors of \cite{FINR04} equipped this space with a topology that can roughly be described as locally uniform convergence of functions together with convergence of their starting times. They proved that this space is Polish and derived a compactness criterion, both of which are important for proving weak convergence in law of random paths, in view of Prokhorov's theorem \cite[Thm~5.2]{Bil99}.

Inspired by this, we introduce the \emph{path space} over a general metrisable space. Elements of this space are \emph{paths}, which are pairs consisting of a closed subset of the real line and a cadlag function that is defined on that subset and takes values in the metrisable space. We equip the path space with topologies that generalise Skorokhod's J1 and M1 topologies, prove that these topologies are Polish, and derive compactness criteria. Our work is motivated by collections of cadlag paths that approximate the Brownian web \cite{NRS05,BGS15,EFS17} and more recent work on extensions of the Brownian web where the limiting objects also consist of cadlag paths \cite{MRV19,FS23}. In our follow-up paper \cite{FS25} we use the results of the present paper to derive tightness criteria for sequences of random compact sets of cadlag paths. In the remainder of Section~\ref{S:intro} we state our main results. Proofs are postponed till Sections \ref{S:cadfun}--\ref{S:mainproof}.

\subsection{Paths and cadlag functions}\label{S:paths}

Let $\ov\R:=[-\infty,\infty]$ be the extended real line. For any metrisable topological space $\Xc$, we set
\be
\Ri(\Xc):=(\Xc\times\R)\cup\big\{(\ast,-\infty),(\ast,\infty)\big\}.
\ee
We can think of $\Ri(\Xc)$ as being obtained from the product space $\Xc\times\ov\R$ by squeezing the sets $\Xc\times\{-\infty\}$ and $\Xc\times\{+\infty\}$ into single points, which we denote by $(\ast,-\infty)$ and $(\ast,+\infty)$. For this reason, we call $\Ri(\Xc)$ the \emph{squeezed space} associated with $\Xc$. We will always equip $\Ri(\Xc)$ with the topology described by the following lemma. We use the (nowadays fairly standard) convention of calling a topological space $\Xc$ \emph{Polish} if $\Xc$ is separable and there exists a complete metric generating the topology on $\Xc$.

\bl[Squeezed space]
Let\label{L:squeeze} $\Xc$ be a metrisable topological space. Then there exists a unique metrisable topology on $\Ri(\Xc)$ such that a sequence $(x_n,t_n)\in\Ri(\Xc)$ converges to a limit $(x,t)\in\Ri(\Xc)$ if and only if:
\[
{\rm(i)}\quad t_n\to t,\qquad{\rm(ii)}\quad\mbox{if $t\in\R$, then }x_n\to x.
\]
If $\Xc$ is Polish, then so is $\Ri(\Xc)$, and if $\Xc$ is compact, then so is $\Ri(\Xc)$.
\el

The squeezed space $\Ri(\ov\R)$ plays an important role in the theory of the Brownian web, see \cite[Figure~6.2]{SSS17}. The basic idea goes back to \cite[formula (3.4)]{FINR04} which defines a metric generating the topology on $\Ri(\ov\R)$.

We let $\Pi(\Xc)$ denote the space of all pairs $(\pi,\pre)$ where $\pi$ is a nonempty compact subset of $\Ri(\Xc)$ and $\pre$ is a total order on $\pi$, such that
\begin{enumerate}
\item $\pi^{\li 2\re}:=\big\{(z,z')\in\pi^2:z\pre z'\big\}$ is a closed subset of $\pi^2:=\pi\times\pi$, where we equip $\pi$ with the induced topology from $\Ri(\Xc)$ and $\pi^2$ with the product topology.
\item If $(x,s),(x',s')\in\pi$ satisfy $s<s'$, then $(x,s)\pre(x',s')$.
\item For each $t\in\R$, the set $\{x\in\Xc:(x,t)\in\pi\}$ has at most two elements.
\end{enumerate}
We call $\Pi(\Xc)$ the \emph{path space} and call its elements \emph{paths}. We will often be sloppy and simply denote a path by $\pi$. We call the closed sets $I(\pi)\sub\R$ and $\hat I(\pi)\sub\ov\R$ defined as
\be
I(\pi):=\hat I(\pi)\cap\R\quad\mbox{with}\quad
\hat I(\pi):=\big\{t\in\ov\R:(x,t)\in\pi\mbox{ for some }x\in\Xc\cup\{\ast\}\big\}
\ee
the \emph{domain} and \emph{extended domain} of $\pi$, respectively. For $t\in I(\pi)$ we uniquely define $\pi(t-)$ and $\pi(t+)$ by
\be\label{pit}
\big\{x\in\Xc:(x,t)\in\pi\big\}=\big\{\pi(t-),\pi(t+)\big\}\quad\mbox{with}\quad\pi(t-)\pre\pi(t+)\quad(t\in I(\pi)).
\ee
If $\pm\infty\in\hat I(\pi)$, then we also set $\pi(\pm\infty):=\ast$. We call
\be\label{Pc}
\Pi_{\rm c}(\Xc):=\big\{\pi\in\Pi(\Xc):\pi(t-)=\pi(t+)\ \forall t\in I(\pi)\big\},
\ee
the space of \emph{continuous} paths. For paths $\pi\in\Pi_{\rm c}$ we write $\pi(t):=\pi(t-)=\pi(t+)$ $(t\in I(\pi))$. We let
\be\label{Pint}
\Pi^|(\Xc):=\big\{\pi\in\Pi(\Xc):\hat I(\pi)\mbox{ is an interval}\big\}
\ee
and naturally write $\Pi^|_{\rm c}(\Xc):=\Pi^|(\Xc)\cap\Pi_{\rm c}(\Xc)$.

The notation $\pi(t-)$ and $\pi(t+)$ is suggestive of left and right limits. Indeed, there is a close connection between paths and cadlag functions, as we now explain. We adopt the following definition. A \emph{cadlag function} with values in a metrisable space $\Xc$ is a triple $(I,f^-,f^+)$ such that $I\sub\R$ is a closed set, $f^-\cn I\to\Xc$ is left continuous, $f^+\cn I\to\Xc$ is right continuous, and
\be\ba{ll}\label{cadlag}
\dis\lim_{s\up t}f^+(s)=f^-(t)\quad\forall t\in I^-,\qquad&\dis\lim_{s\down t}f^-(s)=f^+(t)\quad\forall t\in I^+,\\[5pt]
\dis I^-:=\big\{t\in I:(t-\eps,t)\cap I\neq\emptyset\ \forall\eps>0\big\}\quad&\dis I^+:=\big\{t\in I:(t,t+\eps)\cap I\neq\emptyset\ \forall\eps>0\big\}.
\ec
Note that in the special case that $I=\half$, this says that $f^+\cn I\to\Xc$ is right continuous with left limits and $f^-$ is its left continuous modification, except that we allow for the case that $f^-(0)\neq f^+(0)$, that is, in our formulation a cadlag function defined on $I=\half$ can jump at its initial time, contrary to the usual conventions. The following lemma explains the connection between cadlag functions and paths.

\bl[Cadlag functions]
If\label{L:fpi} $\pi\in\Pi(\Xc)$ is a path, then setting $I:=I(\pi)$ and $f^\pm(t):=\pi(t\pm)$ $(t\in I)$ defines a cadlag function $(I,f^-,f^+)$. Conversely, if $(I,f^-,f^+)$ is a cadlag function, then there exists a path $\pi\in\Pi(\Xc)$ such that $I=I(\pi)$ and $f^\pm(t)=\pi(t\pm)$ $(t\in I)$. This path can be chosen such that $\{-\infty,\infty\}\sub\hat I(\pi)$, and with this choice, it is unique.
\el

Lemma~\ref{L:fpi} can roughly be described by saying that paths are the graphs of cadlag functions. The total order $\pre$ on a path $\pi$ plays a crucial role here. Indeed, the set $\pi\sub\Ri(\R)$ defined as $\pi:=\{(0,t):-1\leq t\leq 1\}\cup\{(1,0)\}$ is compact and satisfies condition~(iii) of our definition of a path, but it is not possible to equip $\pi$ with a total order satisfying conditions (i) and (ii) and indeed $\pi$ is not the graph of a cadlag function.

Let $I\sub\R$ be a closed real interval of positive length and let $\pa I$ be the set of its finite endpoints. Then we let $\Di_I(\Xc)$ denote the space of pairs $(f^-,f^+)$ such that $(I,f^-,f^+)$ is a cadlag function with values in $\Xc$ and $f^-(t)=f^+(t)$ for all $t\in\pa I$. By Lemma~\ref{L:fpi}, we have the natural identification
\be\label{DiI}
\Di_I(\Xc)\cong\big\{\pi\in\Pi(\Xc):\hat I(\pi)=I\cup\{-\infty,\infty\},\ \pi(t-)=\pi(t+)\ \forall t\in\pa I\big\}.
\ee
An element $(f^-,f^+)$ of $\Di_I(\Xc)$ is uniquely determined by either of the functions $f^-$ or $f^+$. The convention is to choose $f^+$. Thus, we can also naturally identify $\Di_I(\Xc)$ with the space of right continuous functions $f\cn I\to\Xc$ such that
\be\ba{c}\label{DiI2}
\dis f(t-):=\lim_{s\up t}f(s)\mbox{ exists for all }t\in I\beh\{l\}
\quand
f(r-)=f(r)\quad\mbox{if }r<\infty,\\[5pt]
\dis\mbox{where}\quad l:=\inf I\mbox{ and }r:=\sup I.
\ec

\subsection{Skorokhod's J1 topology on path space}

A \emph{correspondence} between two sets $A,B$ is a set $C\sub A\times B$ such that
\be\label{corresp}
\forall a\in A\ \exists b\in B\mbox{ s.t.\ }(a,b)\in C
\quand
\forall b\in B\ \exists a\in A\mbox{ s.t.\ }(a,b)\in C.
\ee
We let ${\rm Corr}(A,B)$ denote the set of all correspondences between $A$ and $B$. If $\Xc$ is a metrisable topological space, then we let $\Ki_+(\Xc)$ denote the space of nonempty compact subsets of $\Xc$. If $d$ is a metric generating the topology on $\Xc$, then the corresponding \emph{Hausdorff metric} on $\Ki_+(\Xc)$ is defined as
\be\label{Haus}
d_{\rm H}(K_1,K_2):=\sup_{z_1\in K_1}d(z_1,K_2)\vee\sup_{z_2\in K_2}d(z_2,K_1)
=\inf_{C\in{\rm Corr}(K_1,K_2)}\sup_{(z_1,z_2)\in C}d(z_1,z_2),
\ee
where as usual $d(z,K):=\inf_{z'\in K}d(z,z')$ denotes the distance of a point to a set. We call the topology generated by $d_{\rm H}$ the \emph{Hausdorff topology} on $\Ki_+(\Xc)$. It is closely related to the more well-known Vietoris topology, see Lemma~\ref{L:Vietoris} below. The following facts are proved in \cite[Lemmas B.1--B.3]{SSS14}.

\bl[The Hausdorff topology]
The\label{L:Haus} topology on $\Ki_+(\Xc)$ generated by $d_{\rm H}$ only depends on the topology on $\Xc$ and not on the choice of the metric $d$. The space $\Ki_+(\Xc)$ is Polish if $\Xc$ is Polish and compact if $\Xc$ is compact.
\el

We fix a metric $d_{\rm sqz}$ generating the topology on the squeezed space $\Ri(\Xc)$ and let
\be\label{Hpi}
d_{\rm H}(\pi_1,\pi_2):=\inf_{C\in{\rm Corr}(\pi_1,\pi_2)}\sup_{(z_1,z_2)\in C}d_{\rm sqz}(z_1,z_2)\qquad\big(\pi_1,\pi_2\in\Pi(\Xc)\big)
\ee
denote the corresponding distance in Hausdorff metric between two paths $\pi_1,\pi_2\in\Pi(\Xc)$. Then $d_{\rm H}$ is a metric on the space $\Pi_{\rm c}(\Xc)$ of continuous paths defined in (\ref{Pc}), but only a pseudometric on the larger space $\Pi(\Xc)$, since in general to characterise a path one needs both the set $\pi$ and the total order $\pre$ on $\pi$. There are two natural ways to fix this, which we describe now.

For paths $\pi_1,\pi_2\in\Pi(\Xc)$, we let ${\rm Cor}_+(\pi_1,\pi_2)$ denote the set of correspondences $C$ between $\pi_1$ and $\pi_2$ that are \emph{monotone} in the sense that
\be\label{moncor}
\mbox{there are no $(z_1,z_2),(z'_1,z'_2)\in C$ such that $z_1\prec_1 z'_1$ and $z'_2\prec_2 z_2$},
\ee
where $z\prec z'$ means $z\pre z'$ and $z\neq z'$. If $d_{\rm sqz}$ is a metric generating the topology on the squeezed space $\Ri(\Xc)$, then we define a corresponding distance $d_{\rm tot}$ on the path space $\Pi(\Xc)$ by
\be\label{dpath}
d_{\rm tot}(\pi_1,\pi_2):=\inf_{C\in{\rm Cor}_+(\pi_1,\pi_2)}\sup_{(z_1,z_2)\in C}d_{\rm sqz}(z_1,z_2)\qquad\big(\pi_1,\pi_2\in\Pi(\Xc)\big),
\ee
which differs from (\ref{Hpi}) only in the fact that we restrict the infimum to monotone correspondences. We will shortly see that $d_{\rm tot}$ is, indeed, a metric.

Another natural way to define a metric on path space is as follows. First, we define a metric $d^2_{\rm sqz}$ generating the product topology on $\Ri(\Xc)^2:=\Ri(\Xc)\times\Ri(\Xc)$ by
\be
d^2_{\rm sqz}\big((z_1,z'_1),(z_2,z'_2)\big):=d_{\rm sqz}(z_1,z_2)\vee d_{\rm sqz}(z'_1,z'_2)\qquad\big(z_1,z'_1,z_2,z'_2\in\Ri(\Xc)\big),
\ee
and we let $d^2_{\rm H}$ denote the corresponding Hausdorff metric on $\Ki_+(\Ri(\Xc)^2)$. It is then easy to see that we can define a metric $d_{\rm part}$ on the path space $\Pi(\Xc)$ by setting
\be\label{dpath2}
d_{\rm part}(\pi_1,\pi_2):=d^2_{\rm H}(\pi^{\li 2\re}_1,\pi^{\li 2\re}_2)\qquad\big(\pi_1,\pi_2\in\Pi(\Xc)\big),
\ee
where $\pi^{\li 2\re}$ is defined below Lemma~\ref{L:squeeze}. Below is our first main result. We recall that the statement that $\Pi(\Xc)$ is Polish merely means that $\Pi(\Xc)$ is separable and there exists a complete metric generating the topology. This does not imply that the metrics $d_{\rm tot}$ and $d_{\rm part}$ are complete, which, indeed, they are not.

\bt[Skorokhod's J1 topology on path space]
Let\label{T:J1} $\Xc$ be a metrisable space and let $d_{\rm H}$, $d_{\rm tot}$, and $d_{\rm part}$ be defined as in (\ref{Hpi}), (\ref{dpath}), and (\ref{dpath2}) in terms of some metric $d_{\rm sqz}$ on the squeezed space $\Ri(\Xc)$. Then
\be\label{dcomp}
d_{\rm H}(\pi_1,\pi_2)\leq d_{\rm part}(\pi_1,\pi_2)\leq d_{\rm tot}(\pi_1,\pi_2)
\qquad\big(\pi_1,\pi_2\in\Pi(\Xc)\big).
\ee
The functions $d_{\rm tot}$ and $d_{\rm part}$ are metrics on $\Pi(\Xc)$ that generate the same topology. This topology does not depend on the choice of the metric $d_{\rm sqz}$ on the squeezed space $\Ri(\Xc)$. On the space $\Pi_{\rm c}(\Xc)$, all three metrics $d_{\rm H}$, $d_{\rm tot}$, and $d_{\rm part}$ generate the same topology. The space $\Pi^|_{\rm c}(\Xc)$ defined in (\ref{Pint}) is a closed subset of $\Pi(\Xc)$. If $\Xc$ is Polish, then so are $\Pi(\Xc)$, $\Pi_{\rm c}(\Xc)$, and the space $\Di_I(\R)$ defined in (\ref{DiI}). On $\Di_I(\Xc)$, the metrics $d_{\rm tot}$ and $d_{\rm part}$ generate Skorokhod's J1 topology while $d_{\rm H}$ is a metric that generates the J2 topology.
\et

We call the topology on $\Pi(\Xc)$ generated by $d_{\rm tot}$ or $d_{\rm part}$ the \emph{J1 topology} on $\Pi(\Xc)$. We will derive a compactness criterion for this topology. Before we do this, we discuss a more general class of topologies on path space which includes (a generalisation of) Skorokhod's M1 topology. This will allow us to formulate a compactness criterion that applies to all these topologies at once.

\subsection{M1 and other topologies}

There are great similarities between Skorokhod's J1 and M1 topologies. In fact, it is possible to treat them in a unified framework \cite{Pom76}. We will introduce a general class of topologies on path space that includes the J1 and M1 topologies. To do this, we introduce a natural concept that we will call ``betweenness'' and that seems to be new in this context. It seems quite conceivable it may have been invented in other contexts before, but we have been unable to find a reference. If $\Xc$ is a set, then we define a \emph{betweenness} on $\Xc$ to be a function that assigns to each pair $x,z$ of elements of $\Xc$ a subset $\li x,z\re$ of $\Xc$, such that the following axioms hold for all $x,y,z\in\Xc$:
\begin{enumerate}
\item $\li x,z\re=\li z,x\re$,\label{A1}
\item $x\in\li x,z\re$,\label{A2}
\item $y\in\li x,z\re\ \volgt\ \li x,y\re\cap\li y,z\re=\{y\}$,\label{A3}
\item $y\in\li x,z\re\ \volgt\ \li x,y\re\cup\li y,z\re=\li x,z\re$.\label{A4}
\end{enumerate}
If $y\in\li x,z\re$, then we say that $y$ lies \emph{between} $x$ and $z$. We call $\li x,z\re$ the \emph{segment} with \emph{endpoints} $x$ and $z$. The following lemma says that segments of a betweenness are equipped with a natural total order.

\bl[Total order on segments]
For\label{L:partdef} any betweennness on $\Xc$, for each $x,z\in\Xc$ and $y,y'\in\li x,z\re$, one has
\be\label{partdef}
\li x,y\re\sub\li x,y'\re\ \desd\ y\in\li x,y'\re\ \desd\ y'\in\li y,z\re\ \desd\ \li y,z\re\supset\li y',z\re.
\ee
Setting $y\leq_{x,z}y'$ if any of these equivalent conditions holds defines a total order on $\li x,z\re$.
\el

We next give some examples of betweennesses. For any set $\Xc$, it is straightforward to check that setting $\li x,z\re:=\{x,z\}$ defines a betweenness. We call this the \emph{trivial betweenness}. If $\Xc$ is a linear space, then it is easy to see that
\be\label{linbet}
\li x,z\re:=\big\{(1-p)x+pz:p\in[0,1]\big\}\qquad(x,z\in\Xc)
\ee
defines a betweenness on $\Xc$. We call this the \emph{linear betweenness}. If $(\Xc,\leq)$ is a totally ordered space, then one can check that setting
\be\label{ordbet}
\li x,z\re:=\big\{y\in\Xc:x\leq y\leq z\mbox{ or }z\leq y\leq x\big\}\qquad(x,z\in\Xc)
\ee
defines a betweenness on $\Xc$. We call this the \emph{order betweenness}. If $(\Xc,d)$ is a metric space, then we recall that a \emph{geodesic} in $(\Xc,d)$ is a subset $\Ga$ of $\Xc$ that is isometric to a compact real interval, i.e., there exists an isometry $\ga\cn[s,u]\to\Xc$ (with $s,u\in\R,\ s\leq u$) such that $\Ga$ is the image of $[s,u]$ under $\ga$. Clearly, $\ga$ is uniquely determined by $\Ga$ up to translations and mirror images of the interval $[s,u]$. The points $\ga(s),\ga(u)$ are called the \emph{endpoints} of the geodesic $\Ga$. We say that a metric space has \emph{unique geodesics} if for each $x,z\in\Xc$, there exists a unique geodesic $\Ga$ with endpoints $x,z$. We call the betweenness defined in the following lemma the \emph{geodesic betweennesss}.

\bl[Geodesic betweenness]
Let\label{L:metbet} $(\Xc,d)$ be a metric space with unique geodesics. Then letting $\li x,z\re$ denote the unique geodesic with endpoints $x,z$ defines a betweenness on $\Xc$.
\el

As an example of spaces without a linear structure where Lemma~\ref{L:metbet} is applicable we mention real-trees \cite{DT96}. If $\Xc$ is a metrisable space, then we say that a betweenness on $\Xc$ is \emph{proper} if $\li x,z\re$ is compact for each $x,z\in\Xc$, and the map $\Xc^2\ni(x,z)\mapsto\li x,z\re\in\Ki_+(\Xc)$ is continuous with respect to the product topology on $\Xc^2$ and the Hausdorff topology on $\Ki_+(\Xc)$. We recall that a metric space is called \emph{proper} if for each $x\in\Xc$ and $r\geq 0$, the closed ball $\{y\in\Xc:d(x,y)\leq r\}$ is a compact subset of $\Xc$. We do not know if the properness assumption in part~(iv) of the following lemma is needed, but it is certainly sufficient.

\bl[Proper betweennesses]
The\label{L:compbet} following betweennesses are proper:
\begin{enumerate}
\item The trivial betweenness on any metrisable space $\Xc$,
\item The linear betweenness on any normed linear space $\Xc$,
\item The order betweenness on any closed subset $\Xc$ of $\R$,
\item The geodesic betweennesss on any proper metric space $\Xc$.
\end{enumerate}
\el

If $\Xc$ is a metrisable space that is equipped with a betweenness, then for each path $\pi\in\Pi(\Xc)$ we define a \emph{filled path} by
\be\label{fillpath}
\ov\pi:=\big\{(x,t):t\in I(\pi),\ x\in\li\pi(t-),\pi(t+)\re\big\}\cup\big\{(\ast,t):t\in\hat I(\pi)\beh I(\pi)\big\}.
\ee
Note that for the trivial betweenness, $\ov\pi=\pi$. Also, for continuous paths $\pi\in\Pi_{\rm c}(\Xc)$ one has $\ov\pi=\pi$ regardless of the choice of the betweenness. We equip $\ov\pi$ with the total order defined by setting $(x,s)\pre(y,t)$ if either $s<t$ and $x,y$ are arbitrary, or $s=t\in\R$ and $x\leq_{\pi(t-),\pi(t+)}y$, where $\leq_{\pi(t-),\pi(t+)}$ is the total order on the segment $\li\pi(t-),\pi(t+)\re$ defined in Lemma~\ref{L:partdef}. Informally, the total order $\pre$ corresponds to the direction of time. We need the following simple lemma.

\bl[Compatibility with the topology]
Let\label{L:cogra} $\Xc$ be a metrisable space that is equipped with a proper betweenness. Then for each $\pi\in\Pi(\Xc)$, the set $\ov\pi^{\li 2\re}:=\big\{(z,z')\in\ov\pi^2:z\pre z'\big\}$ is a closed subset of $\ov\pi^2$.
\el

It is easy to see that a path $\pi$ is uniquely determined by the corresponding filled path $\ov\pi$. We generalise (\ref{Hpi}) by setting
\be\label{Hpi2}
d_{\rm H}(\pi_1,\pi_2):=\inf_{C\in{\rm Corr}(\ov\pi_1,\ov\pi_2)}\sup_{(z_1,z_2)\in C}d_{\rm sqz}(z_1,z_2).
\ee
In general, this is only a pseudometric. We define monotone correspondences between filled paths in the same was as for paths. Generalising (\ref{dpath}) and (\ref{dpath2}), we define metrics on $\Pi(\Xc)$ by
\bc\label{dSkor}
\dis d_{\rm tot}(\pi_1,\pi_2)&:=&\dis\inf_{C\in{\rm Cor}_+(\ov\pi_1,\ov\pi_2)}\sup_{(z_1,z_2)\in C}d_{\rm sqz}(z_1,z_2),\\[5pt]
\dis d_{\rm part}(\pi_1,\pi_2)&:=&\dis d^2_{\rm H}(\ov\pi^{\li 2\re}_1,\ov\pi^{\li 2\re}_2).
\ec
In particular, for the trivial betweenness, these definitions coincide with the definitions in (\ref{dpath}) and (\ref{dpath2}). The following theorem is our second main result. For the trivial betweenness, the statements of the theorem have already been proved in Theorem~\ref{T:J1}.


\bt[Skorokhod-type topologies on path space]
Assume\label{T:M1} that $\Xc$ is a metrisable space that is equipped with a proper betweenness. Let $d_{\rm H}$, $d_{\rm tot}$, and $d_{\rm part}$ be defined as in (\ref{Hpi2}) and (\ref{dSkor}). Then
\be\label{dcomp2}
d_{\rm H}(\pi_1,\pi_2)\leq d_{\rm part}(\pi_1,\pi_2)\leq d_{\rm tot}(\pi_1,\pi_2)
\qquad\big(\pi_1,\pi_2\in\Pi(\Xc)\big).
\ee
The functions $d_{\rm tot}$ and $d_{\rm part}$ are metrics on $\Pi(\Xc)$ that generate the same topology. This topology does not depend on the choice of the metric $d_{\rm sqz}$ on the squeezed space $\Ri(\Xc)$. If $\Xc$ is Polish, then so is $\Pi(\Xc)$ and the same is true for the space $\Di_I(\R)$ from (\ref{DiI}). In the special case that $\Xc=\R$ equipped with the linear betweenness, the metrics $d_{\rm tot}$ and $d_{\rm part}$ generate Skorokhod's M1 topology on $\Di_I(\R)$ while $d_{\rm H}$ is a metric that generates the M2 topology.
\et

In general, we call the topology generated by $d_{\rm tot}$ or $d_{\rm part}$ from (\ref{dSkor}) the \emph{Skorokhod topology} corresponding to the betweenness on $\Xc$. Choosing the trivial betweenness then yields Skorokhod's J1 topology while the linear betweenness on $\R$ yields Skorokhod's M1 topology. Various ways to extend the M1 topology to more general spaces than $\Xc=\R$ have been proposed in the literature, not all of which are covered by Theorem~\ref{T:M1}. In particular, it is natural to identify $\Di_I(\R^n)$ with the product space $\Di_I(\R)^n$, equipped with the product topology, but the topology one obtains in this way does not come from a betweenness. In view of this, we will reserve the term ``M1 topology'' for the one-dimensional setting only. We note that the Skorokhod topology corresponding to the linear betweenness on a Banach space $\Xc$ is called the \emph{strong M1 topology} in \cite{PM15}.

The following lemma implies that the J1 topology is stronger than any other Skorokhod topology, defined by a betweenness other than the trivial betweenness. In particular, convergence in the J1 topology implies convergence in the M1 topology.

\bl[Comparison of Skorokhod topologies]\hspace{-7.05pt}
Let\label{L:J1vsM1} $\Xc$ be a metrisable space on which there are defined two proper betweennesses $\li\,\cdot\,,\,\cdot\,\re$ and $\li\,\cdot\,,\,\cdot\,\re'$ such that $\li x,z\re\subset\li x,z\re'$ for all $x,z\in\Xc$. Then the Skorokhod topology defined by the first betweenness is stronger than the Skorokhod topology defined by the second betweenness.
\el

\subsection{Compactness criteria}\label{S:compcrit}

In this subsection, we give a compactness criterion for subsets of the path space $\Pi(\Xc)$ equipped with the Skorokhod topology corresponding to a betweenness. We also give a compactness criterion for sets of continuous paths. Our criteria are similar to well-known criteria for spaces of functions defined on a fixed domain. Let $(\Xc,d)$ be a metric space. We say that a set $\Ai\sub\Pi(\Xc)$ satisfies the \emph{compact containment condition} if
\be\label{ccc}
\forall T<\infty\ \exists\mbox{ compact }C\sub\Xc\mbox{ s.t.\ }\pi(t\pm)\in C\ \forall\pi\in\Ai\mbox{ and }t\in I(\pi)\cap[-T,T].
\ee
For each $0<T<\infty$ and $\de>0$, we define the (traditional) \emph{modulus of continuity} of a continuous path $\pi\in\Pi_{\rm c}(\Xc)$ as
\be\label{modul}
m_{T,\de}(\pi):=\sup\big\{d\big(\pi(t_1),\pi(t_2)\big):t_1,t_2\in I(\pi),\ -T\leq t_1<t_2\leq T,\ t_2-t_1\leq\de\big\}.
\ee
We say that a set $\Ai\sub\Pi_{\rm c}(\Xc)$ is \emph{equicontinuous} if
\be
\lim_{\de\to 0}\sup_{\pi\in\Ai}m_{T,\de}(\pi)=0\qquad\forall T<\infty.
\ee
The following theorem generalises the classical Arzela--Ascoli theorem to sets of functions that are not necessarily all defined on the same domain, which moreover does not need to be an interval.

\bt[Arzela--Ascoli]
Let\label{T:ArzAsc} $\Xc$ be a metric space. Then a set $\Ai\sub\Pi_{\rm c}(\Xc)$ is precompact if and only if it is equicontinuous and satisfies the compact containment condition.
\et

For paths with jumps, it is possible to give a very similar compactness criterion. We first need to introduce a bit of notation. Let $\R_\sg$ be the space that consists of all words of the form $t\star$ where $t\in\R$ is a real number and $\star\in\{-,+\}$ is a sign. We think of $\R_\sg$ as obtained by cutting each point of the real line into two. Consequently, we call $\R_\sg$ the \emph{split real line} and call elements of $\R_\sg$ \emph{split real numbers}. We denote split real numbers either by words $t\star$ consisting of a Roman letter and a sign, or by a single Greek letter. In this case, if $\tau=t\star$, then we call $\un\tau:=t$ the \emph{real part} of $\tau$ and we call $\sg(\tau):=\star$ its \emph{sign}. We equip $\R_\sg$ with the lexicographic order, i.e., we set $\sig\leq\tau$ if and only if either $\un\sig<\un\tau$ or $\un\sig=\un\tau$ and $\sg(\sig)\leq\sg(\tau)$, where $\{-,+\}$ is equipped with the natural total order in which $-\leq +$. We write $\sig<\tau$ if $\sig\leq\tau$ and $\sig\neq\tau$. Although we do not need this here, we note that Kolmogorov \cite{Kol56} observed that it is possible to equip $\R_\sg$ with a topology such that continuous functions on $\R_\sg$ correspond to cadlag functions on $\R$, see Subsection~\ref{S:split} below.

Let $\Xc$ be a metrisable space that is equipped with a proper betweenness. For any path $\pi\in\Pi(\Xc)$, we define $I_\sg(\pi)\sub\R_\sg$ by
\be\label{Isgnw}
I_\sg(\pi):=\big\{t\star:t\in I(\pi),\ \star\in\{-,+\}\big\}.
\ee
Using the notation for split real numbers we have just introduced, for each $0<T<\infty$ and $\de>0$, we define the \emph{Skorokhod modulus of continuity} as
\be\ba{r@{\,}l}\label{mJM1}
\dis m^{\rm S}_{T,\de}(\pi):=\sup\big\{d\big(\pi(\tau_2),\li\pi(\tau_1),\pi(\tau_3)\re\big):&\dis\tau_1,\tau_2,\tau_3\in I_\sg(\pi),\ \tau_1\leq\tau_2\leq\tau_3,\\[5pt]
&\dis-T\leq\un\tau_1,\ \un\tau_3\leq T,\ \un\tau_3-\un\tau_1\leq\de\big\},
\ec
where as before $d(x,A)$ denotes the distance of a point $x$ to a set $A$ and $\li x,y\re$ is the segment with endpoints $x$ and $y$. We say that a set $\Ai\sub\Pi(\Xc)$ is \emph{Skorokhod-equicontinuous} if
\be
\lim_{\de\to 0}\sup_{\pi\in\Ai}m^{\rm S}_{T,\de}(\pi)=0\qquad\forall T<\infty.
\ee
Specialising these definitions to the trivial betweenness, for which $\li\pi(\tau_1),\pi(\tau_3)\re=\{\pi(\tau_1),\pi(\tau_3)\}$, yields the definitions of the \emph{J1-modulus of continuity} and \emph{J1-equicontinuity}. If $\Xc=\R$, equipped with the linear betweenness, then we speak of the \emph{M1-modulus of continuity} and \emph{M1-equicontinuity}. 

\bt[Compactness criterion]
Let\label{T:ArzAsc2} $(\Xc,d)$ be a metric space that is equipped with a proper betweenness. Then a set $\Ai\sub\Pi(\Xc)$ is precompact in the Skorokhod topology corresponding to this betweenness if and only if it is Skorokhod-equicontinuous and satisfies the compact containment condition.
\et

For functions $f\in\Di_I(\Xc)$ as defined in (\ref{DiI2}), one can check that the Skorokhod modulus of continuity can more simply be defined as
\be\ba{r@{\,}l}\label{fmJM1}
\dis m^{\rm S}_{T,\de}(f)=\sup\big\{d\big(f(t_2),\li f(t_1),f(t_3)\re\big):&\dis t_1,t_2,t_3\in I,\ t_1\leq t_2\leq t_3,\\[5pt]
&\dis-T\leq t_1,\ t_3\leq T,\ t_3-t_1\leq\de\big\}.
\ec
As a fairly simple consequence of Theorem~\ref{T:ArzAsc2}, we can derive the following compactnes criterion for $\Di_I(\Xc)$. As before $\pa I$ denotes the set of finite endpoints of a closed real interval $I$.

\bt[Functions on a fixed interval]
Let\label{T:ArzI} $(\Xc,d)$ be a metric space that is equipped with a proper betweenness and let $I$ be a closed real interval of positive length. Then a set $\Fi\sub\Di_I(\Xc)$ is precompact in the Skorokhod topology if and only if:
\begin{enumerate}
\item the compact containment condition holds,
\item $\dis\lim_{\de\to 0}\sup_{f\in\Fi}m^{\rm S}_{T,\de}(f)=0$ for all $T<\infty$,
\item $\dis\lim_{\de\to 0}\sup_{f\in\Fi}\;\sup\big\{d\big(f(s),f(t)\big):s\in I,\ |s-t|\leq\de\big\}=0$ for all $t\in\pa I$.
\end{enumerate}
\et

Note that compared to Theorem~\ref{T:ArzAsc2}, we need the extra condition~(iii) to rule out the scenario that a sequence of functions in $\Fi$ converges to a function with a discontinuity at a time $t\in\pa I$. For the J1 topology on $\Di_{[0,1]}$, Theorem~\ref{T:ArzI} was first proved by Kolmogorov in \cite[Thm~IV]{Kol56}. The analogous statement for the M1 topology was proved by Skorokhod in \cite[2.7.3]{Sko56}.

\subsection{The ordered Hausdorff topology}\label{S:orHa}

The metrics $d_{\rm tot}$ and $d_{\rm part}$ defined first in (\ref{dpath}) and (\ref{dpath2}) and then more generally in (\ref{dSkor}) are special cases of two metrics that can be defined on an even more general space of totally ordered compact sets. In the present subsection, we investigate the topology generated by $d_{\rm tot}$ and $d_{\rm part}$ on this general space. The results in this subsection are more theoretical in nature and can be skipped by readers who are interested in cadlag functions only.

Throughout this subsection $\Xc$ is a metrisable space. We let $\Ki_{\rm part}(\Xc)$ denote the space of all pairs $(K,\pre)$ where $K$ is a nonempty compact subset of $\Xc$ and $\pre$ is a partial order on $K$ that is \emph{compatible with the topology} in the sense that the set
\be\label{Kli}
K^{\li 2\re}:=\big\{(x,y)\in K^2:x\pre y\big\}
\ee
is a closed subset of $K^2$, where we equip $K$ with the induced topology from $\Ki$ and $K^2$ with the product topology. We let $\Ki_{\rm tot}(\Xc)$ denote the space of pairs $(K,\pre)\in\Ki_{\rm part}(\Xc)$ for which $\pre$ is a total order on $K$. We will sometimes be sloppy and denote elements of $\Ki_{\rm part}(\Xc)$ simply as $K$. We equip the space $\Xc^2$ with the metric
\be\label{d2}
d^2\big((x_1,y_1),(x_2,y_2)\big):=d(x_1,x_2)\vee d(y_1,y_2),
\ee
which generates the product topology and let $d^2_{\rm H}$ denote the associated Hausdorff metric on $\Ki_+(\Xc^2)$. An element $(K,\pre)$ of $\Ki_{\rm part}(\Xc)$ is clearly uniquely determined by the compact set $K^{\li 2\re}\sub\Xc^2$ defined in (\ref{Kli}), so setting
\be\label{part}
d_{\rm part}(K_1,K_2):=d^2_{\rm H}(K^{\li 2\re}_1,K^{\li 2\re}_2)\qquad\big(K_1,K_2\in\Ki_{\rm part}(\Xc)\big)
\ee
defines a metric $d_{\rm part}$ on $\Ki_{\rm part}(\Xc)$. For $K_1,K_2\in\Ki_{\rm tot}(\Xc)$, we let ${\rm Cor}_+(K_1,K_2)$ denote the set of monotone correspondences between $K_1$ and $K_2$, defined as in (\ref{moncor}), and define a metric $d_{\rm tot}$ on $\Ki_{\rm tot}(\Xc)$ by
\be\label{tot}
d_{\rm tot}(K_1,K_2):=\inf_{R\in{\rm Cor}_+(K_1,K_2)}\sup_{(x_1,x_2)\in R}d(x_1,x_2)
\qquad\big(K_1,K_2\in\Ki_{\rm tot}(\Xc)\big).
\ee
The following theorem is a major ingredient in the proof of Theorems \ref{T:J1} and \ref{T:M1}. Formula (\ref{partot}) generalises (\ref{dcomp}) and (\ref{dcomp2}). The final statement shows that in general no inequality of the form $d_{\rm tot}\leq Cd_{\rm part}$ holds for any $C<\infty$.


\bt[The ordered Hausdorff topology]
Let\label{T:partot} $\Xc$ be a metrisable space, let $d$ be a metric on $\Xc$ and let $d_{\rm part}$ and $d_{\rm tot}$ be defined in terms of $d$ as in (\ref{part}) and (\ref{tot}). Then $d_{\rm part}$ and $d_{\rm tot}$ are metrics that generate the same topology on $\Ki_{\rm tot}(\Xc)$, and this topology does not depend on the choice of the metric $d$ on $\Xc$. If $\Xc$ is Polish, then so is $\Ki_{\rm tot}(\Xc)$. One has
\be\label{partot}
d_{\rm H}(K_1,K_2)\leq d_{\rm part}(K_1,K_2)\leq d_{\rm tot}(K_1,K_2)
\qquad\big(K_1,K_2\in\Ki_{\rm tot}(\Xc)\big).
\ee
If $\Xc=[0,1]$, then for each $\eps>0$, there exist $K_1,K_2\in\Ki_{\rm tot}(\Xc)$ such that $d_{\rm part}(K_1,K_2)\leq\eps d_{\rm tot}(K_1,K_2)$.
\et

We call the topology on $\Ki_{\rm tot}(\Xc)$ generated by $d_{\rm part}$ or $d_{\rm tot}$ the \emph{ordered Hausdorff topology.} Since $\Ki_{\rm tot}(\Xc)$ is Polish, there exists a complete metric generating the topology, but as we will show in Lemma~\ref{L:noncompl} below, if $(\Xc,d)$ is complete, then $\Ki_{\rm tot}(\Xc)$ is not in general complete in the metrics $d_{\rm part}$ or $d_{\rm tot}$.

We also have a compactness criterion for subsets of $\Ki_{\rm tot}(\Xc)$. For $K\in\Ki_{\rm tot}(\Xc)$ and $\eps>0$, we define the \emph{mismatch modulus} $m_\eps(K)$ as
\be\ba{r@{\,}l}\label{mismatch}
\dis m_\eps(K):=\sup\big\{d(x_1,y_1)\vee d(x_2,y_2):
&\dis x_1,y_1,x_2,y_2\in K\\[5pt]
&\dis d(x_1,x_2)\vee d(y_1,y_2)\leq\eps,\ x_1\pre y_1,\ y_2\pre x_2\big\}.
\ec
The following theorem is a major ingredient in the proof of Theorems \ref{T:ArzAsc}, \ref{T:ArzAsc2}, and \ref{T:ArzI}.

\bt[Compact subsets]
Let\label{T:totcomp} $(\Xc,d)$ be a metric space and let $\Ki_{\rm tot}(\Xc)$ be equipped with the ordered Hausdorff topology. Then a set $\Ai\sub\Ki_{\rm tot}(\Xc)$ is precompact if and only if
\be\label{totcomp}
{\rm(i)}\ \exists\mbox{ compact }C\sub\Xc\mbox{ s.t.\ }K\sub C\ \forall K\in\Ai\quand{\rm(ii)}\ \lim_{\eps\to 0}\sup_{K\in\Ai}m_\eps(K)=0.
\ee
\et

We need one more result that will help us make the link between our definition of the Skorokhod topologies and the more traditional formulations. Recall the definition of the space $\Di_{[0,1]}(\Xc)$ of cadlag functions $f\cn[0,1]\to\Xc$ in (\ref{DiI2}). A \emph{cadlag parametrisation} of an element $K\in\Ki_{\rm tot}(\Xc)$ is a function $\ga\in\Di_{[0,1]}(\Xc)$ such that
\be
K=\big\{\ga(t),\ga(t-):t\in[0,1]\big\}
\quand
\ga(s)\prec\ga(t)\ \forall 0\leq s<t\leq 1.
\ee
Clearly, not every element of $\Ki_{\rm tot}(\Xc)$ has a cadlag parametrisation. For those that do, the following proposition gives an expression for the metrics $d_{\rm H}$ and $d_{\rm tot}$ that will later help us make the connection between our definitions and the classical definitions of the J1 and M1 topologies. Let $\La$ be the space of all bijections $\la\cn[0,1]\to[0,1]$ and let $\La_+$ be the subset consisting of all bijections $\la$ that are monotone in the sense that $s\leq t$ implies $\la(s)\leq\la(t)$. Note that each $\la\in\La_+$ is continuous and strictly increasing with $\la(0)=0$ and $\la(1)=1$.

\bp[Distance between cadlag curves]
Let\label{P:curve} $(\Xc,d)$ be a metric space, and assume that $K_1,K_2\in\Ki_{\rm tot}(\Xc)$ have cadlag parametrisations $\ga_1,\ga_2$, respectively. Then
\bc\label{curcor}
\dis d_{\rm H}(K_1,K_2)&=&\dis\inf_{\la\in\La}\sup_{t\in[0,1]}d\big(\ga_1(t),\ga_2\big(\la(t)\big)\big),\\[5pt]
\dis d_{\rm tot}(K_1,K_2)&=&\dis\inf_{\la\in\La_+}\sup_{t\in[0,1]}d\big(\ga_1(t),\ga_2\big(\la(t)\big)\big).
\ec
\ep

\subsection{Discussion}

It is instructive to compare our approach with the traditional approach to the Skorokhod topologies. For each $0<t<\infty$, let $\La_+[0,t]$ denote the space of continuous strictly increasing functions $\la\cn[0,t]\to[0,t]$ such that $\la(0)=0$ and $\la(t)=t$. Then for any metric space $(\Xc,d)$, setting
\be\label{timchan}
d_{\rm S}(f_1,f_2):=\inf_{\la\in\La_+[0,t]}\sup_{s\in[0,t]}\big[d\big(f_1(s),f_2(\la(s))\big)\vee|s-\la(s)|\big]\qquad\big(f_1,f_2\in\Di_{[0,t]}(\Xc)\big)
\ee
defines a metric on $\Di_{[0,t]}(\Xc)$ that generates the J1 topology. The definition of the M1 topology is similar, but instead of time changes $\la\in\La_+[0,t]$ one now has to work with parametrisations of the filled graphs of $f_1$ and $f_2$. For brevity, we omit the details. The metric in (\ref{timchan}) works only for compact time intervals. To define the J1 topology on $\Di_\half(\Xc)$, one needs the following result.

\bl[J1 topology on the half-line]
Let\label{L:halfline} $(\Xc,d)$ be a metric space and for each $0<t<\infty$ let $d_{\rm S}$ denote the metric on $\Di_{[0,t]}(\Xc)$ defined in (\ref{timchan}). Then there exists a metric $d_{\rm S}$ on $\Di_\half(\Xc)$ such that for all $f_n,f\in\Di_\half(\Xc)$, one has $d_{\rm S}(f_n,f)\to 0$ if and only if
\be
d_{\rm S}\big(f_n\big|_{[0,t]},f\big|_{[0,t]}\big)\asto{n}0\quad\forall t>0\mbox{ s.t.\ }f(t-)=f(t).
\ee
\el

There is a close connection between the metric $d_{\rm S}$ defined in (\ref{timchan}) and our metric $d_{\rm tot}$. To see this, equip the space $\Xc\times[0,t]$ with the metric $d'\big((x,s),(y,r)\big):=d(x,y)\vee|s-r|$. Fix $f_1,f_2\in\Di_{[0,t]}(\Xc)$ and let $\pi_1,\pi_2$ denote their associated paths, viewed as totally ordered compact subsets of $\Xc\times[0,t]$. Then Proposition~\ref{P:curve} says that $d_{\rm S}(f_1,f_2)=d_{\rm tot}(\pi_1,\pi_2)$, where $d_{\rm tot}$ is the metric on $\Ki_{\rm tot}(\Xc\times[0,t])$ defined in terms of $d'$ as in (\ref{tot}). By contrast, our metric $d_{\rm part}$ does not seem to correspond to a previously known metric.

A difference between our metric $d_{\rm tot}$ and the metric $d_{\rm S}$ defined in (\ref{timchan}) is that the latter only makes sense for paths that are defined on a compact interval and that do not jump at the endpoints of this interval. To deal with unbounded intervals, in the traditional approach, one still needs to construct a metric as in Lemma~\ref{L:halfline}. By contrast, our method to deal with unbounded intervals is based on the squeezed space, which cares less about spatial distances between points whose time coordinates are very large or very negative. Occasionally, one needs functions that can jump at finite endpoints of the interval on which they are defined. In the traditional approach, this is usually solved in an ad hoc manner, by embedding the function space under consideration in a slightly larger one, see, for example, the definitions leading up to \cite[Thm~7]{Jak96}. As our approach shows, allowing jumps at the endpoints of intervals is actually the natural thing to do as it leads to the clean compactness criterion of Theorem~\ref{T:ArzAsc2} compared to the more complicated Theorem~\ref{T:ArzI}. It should be remarked, however, that Aldous' tightness criterion \cite{Ald78} in a stochastic context nicely combines conditions (ii) and (iii) of Theorem~\ref{T:ArzI}.

An attractive aspect of our approach is that in many cases where discrete-time processes approximate a continuous-time process, it becomes unnecessary to interpolate. Let $\Xc$ be equipped with a betweenness and let $\pi\in\Pi_{\rm c}(\Xc)$. Then we say that $\pi'\in\Pi_{\rm c}(\Xc)$ is an \emph{interpolation} of $\pi$ if $I(\pi')$ is the convex hull of $I(\pi)$ and for each $s,u\in I(\pi)$ such that $s<u$ and $(s,u)\cap I(\pi)=\emptyset$, the function $[s,u]\ni t\mapsto\pi'(t)$ takes values in the segment $\li\pi(s),\pi(u)\re$ and is monotone with respect to the total order on $\li\pi(s),\pi(u)\re$. Note that an interpolation with respect to the trivial betweenness is a piecewise constant interpolation. We will prove the following simple lemmas.

\bl[Continuous interpolation]
Let\label{L:interpol1} $\Xc$ be a metrisable space that is equipped with a proper betweenness. Let $\pi_n\in\Pi_{\rm c}(\Xc)$ and $\pi\in\Pi^|_{\rm c}(\Xc)$ and for each $n$, let $\pi'_n$ be an interpolation of $\pi_n$. Assume that for each $t\in I(\pi)$, there exist $t_n\in I(\pi_n)$ such that $t_n\to t$. Then $\pi_n\to\pi$ in the topology on $\Pi_{\rm c}(\Xc)$ if and only if $\pi'_n\to\pi$.
\el

\bl[Piecewise constant interpolation]
Let\label{L:interpol2} $\Xc$ be a metrisable space that is equipped with the trivial betweenness. Let $\pi_n\in\Pi_{\rm c}(\Xc)$ and $\pi\in\Pi^|(\Xc)$ and for each $n$, let $\pi'_n$ be an interpolation of $\pi_n$. Assume that for each $t\in I(\pi)$, there exist $t_n\in I(\pi_n)$ such that $t_n\to t$. Then $\pi_n\to\pi$ in the topology on $\Pi(\Xc)$ if and only if $\pi'_n\to\pi$.
\el

Lemmas \ref{L:interpol1} and \ref{L:interpol2} show that for convergence to a continuous function or convergence in the J1 topology, it is not necessary to interpolate. For convergence in the M1 topology, it may sometimes still make sense to interpolate, however.

By definition, a \emph{laglad} function (from the French ``limite \`a gauche, limite \`a droit'') is a function that has both left- and right- limits in each point, but whose value in a point does not need to be equal to either the left or right limit at that point. Whitt \cite[Chapter~15]{Whi02} introduces a topology on spaces of laglad functions. It seems it should be possible to generalise many of our results to laglad functions, where as a first step, one needs a variant of the split real line in which each each point of the real line is replaced by three, rather than two points. The details are not straightforward, however, so we leave this as an open problem.

\subsection{Outline of the proofs}

In Section~\ref{S:cadfun} we collect some preliminary facts about cadlag functions, the Hausdorff topology, and betweenness. In Subsections \ref{S:split} and \ref{S:splitproof} we develop the observation, originally due to Kolmogorov \cite{Kol56}, that there exists a topology on the split real line $\R_\sg$ introduced below Theorem~\ref{T:ArzAsc} such that cadlag functions on $\R$ corespond to continuous functions on $\R_\sg$. Although we will not use this very often, it will occasionally be handy in our proofs. The facts in this subsection are of some independent interest as they can also be used to give a natural definition of cadlag functions of several variables, that seems simpler than the approach used by other authors such as \cite{BW71,Neu71}. Subsection~\ref{S:squeeze} is devoted to the squeezed space $\Ri(\Xc)$ and contains the proof of Lemmas \ref{L:squeeze} and \ref{L:interpol2}. Subsection~\ref{S:Haus} collects some facts about the Hausdorff topology (such as Lemma~\ref{L:Haus}) that will be needed many times later on. In Subsection~\ref{S:between} we study betweenness and prove Lemmas \ref{L:partdef}, \ref{L:metbet}, and \ref{L:compbet}. Subsection~\ref{S:cad}, finally, contains the proofs of Lemmas \ref{L:fpi} and \ref{L:cogra}.

Section~\ref{S:Hord} is devoted to the ordered Hausdorff topology. Theorem~\ref{T:partot} is proved in Subsections \ref{S:mis} and \ref{S:Pol}, Theorem~\ref{T:totcomp} in Subsection~\ref{S:totcomp}, and Proposition~\ref{P:curve} Subsection~\ref{S:curve}. These results lay the basis for the proofs of our main results, which can be found in Section~\ref{S:mainproof}. The basic statements about the metrics $d_{\rm part}$, $d_{\rm tot}$, and $d_{\rm H}$ from Theorems \ref{T:J1} and \ref{T:M1} are proved in Subsection~\ref{S:Pimet}, while the statements about Polishness are proved in Subsection~\ref{S:JPol}. The compactness criteria Theorems \ref{T:ArzAsc}, \ref{T:ArzAsc2}, and \ref{T:ArzI} are proved in Subsection~\ref{S:cocrit}, which also contains the proof of Lemmas \ref{L:J1vsM1} and \ref{L:interpol1}. In Subsection~\ref{S:classic}, we complete the proofs of Theorems \ref{T:J1} and \ref{T:M1} by showing that our definitions of the J1, J2, M1, and M2 topologies on the space $\Di_I(\Xc)$ coincide with the classical definitions. This subsection also contains an alternative proof of Lemma~\ref{L:halfline}, which is of course well known from the literature.

\section{Paths and cadlag functions}\label{S:cadfun}

\subsection{The split real line}\label{S:split}

Basic notation for the split real line $\R_\sg$ has been introduced in Subsection~\ref{S:compcrit} below Theorem~\ref{T:ArzAsc}. This space can be equipped with a natural topology that is the subject of this subsection. The split unit interval $[0,1]_\sg$ (defined as in (\ref{Isgnw})) as a topological space originates from the work of Alexandroff and Urysohn \cite{AU29}. It is known as the \emph{split interval}, \emph{double arrow space}, or \emph{two arrows space}. An account of many of its properties can be found in \cite[Chapter~41]{Fre03}. The topology on the split unit interval $[0,1]_\sg$ with the points $0-$ and $1+$ removed is called the \emph{weak parallel line topology} in \cite{SS78}, where again several of its properties are proved. These references do not discuss cadlag functions, however. Below, we give an overview of the most important properties of the topology on the split real line. Proofs will be provided in the next subsection.

We use notation for intervals in $\R_\sg$ similar to the usual notation for the real line, i.e.,
\be\ba{ll}\label{interval}
\dis\lo\sig,\rho\ro:=\{\tau\in\R_\sg:\sig<\tau<\rho\},\quad
&\dis\lc\sig,\rho\ro:=\{\tau\in\R_\sg:\sig\leq\tau<\rho\},\\[5pt]
\dis\lo\sig,\rho\rc:=\{\tau\in\R_\sg:\sig<\tau\leq\rho\},\quad
&\dis\lc\sig,\rho\rc:=\{\tau\in\R_\sg:\sig\leq\tau\leq\rho\}.
\ec
Note that there is some redundancy in this notation: for example, $\lo s-,t+\ro=\lc s+,t-\rc$. We equip $\R_\sg$ with the \emph{order topology}, which is the weakest topology to make all sets of the form $\{\tau\in\R_\sg:\tau>\sig\}$ and $\{\tau\in\R_\sg:\tau<\sig\}$ with $\sig\in\R_\sg$ open. A consequence of this is that intervals of the form $\lo\sig,\rho\ro$ are open, those of the form $\lc\sig,\rho\rc$ are closed, and $\lo s-,t+\ro=\lc s+,t-\rc$ is both open and closed. The following lemma describes convergent sequences in this topology.

\bl[Convergence criterion]
For\label{L:concr} $\tau_n\in\R_\sg$ and $t\in\R$, one has
\begin{enumerate}
\item $\tau_n\to t+$ if and only if $\un\tau_n\to t$ and $\tau_n\geq t+$ for
  all $n$ sufficiently large,
\item $\tau_n\to t-$ if and only if $\un\tau_n\to t$ and $\tau_n\leq t-$ for
  all $n$ sufficiently large.
\end{enumerate}
\el

The following lemma lists some elementary properties of $\R_\sg$.

\bl[The split real line]
The\label{L:Rpm} space $\R_\sg$ is first countable, Hausdorff, and separable,
but not second countable and not metrisable. Moreover, $\R_\sg$ is totally
disconnected, meaning that its only connected subsets are singletons.
\el

Let $\Xc$ be a Hausdorff topological space. For any set $\Ii\sub\R_\sg$, we let $\Ci_\Ii(\Xc)$ the space of continuous functions $f\cn\Ii\to\Xc$. For any $I\sub\R$, we define $I_\sg\sub\R_\sg$ by
\be\label{Isg}
I_\sg:=\big\{t\star:t\in I,\ \star\in\{-,+\}\big\}.
\ee
Recall that in Subsection~\ref{S:paths} we defined a cadlag function with values in $\Xc$ to be a triple $(I,f^-,f^+)$ where $I\sub\R$ is closed and $f^-,f^+\cn I\to\Xc$ are left and right continuous functions that satisfy (\ref{cadlag}). The following lemma is a direct consequence of Lemma~\ref{L:concr}.

\bl[Cadlag functions]
Let\label{L:cadlag} $\Xc$ be a Hausdorff topological space, let $I\sub\R$ be closed, let $f\cn I_\sg\to\Xc$ be a function, and let $f^\pm\cn I\to\Xc$ be defined by $f^\pm(t):=f(t\pm)$ $(t\in I)$. Then $f\in\Ci_{I_\sg}(\Xc)$ if and only if $(I,f^-,f^+)$ is a cadlag function.
\el

We equip the product space $\R_\sg^d$ with the product topology. By definition, a subset $A\sub\R_\sg^d$ is \emph{bounded} if $A\sub\lc\sig,\tau\rc^d$ for some $\sig,\tau\in\R_\sg$. The following proposition gives a characterisation of the compact subsets of $\R_\sg^d$, similar to the well-known characterisation of compact subsets of $\R^d$.

\bp[Compact sets]
For\label{P:Rpmcompact} a subset $C\sub\R_\sg^d$, the following three claims are
equivalent: (i) $C$ is compact, (ii) $C$ is sequentially compact, and (iii)
$C$ is closed and bounded.
\ep

We define the \emph{extended split real line} by $\ov\R_\sg:=\{t\star:t\in\ov\R,\ \star\in\{-,+\}\}$ and equip it with a topology so that it is homeomorphic to $\lc 0-,1+\rc$. We also set
\be\label{hatRsg}
\hat\R_\sg:=\lc-\binf,+\binf\rc\quad\mbox{with}\quad-\binf:=-\infty+\mbox{ and }+\binf:=+\infty-.
\ee
Note that $\hat\R_\sg$ is homeomorphic to $\lc 0+,1-\rc$, therefore compact by Lemma~\ref{P:Rpmcompact}, and that $\R_\sg$ is dense in $\hat\R_\sg$. In other words, $\hat\R_\sg$ is a compactification of $\R_\sg$. The notation introduced in (\ref{hatRsg}) provides us with a natural way to denote half infinite intervals in $\R_\sg$; for example, $\lc\sig,\binf\ro=\{\tau\in\R_\sg:\sig\leq\tau\}$.


\subsection{Proofs of the properties of the split real line}\label{S:splitproof}

In this subsection we prove the properties of the split real line stated in the previous subsection.\med

\bpro[of Lemma~\ref{L:concr}]
By symmetry, it suffices to prove (i). A basis for the topology is formed by all intervals of the form $\lo\sig,\rho\ro$ with $\sig,\rho\in\R_\sg$. If $t+\in\lo\sig,\rho\ro$, then $\sig<t+<\rho$ and hence $\lo t-,u+\ro\sub\lo\sig,\rho\ro$ for some $u>t$. It follows that the sets of the form $\lo t-,u+\ro=\lc t+,u-\rc$ with $u\in\{t+n^{-1}:n\geq 1\}$ form a fundamental system of neighbourhoods of $t+$, which is easily seen to imply the claim.
\epro

\bpro[of Lemma~\ref{L:Rpm}]
It is easy to see that $\R_\sg$ has the Hausdorff property. In the proof of Lemma~\ref{L:concr}, we have already seen that each point has a countable fundamental system of neighbourhoods, so $\R_\sg$ is first countable. On the other hand, each basis of the topology must for each $t\in\R$ contain an open set $O$ such that $t\in O\sub\lo t-,(t+1)+\ro=\lc t+,(t+1)-\rc$. These open sets are all distinct, so $\R_\sg$ is not second countable. By Lemma~\ref{L:concr}, the set $\{t+:t\in\Q\}$ is dense so $\R_\sg$ is separable. Since in metric spaces, separability implies second countability, we conclude that $\R_\sg$ is not metrisable. Since for each $t\in\R$, we can write $\R_\sg$ as the union of two disjoint closed sets as $\R_\sg=\lo-\binf,t-\rc\cup\lc t+,\binf\ro$, we see that $\R_\sg$ is totally disconnected.
\epro

\bpro[of Lemma~\ref{L:cadlag}]
If $f\in\Ci_{I_\sg}(\Xc)$, then Lemma~\ref{L:concr} implies that $f^-$ and $f^+$ are left and right continuous, respectively, and satisfy (\ref{cadlag}). Conversely, if $f^-$ and $f^+$ are left and right continuous and satisfy (\ref{cadlag}), then Lemma~\ref{L:concr} implies that $f(\tau_n)\to f(\tau)$ for each $\tau_n,\tau\in I_\sg$ such that $\tau_n\to\tau$ in the topology on $\R_\sg$. Since $\R_\sg$ is first countable by Lemma~\ref{L:Rpm}, this implies that $f\cn I_\sg\to\Xc$ is continuous.
\epro

The next lemma prepares for the proof of Proposition~\ref{P:Rpmcompact}. Even though $\R_\sg$ is not second countable, it has a property that is almost as good. We cite the following property from \cite[Example~96, point~5]{SS78}.

\bl[Strong Lindel\"of property]
Every open cover of a subset of $\R_\sg$ has a countable subcover.
\el

We note that the product space $\R_\sg\times\R_\sg$, equipped with the product topology, does \emph{not} have the strong Lindel\"of property. Indeed, the collection of open sets
\be\label{uncover}
\big\{\lc t+,\binf\ro\times\lc-t+,\binf\ro:t\in\R\big\}
\ee
covers the set $\{\lo t+,-t+\ro:t\in\R\}$, but no countable subset of (\ref{uncover}) has this property. The set $\{(s+,t+):s,t\in\R\}$ with the induced topology from $\R^2_\sg$ is a well-known counterexample in topology, known as the \emph{Sorgenfrey plane}.\med

\bpro[of Proposition~\ref{P:Rpmcompact}]
We first prove the statement for $\R_\sg$. In any first countable space, being sequentially compact is equivalent to being countably compact, which means that every countable open covering has a finite subcovering. By the strong Lindel\"of property, a subset of $\R_\sg$ is compact if and only if it is countably compact, proving that (i) and (ii) are equivalent.

Since the map $\tau\mapsto\un\tau$, that assigns to a split real number $\tau$ its real part, is continuous, and since the continuous image of a compact set is compact, (i) implies that $\un C:=\{\un\tau:\tau\in C\}$ is closed and bounded. Since moreover a compact subset of a Hausdorff space is closed, (i) implies (iii).

Property (iii) implies that each sequence $\tau_n\in C$ has a subsequence $\tau'_n$ such that $\un\tau'_n$ converges to a limit $t\in\un C$. The $\tau'_n$ must then contain a further subsequence $\tau''_n$ such that one of the following three cases occurs: 1. $\un\tau''_n<t$ for all $n$, 2.\ $\un\tau''_n>t$ for all $n$, or 3.\ $\un\tau''_n$ is constant. In either case, the fact that $C$ is closed implies that $\tau''_n$ converges to a limit in $C$, proving the implication (iii)$\volgt$(ii). This completes the proof for $\R_\sg$.

We saw before that the strong Lindel\"of property does not hold for $\R_\sg^d$ in dimensions $d\geq 2$, so to prove the statement for these spaces we have to proceed differently. Property~(i) implies countable compactness which by the fact that $\R_\sg$ and hence also $\R_\sg^d$ are first countable is equivalent to (ii). The continuous image of a countably compact set is countably compact. Applying this to the coordinate projections and using what we already know for $\R_\sg$, we see that (ii) implies that $C$ is bounded. Since, moreover, in any first countable Hausdorff space, being sequentially compact implies being closed, we see that (ii) implies (iii). By Tychonoff's theorem and what we already know for $\R_\sg$, the set $\lc s-,t+\rc^d$ is compact for each $-\infty<s<t<\infty$. Since a closed subset of a compact set is compact, (iii) implies (i).
\epro

\subsection{Squeezed space}\label{S:squeeze}

In this subsection we prove Lemmas \ref{L:squeeze} and \ref{L:interpol2}. A metric generating the topology on $\Ri(\ov\R)$ is given in \cite[formula (3.4)]{FINR04} but this formula does not easily extend to more general metric spaces $(\Xc,d)$, so we will proceed a bit differently. Let $d_{\ov\R}$ be a metric that generates the topology on the extended real line $\ov\R$ and let $\phi\cn\ov\R\to\half$ be a continuous function such that $\phi(\pm\infty)=0$ and $\phi(t)>0$ for all $t\in\R$. We define $d_{\rm sqz}\cn\Ri(\Xc)^2\to\half$ by
\be\label{sqz}
d_{\rm sqz}\big((x,s),(y,t)\big):=\big(\phi(s)\wedge\phi(t)\big)\big(d(x,y)\wedge 1\big)+\big|\phi(s)-\phi(t)\big|+d_{\ov\R}(s,t),
\ee
where naturally the first term is zero if $(x,s)$ or $(y,t)$ are elements of $\big\{(\ast,-\infty),(\ast,\infty)\big\}$ (even though $d(x,y)$ is not defined in this case).

\bl[Squeezed space]
Let \label{L:sqz} $(\Xc,d)$ be a metric space. Then $d_{\rm sqz}$ is a metric on $\Ri(\Xc)$. One has $d_{\rm sqz}\big((x_n,t_n),(x,t)\big)\to 0$ if and only if:
\[
{\rm(i)}\quad t_n\to t,\qquad{\rm(ii)}\quad\mbox{if $t\in\R$, then }x_n\to x.
\]
\el

\bpro
We first prove that $d_{\rm sqz}$ is a metric on $\Ri(\Xc)$. For brevity, we write $d'(x,y):=d(x,y)\wedge 1$. Then $d'$ is a metric on $E$. The only nontrivial statement that we have to prove is the triangle inequality, and it suffices to prove this for the function
\[
\rho\big((x,s),(y,t)\big):=\big(\phi(s)\wedge\phi(t)\big)d'(x,y)+\big|\phi(s)-\phi(t)\big|.
\]
We estimate
\be\label{rhop}
\rho\big((x,s),(z,u)\big)
\leq\big(\phi(s)\wedge\phi(u)\big)\big(d'(x,y)+d'(y,z)\big)
+\big|\phi(s)-\phi(u)\big|.
\ee
If $\phi(t)\geq\phi(s)\wedge\phi(u)$, then $\phi(s)\wedge\phi(u)$ is less than $\phi(s)\wedge\phi(t)$ and also less than $\phi(t)\wedge\phi(u)$, so we can simply estimate the expression in (\ref{rhop}) from above by
\[
\big(\phi(s)\wedge\phi(t)\big)d'(x,y)+\big(\phi(t)\wedge\phi(u)\big)d'(y,z)\big)+\big|\phi(s)-\phi(t)\big|+\big|\phi(t)-\phi(u)\big|
\]
and we are done. On the other hand, if $\phi(t)<\phi(s)\wedge\phi(u)$, then
\[
\big|\phi(s)-\phi(t)\big|+\big|\phi(t)-\phi(u)\big|
=\big|\phi(s)-\phi(u)\big|+2\big(\phi(s)\wedge\phi(u)-\phi(t)\big).
\]
Using the fact that $d'\leq 1$, we can now estimate the right-hand side of (\ref{rhop}) from above by
\[\ba{l}
\dis\phi(t)\big(d'(x,y)+d'(y,z)\big)+2\big(\phi(s)\wedge\phi(u)-\phi(t)\big)
+\big|\phi(s)-\phi(u)\big|\\[5pt]
\dis\quad=\big(\phi(s)\wedge\phi(t)\big)d'(x,y)
+\big(\phi(t)\wedge\phi(u)\big)d'(y,z)\\[5pt]
\dis\quad\phantom{=}+\big|\phi(s)-\phi(t)\big|+\big|\phi(t)-\phi(u)\big|,
\ea\]
and again we are done. This completes the proof that $d_{\rm sqz}$ is a metric on $\Ri(\Xc)$.

It remains to prove that
\be\label{sqz0}
\big(\phi(t_n)\wedge\phi(t)\big)\big(d(x_n,x)\wedge 1\big)+\big|\phi(t_n)-\phi(t)\big|+d_{\ov\R}(t_n,t)\asto{n}0
\ee
if and only if conditions (i) and (ii) of the lemma are satisfied. Because of the third term on the left-hand side, a necessary condition for (\ref{sqz}) is that $t_n\to t$, and this condition also guarantees that the second term tends to zero. If $t\in\{-\infty,+\infty\}$, then this is all one needs since the first term now tends to zero regardless of the values of $x_n$ and $x$, but if $t\in\R$, then one needs in addition that $d(x_n,x)\to 0$.
\epro

The following lemma shows that $\Ri(\Xc)$ is Polish if $\Xc$ is.

\bl[Preservation of Polishness]\ \\[5pt]
{\bf(a)}\label{L:pointprop} If $(\Xc,d)$ is separable, then so is $(\Ri(\Xc),d_{\rm sqz})$.\med

\noi
{\bf(b)} If $(\Xc,d)$ is complete, then so is $(\Ri(\Xc),d_{\rm sqz})$.
\el

\bpro
If $D$ is a countable dense subset of $(\Xc,d)$, then $D\times\Q$ is a countable dense subset of $(\Ri(\Xc),d_{\rm sqz})$, proving (a).

To prove (b), let $(x_n,t_n)$ be a Cauchy sequence in $(\Ri(\Xc),d_{\rm sqz})$. Then by (\ref{sqz}) $t_n$ is a Cauchy sequence in $\ov\R$ and hence $t_n\to t$ for some $t\in\ov\R$. If $t\in\R$, then by (\ref{sqz}) $x_n$ is a Cauchy sequence in $(\Xc,d)$ so by the completeness of the latter, $x_n\to x$ for some $x\in\Xc$. By Lemma~\ref{L:sqz}, it follows that $(x_n,t_n)$ converges, proving the completeness of $(\Ri(\Xc),d_{\rm sqz})$.
\epro

The following lemma identifies the compact subsets of $\Ri(\Xc)$. In particular, the lemma shows that $\Ri(\Xc)$ is compact if $\Xc$ is compact.

\bl[Compactness criterion]
A\label{L:REcomp} set $A\sub\Ri(\Xc)$ is precompact if and only if for each $T<\infty$, there exists a compact set $K\sub\Xc$ such that $\{x\in\Xc:\exists t\in[-T,T]\mbox{ s.t.\ }(x,t)\in A\}\sub K$.
\el

\bpro
Assume that $A\sub\Ri(\Xc)$ has the property that for each $T<\infty$, there exists a compact set $K\sub\Xc$ such that $\{x\in\Xc:(x,t)\in A,\ t\in[-T,T]\}\sub K$. To show that $A$ is precompact, we will show that each sequence $(x_n,t_n)\in A$ has a convergent subsequence. By the compactness of $\ov\R$, we can select a subsequence $(x'_n,t'_n)$ such that $t'_n\to t$ for some $t\in\ov\R$. If $t=\pm\infty$, then by Lemma~\ref{L:sqz} $(x'_n,t'_n)\to(\ast,\pm\infty)$ and we are done. Otherwise, there exists a $T<\infty$ such that $t'_n\in[-T,T]$ for all $n$ large enough. By assumption, there then exists a compact set $K\sub\Xc$ such that $x'_n\in K$ for all $n$ large enough, so we can select a further subsequence such that $(x''_n,t''_n)$ converges to a limit $(x,t)\in\Xc\times\R$.

Assume, on the other hand, that $A\sub\Ri(\Xc)$ has the property that for some $T<\infty$, there does not exist a compact set $K\sub\Xc$ such that $\{x\in\Xc:(x,t)\in A,\ t\in[-T,T]\}\sub K$. Set
\[
B:=\big\{x\in\Xc:(x,t)\in A\mbox{ for some }t\in[-T,T]\big\}
\]
The closure of $B$ cannot be compact, since this would contradict our assumption. It follows that there exists a sequence $x_n\in B$ that does not contain a convergent subsequence, and there exist $t_n\in[-T,T]$ such that $(x_n,t_n)\in A$. But then, in view of Lemma~\ref{L:sqz}, the sequence $(x_n,t_n)$ cannot contain a convergent subsequence either, proving that $A$ is not precompact.
\epro

\bpro[of Lemma~\ref{L:squeeze}]
Immediate from Lemmas \ref{L:sqz}, \ref{L:pointprop}, and \ref{L:REcomp}.
\epro

\bpro[of Lemma~\ref{L:interpol2}]
Let $d_{\rm sqz}$ be defined as in (\ref{sqz}) and let $d_{\rm tot}$ be defined as in (\ref{dpath}). To prove the lemma, it suffices to show that $d_{\rm tot}(\pi_n,\pi'_n)\to 0$. Let us say that $(s,u)$ is a \emph{gap} of $\pi_n$ if $s<u$, $s,u\in I(\pi_n)$, and $(s,u)\cap I(\pi_n)=\emptyset$. We call $(\pi_n(s),s)$ and $(\pi_n(u),u)$ the \emph{endpoints} of the gap. Let $C_n$ be the correspondence between $\pi_n$ and $\pi'_n$ consisting of all pairs $(z,z')$ with $z\in\pi_n$ and $z'\in\pi'_n$ such that either $z=z'$ or there exists a gap $(s,u)$ of $\pi_n$ such that $z=(x,t)$ is an endpoint of the gap while $z'=(x',t')$ satisfies $s<t'<u$ and $x'=x$. The definition of an interpolation implies that this correspondence is monotone, so using (\ref{sqz}) we obtain
\be
d_{\rm tot}(\pi_n,\pi'_n)\leq\sup_{(z,z')\in C}d_{\rm sqz}(z,z')
\leq\sup\big\{|\phi(s)-\phi(u)|+d_{\ov\R}(s,u):(s,u)\mbox{ is a gap of }\pi_n\big\}.
\ee
Since $\pi\in\Pi^|(\Xc)$ and for each $t\in I(\pi)$, there exist $t_n\in I(\pi_n)$ such that $t_n\to t$, the maximal length of gaps of $\pi_n$ as measured by $d_{\ov\R}$ tends to zero completing the proof.
\epro

\subsection{The Hausdorff metric}\label{S:Haus}

If $(\Xc,d)$ is a metric space, then the Hausdorff metric $d_{\rm H}$ on the space $\Ki_+(\Xc)$ of nonempty compact subsets of $\Xc$ has been defined in (\ref{Haus}). We have not actually given a proof or reference for the equality of the two formulas there. The following lemma fills this gap and shows that the infimum in the second formula is attained.

\bl[Hausdorff metric and correspondences]
Let\label{L:Haus2} $(\Xc,d)$ be a metric space. Then
\be\label{Haus2}
d_{\rm H}(K_1,K_2):=\sup_{z_1\in K_1}d(z_1,K_2)\vee\sup_{z_2\in K_2}d(z_2,K_1)
=\inf_{R\in{\rm Cor}(K_1,K_2)}\sup_{(x_1,x_2)\in R}d(x_1,x_2).
\ee
Moreover, there exists an $R\in{\rm Cor}(K_1,K_2)$ such that $d_{\rm H}(K_1,K_2)=\max_{(x_1,x_2)\in R}d(x_1,x_2)$.
\el

\bpro
Let $R\in{\rm Cor}(K_1,K_2)$ and let $D:=\sup_{(x_1,x_2)\in R}d(x_1,x_2)$. Then $d(x_1,K_2)\leq D$ and $d(x_2,K_1)\leq D$ for each $x_1\in K_1$, $x_2\in K_2$, and hence $d_{\rm H}(K_1,K_2)\leq D$. On the other hand, by the compactness of $K_2$ and the continuity of the function $d(x_1,\,\cdot\,)$, for each $x_1\in K_1$, there exists an $x_2\in K_2$ such that $d(x_1,K_2)=d(x_1,x_2)$. The same statement holds with the roles of $K_1$ and $K_2$ interchanged, so setting
\be
R:=\big\{(x_1,x_2)\in K_1\times K_2:d(x_1,x_2)\in\{d(x_1,K_2),d(x_2,K_1)\}\big\}
\ee
defines a correspondence between $K_1$ and $K_2$. By the compactness of $K_1$ and the continuity of the map $d(\,\cdot\,,K_2)$, there exists an $x'_1\in K_1$ such that $d(x'_1,K_2)=\max_{x_1\in K_1}d(x_1,K_2)$, and similarly there exists an $x''_2\in K_2$ such that $d(x''_2,K_1)=\max_{x_2\in K_2}d(x_2,K_1)$. By our earlier arguments, there exist $x'_2\in K_2$ and $x''_1\in K_1$ such that $d(x'_1,K_2)=d(x'_1,x'_2)$ and $d(x''_2,K_1)=d(x''_1,x''_2)$. Then
\be
d_{\rm H}(K_1,K_2)=d(x'_1,x'_2)\vee d(x''_1,x''_2)=\max_{(x_1,x_2)\in R}d(x_1,x_2).
\ee
\epro

We will often need the following lemma which we cite from \cite[Lemma~B.1]{SSS14}. Note that this lemma implies that the topology generated by $d_{\rm H}$ does not depend on the choice of the metric $d$ on $\Xc$.

\bl[Convergence criterion]
Let\label{L:Hauconv} $K_n,K\in\Ki_+(\Xc)$ $(n\geq 1)$. Then $K_n\to K$ in the Hausdorff topology if and only if there exists a $C\in\Ki_+(\Xc)$ such that $K_n\sub C$ for all $n\geq 1$ and
\bc\label{Haulim}
K&=&\dis\{x\in \Xc:\exists x_n\in K_n\mbox{ s.t.\ }x_n\to x\}\\[5pt]
&=&\dis\{x\in \Xc:\exists x_n\in K_n
\mbox{ s.t.\ $x$ is a cluster point of } (x_n)_{n\in\N}\}.
\ec
\el

The following lemma is \cite[Lemma~B.2]{SSS14}. In particular, it shows that $\Ki_+(\Xc)$ is Polish if $\Xc$ is.

\bl[Properties of the Hausdorff metric]\ \\[5pt]
{\bf(a)}\label{L:Hauprop} If $(\Xc,d)$ is separable, then so is $(\Ki_+(\Xc),d_{\rm H})$.\med

\noi
{\bf(b)} If $(\Xc,d)$ is complete, then so is $(\Ki_+(\Xc),d_{\rm H})$.
\el

Recall that a set is called precompact if its closure is compact. The following lemma is \cite[Lemma~B.3]{SSS14}. In particular, it shows that $\Ki_+(\Xc)$ is compact if $\Xc$ is.

\bl[Compactness in the Hausdorff topology]
A\label{L:Haucomp} set $\Ai\sub\Ki_+(\Xc)$ is precompact if and only if there exists a $C\in\Ki_+(\Xc)$ such that $K\sub C$ for each $K\in\Ai$.
\el

The following lemma says connectedness is a property of compact sets that is preserved under limits.

\bl[Preservation of connectedness]
The\label{L:Hconnec} set $\Ki_{\rm c}(\Xc)$ of all connected nonempty compact subsets of $\Xc$ is a closed subset of $\Ki_+(\Xc)$.
\el

\bpro
Imagine that $K_n,K\in\Ki_+(\Xc)$ satisfy $K_n\to K$. If $K$ is not connected, then there exist disjoint nonempty compact sets $C_1,C_2$ such that $K=C_1\cup C_2$. Let $\eps:=d(C_1,C_2)=\inf\{d(x_1,x_2):x_1\in C_1,\ x_2\in C_2\}$. By the compactness of $C_1$ and $C_2$, the infimum is attained and $\eps>0$. Let $U_i:=\{x\in\Xc:d(x,C_i)\leq\eps/3\}$ $(i=1,2)$. Then $U_1,U_2$ are disjoint closed sets. For all $n$ large enough such that $d_{\rm H}(K_n,K)\leq\eps/3$, one has $K_n\sub U_1\cup U_2$ while $K_n\cap U_1$ and $K_n\cap U_2$ are both nonempty, which proves that $K_n$ is not connected.
\epro

We recall that the image of a compact set under a continuous map is compact. In what follows, we will need the following simple observation.

\bl[Continuous image]
Let\label{L:contim} $\Xc,\Yi$ be metrisable spaces and let $\psi\cn\Xc\to\Yi$ be continuous. If $K_n,K\in\Ki_+(\Xc)$ satisfy $K_n\to K$, then their images under $\psi$ satisfy $\psi(K_n)\to\psi(K)$ in $\Ki_+(\Xc)$.
\el

\bpro
This follows easily from Lemma~\ref{L:Hauconv}. Since $K_n\to K$, there exists a compact set $C\sub\Xc$ such that $K_n\sub C$ for all $n$, and now $\psi(C)\sub\Yi$ is a compact set such that $\psi(K_n)\sub\psi(C)$ for all $n$. By (\ref{Haulim}), it now suffices to check that:
\be\ba{rl}\label{suff}
{\rm(i)}&\dis\psi(K)\sub\{y\in\Yi:\exists x_n\in K_n\mbox{ s.t.\ }\psi(x_n)\to y\},\\[5pt]
{\rm(ii)}&\dis\{y\in\Yi:\exists x_n\in K_n
\mbox{ s.t.\ $y$ is a subsequential limit of } (\psi(x_n))_{n\in\N}\}\sub\psi(K).
\ec
Here (i) follows from (\ref{Haulim}) and the continuity of $\psi$. To prove (ii), if $\psi(x'_n)\to y$ for some subsequence $x'_n$, then since $K_n\sub C$ for all $n$ there exists a further subsequence $x''_n$ such that $x''_n\to x$ for some $x\in\Xc$. Then $x\in K$ by (\ref{Haulim}) and hence $\psi(x)=y$ by the continuity of $\psi$.
\epro

Although we do not need this for our main results, it is useful to point out the relation between the Hausdorff topology on $\Ki_+(\Xc)$ and the much more widely known Vietoris topology on the space ${\rm Clos}(\Xc)$ of closed subsets of $\Xc$. By definition, the latter is generated by all sets of the form
\be
\big\{A:A\cap O\neq\emptyset\big\}
\quand
\big\{A:A\sub O\big\}
\ee
where $O\sub\Xc$ is open, see \cite[Appendix~4A2]{Fre03}. In the special case that $\Xc$ is compact the following lemma is proved in \cite[4A2Tg(ii)]{Fre03} but we have not found a reference for the general statement because most authors only seem interested in ${\rm Clos}(\Xc)$ and not in $\Ki_+(\Xc)$.

\bl[Hausdorff and Vietoris topologies]
Let\label{L:Vietoris} $\Xc$ be a metrisable space. Then the Vietoris topology on ${\rm Clos}(\Xc)$ induces on $\Ki_+(\Xc)$ the Hausdorff topology.
\el

\bpro
Since $\big\{A:A\sub O_1\big\}\cap\big\{A:A\sub O_2\big\}=\big\{A:A\sub O_1\cap O_2\big\}$, a basis for the Vietoris topology is formed by all sets of the form
\be\label{Vbasis}
\big\{A:A\sub O_0,\ A\cap O_i\neq\emptyset\ \forall i=1,\ldots,n\big\}
\ee
where $O_0,\ldots,O_n$ are open subsets of $\Xc$. Let $\Ob$ denote the set of all sets of the form (\ref{Vbasis}), intersected with $\Ki_+(\Xc)$. Fix a metric $d$ on $\Xc$ and let
\be
B_\eps(K):=\big\{K'\in\Ki_+(\Xc):d_{\rm H}(K,K')<\eps\big\}
\ee
denote the open ball of radius $\eps>0$ around a set $K\in\Ki_+(\Xc)$ in the Hausdorff metric corresponding to $d$. We need to show that for each $K\in\Ki_+(\Xc)$:
\begin{itemize}
\item[{\rm I.}] For each $\eps>0$, there exists a set $\Oi\in\Ob$ such that $K\in\Oi\sub B_\eps(K)$.
\item[{\rm II.}] For each $\Oi\in\Ob$ such that $K\in\Oi$, there exists an $\eps>0$ such that $B_\eps(K)\sub\Oi$.
\end{itemize}
To prove I., we use the fact that $K$ is compact and hence totally bounded to find a finite set $S\sub K$ such that $d(x,S)<\eps/2$ for all $x\in K$. We also set $O:=\{x\in\Xc:d(x,K)<\eps/2\}$ and define $\Oi\in\Ob$ by
\be
\Oi:=\big\{K\in\Ki_+(\Xc):K\sub O,\ K\cap B_{\eps/2}(x)\neq\emptyset\ \forall x\in S\big\}.
\ee
Then clearly $K\in\Oi\sub B_\eps(K)$, proving I. To prove also II, let $O_0,\ldots,O_n$ be open subsets of $\Xc$ and assume that $K\in\Ki_+(\Xc)$ satisfies $K\sub O_0$ and $K\cap O_i\neq\emptyset$ for all $1\leq i\leq n$. Using the compactness of $K$ it is easy to see that $\eps_0:=d(K,O_0^{\rm c})>0$, where $O_0^{\rm c}$ denotes the complement of $O_0$.
For $1\leq i\leq n$ we can moreover choose $x_i\in K\cap O_i$ and $\eps_i>0$ such that $B_{\eps_i}(x_i)\sub O_i$, where $B_\eps(x)$ denotes the open ball of radius $\eps$ around $x$ in $\Xc$. Then II is satisfied for $\eps:=\eps_0\wedge\cdots\wedge\eps_n$.
\epro

\subsection{Betweenness}\label{S:between}

In this subsection we prove Lemmas \ref{L:partdef}, \ref{L:metbet}, and \ref{L:compbet}. We start by proving some elementary consequences of the axioms (\ref{A1})--(\ref{A4}).

\bl[Elementary properties]
Each\label{L:betwel} betweenness satisfies, for each $x,y,y',z\in\Xc$:
\begin{enumerate}\addtocounter{enumi}{4}
\item $\li x,x\re=\{x\}$,\label{A5}
\item $y\in\li x,z\re\ \volgt\ \li x,y\re\sub\li x,z\re$,\label{A6}
\item $x\in\li y,z\re$ and $y\in\li x,z\re\ \volgt\ x=y$.\label{A7}
\item $y,y'\in\li x,z\re\mbox{ and }y'\in\li x,y\re\ \volgt\ y\in\li y',z\re$.\label{A8}
\end{enumerate}
\el

\bpro
Clearly, (\ref{A2}) and (\ref{A3}) imply (\ref{A5}) and (\ref{A4}) implies (\ref{A6}). To prove (\ref{A7}), we first observe that $x\in\li y,z\re$ and $y\in\li x,z\re$ imply by (\ref{A6}) $\li x,z\re\sub\li y,z\re\sub\li x,z\re$ and hence $\li x,z\re=\li y,z\re$. Using this as well as the assumptions $y\in\li x,z\re$ and $x\in\li y,z\re$ we can conclude by (\ref{A1}) and (\ref{A3}) that $\{y\}=\li x,y\re\cap\li y,z\re=\li y,x\re\cap\li x,z\re=\{x\}$ and hence $x=y$.

To prove (\ref{A8}), assume that $y,y'\in\li x,z\re$ and $y'\in\li x,y\re$. The statement is trivial if $y=y'$ so without loss of generality we assume that $y\neq y'$. Since $y'\in\li x,z\re$ we have by (\ref{A4}) that $\li x,z\re=\li x,y'\re\cup\li y',z\re$. Since also $y\in\li x,z\re$ we must have either $y\in\li x,y'\re$, or $y\in\li y',z\re$, or both. The first possibility would by (\ref{A1}) and (\ref{A7}) and the fact that $y'\in\li x,y\re$ imply that $y=y'$, which contradicts our assumptions, so we conclude that $y\in\li y',z\re$.
\epro

\bpro[of Lemma~\ref{L:partdef}]
The first implication $\volgt$ in (\ref{partdef}) follow from the fact that $y\in\li x,y\re$ by (\ref{A1}) and (\ref{A2}), while the reverse implication follows from (\ref{A6}). The second equivalence in (\ref{partdef}) follows from (\ref{A8}) and the third equivalence follows from the first one, by the symmetry (\ref{A1}). It is clear that (\ref{partdef}) defines a partial order $\leq_{x,z}$ on $\li x,z\re$. By (\ref{A4}), if $y,y'\in\li x,z\re$, then at least one of the conditions $y'\in\li x,y\re$ and $y'\in\li y,z\re$ must hold, which shows that $\leq_{x,z}$ is a total order.
\epro

Lemma~\ref{L:partdef} has a useful consequence.

\bl[Subsegments]
For\label{L:subseg} any betweenness, if $y\in\li x,z\re$, then $\li y,z\re=\{y'\in\li x,z\re:y\leq_{x,z}y'\}$, and $y',y''\in\li y,z\re$ satisfy $y'\leq_{y,z}y''$ if and only if $y'\leq_{x,z}y''$.
\el

\bpro
By Lemma~\ref{L:partdef} $y'\in\li x,z\re$ satisfies $y\leq_{x,z}y'$ if and only if $y'\in\li y,z\re$. This proves that $\li y,z\re=\{y'\in\li x,z\re:y\leq_{x,z}y'\}$. Again by Lemma~\ref{L:partdef}, $y',y''\in\li y,z\re$ satisfy $y'\leq_{y,z}y''$ if and only if $\li y'',z\re\sub\li y',z\re$, which in turn is equivalent to $y'\leq_{x,z}y''$.
\epro

\bpro[of Lemma~\ref{L:metbet}]
We need to check that our definition satisfies axioms (\ref{A1})--(\ref{A4}) of a betweenness. Axioms (\ref{A1}) and (\ref{A2}) are trivial. To prove (\ref{A3}) and (\ref{A4}), set $r:=d(x,z)$ and let $\ga\cn[0,r]\to\Xc$ be the unique isometry such that $\ga(0)=x$ and $\ga(r)=z$. Since an isometry is one-to-one, there exists a unique $p\in[0,r]$ such that $\ga(p)=y$. Clearly, the restrictions of $\ga$ to $[0,p]$ and $[p,r]$ are isometries, so $\li x,y\re=\{\ga(t):0\leq t\leq p\}$ and $\li y,z\re=\{\ga(t):p\leq t\leq r\}$. From these observations, axioms (\ref{A3}) and (\ref{A4}) follow immediately.
\epro

We say that a betweenness on a metrisable space $\Xc$ is \emph{generated by an interpolation function} if there exists a continuous function $\vhi\cn\Xc^2\times[0,1]\to\Xc$ that satisfies $\vhi(x,z,0)=x$, $\vhi(x,z,1)=z$, and 
\be\label{intfun}
\li x,z\re=\big\{\vhi(x,z,p):p\in[0,1]\big\}\qquad(x,z\in\Xc).
\ee
The following lemma prepares for the proof of Lemma~\ref{L:compbet}.

\bl[Proper betweennesses]
If\label{L:lincomp} $\Xc$ is a metrisable space, then the trivial betweenness is proper. The same is true for any betweenness that is generated by an interpolation function. If $\Xc$ is a closed subset of $\R$, then the order betweenness on $\Xc$ is proper.
\el

\bpro
For the trivial betweenness, $\li x,z\re=\{x,z\}$ is clearly compact for each $x,z\in\Xc$, and the continuity of the map $(x,z)\mapsto\{x,z\}$ with respect to the Hausdorff topology follows immediately from Lemma~\ref{L:Hauconv}.

If a betweenness is generated by an interpolation function, then $\li x,z\re$, being the image of $[0,1]$ under the continuous map $p\mapsto\vhi(x,z,p)$, is clearly compact for all $x,z\in\Z$. Let $x_n\to x$ and $z_n\to z$. To show that $\li x_n,z_n\re\to\li x,z\re$ in the Hausdorff topology, we check the conditions of Lemma~\ref{L:Hauconv}. Since $x_n\to x$ and $z_n\to z$, the sets $A:=\{x\}\cup\{x_n:n\in\N\}$ and $B:=\{z\}\cup\{z_n:n\in\N\}$ are compact. Let $\vhi(A\times B\times[0,1])$ denote the image of $A\times B\times[0,1]$ under $\vhi$, which is compact. Clearly $\li x_n,z_n\re\sub\vhi(A\times B\times[0,1])$ for all $n$. To complete the argument, it suffices to show that
\be\ba{l}\label{liconv}
\dis\big\{y\in\Xc:\exists y_n\in\li x_n,z_n\re\mbox{ s.t.\ $y$ is a cluster point of }(y_n)_{n\in\N}\big\}\\[5pt]
\dis\quad\sub\li x,z\re\sub\big\{y\in\Xc:\exists y_n\in\li x_n,z_n\re\mbox{ s.t.\ }y_n\to y\big\}.
\ec
For the first inclusion, assume that $y$ is a cluster point of $y_n=\vhi(x_n,z_n,p_n)$. By going to a subsequence, we can assume that $p_n\to p$ for some $p\in[0,1]$. Then $y=\vhi(x,z,p)\in\li x,z\re$. For the second inclusion, assume that $y=\vhi(x,z,p)$ for some $p\in[0,1]$. Then $y_n:=\vhi(x_n,z_n,p)\in\li x_n,z_n\re$ converge to $y$, completing the proof that each betweenness that is generated by an interpolation function is proper.

For the final statement of the lemma, assume that $\Xc$ is a closed subset of $\R$. Then clearly $\li x,z\re:=[x,z]\cap\Xc$ is compact for each $x,z\in\R$. Let $x_n\to x$ and $z_n\to z$. To show that $\li x_n,z_n\re\to\li x,z\re$ in the Hausdorff topology, we again check the conditions of Lemma~\ref{L:Hauconv}. Clearly, $\li x_n,z_n\re\sub[S,T]\cap\Xc$ for each $n$, where $S:=\inf_nx_n$ and $T:=\sup_nz_n$, so to complete the argument, it again suffices to check (\ref{liconv}). For the first inclusion, assume that $y\in\Xc$ is a cluster point of $y_n\in\li x_n,z_n\re$. Since $x_n\leq y_n\leq z_n$ for each $n$, taking the limit, we see that $x\leq y\leq z$ and hence $y\in\li x,z\re$. For the second inclusion, assume that $y\in\li x,z\re$. If $x<y<z$, then $x_n<y<z_n$ for all $n$ large enough, so setting $y_n:=y$ for $n$ large enough and $y_n:=x_n$ otherwise proves that $y$ is al element of the set on the right-hand side of (\ref{liconv}). If $y\in\{x,z\}$, then setting $y_n:=x_n$ or $:=z_n$ proves the same claim, so the proof is complete.
\epro

To prove Lemma~\ref{L:compbet} we need one more result.

\bl[Interpolation functions]
If\label{L:interpol} $\Xc$ is a normed linear space, then the linear betweenness is generated by an interpolation function. If $\Xc$ is a proper metric space with unique geodesics, then the same is true for the geodesic betweenness.
\el

We first show how this implies Lemma~\ref{L:compbet} and then prove Lemma~\ref{L:interpol}.\med

\bpro[of Lemma~\ref{L:compbet}]
Parts (i) and (iii) follow directly from Lemma~\ref{L:lincomp} while parts (ii) and (iv) follow by combining Lemma~\ref{L:lincomp} with Lemma~\ref{L:interpol}.
\epro

We next set out to prove Lemma~\ref{L:interpol}. The statement about the linear betweenness is trivial, but before we can prove the statement about the geodesic betweenness, we first need a better understanding of metric spaces with unique geodesics, which is provided by Proposition~\ref{P:geouni} below. In any metric space $(\Xc,d)$, for all $x,z\in\Xc$ and $\eps\geq 0$, we define
\be
\eta_{x,z}(\eps):=\sup\big\{d(y_1,y_2):\big[d(x,y_1)\wedge d(x,y_2)\big]+\big[d(y_1,z)\wedge d(y_2,z)\big]\leq d(x,z)+\eps\big\}.
\ee
In other words, this is the largest distance between two points $y_1,y_2\in\Xc$ for which there exist constants $r,r'\geq 0$ with $r+r'\leq d(x,z)+\eps$ such that $d(x,y_i)\leq r$ and $d(y_i,z)\leq r'$ $(i=1,2)$.

\bp[Unique geodesics]
Let\label{P:geouni} $(\Xc,d)$ be a metric space. Consider the following conditions.
\begin{itemize}
\item[{\rm(i)}] For all $x,z\in\Xc$ and $r,r'\geq 0$ with $r+r'=d(x,z)$, there exists an $y\in\Xc$ such that $d(x,y)=r$ and $d(y,z)=r'$.
\item[{\rm(ii)}] $\dis\eta_{x,z}(0)=0$ for all $x,z\in\Xc$.
\item[{\rm(ii)'}] $\dis\lim_{\eps\to 0}\eta_{x,z}(\eps)=0$ for all $x,z\in\Xc$.
\end{itemize}
Then $(\Xc,d)$ has unique geodesics if and only if (i) and (ii) hold. Moreover, (ii)' implies (ii), and if $(\Xc,d)$ is a proper metric space, then (ii)' implies (ii).
\ep

\bpro
In any metric space $(\Xc,d)$, let us introduce the notation
\be
\li x,z\re:=\big\{y\in\Xc:d(x,y)+d(y,z)=d(x,z)\big\}.
\ee
We claim that
\be\label{yyy}
y\in\li x,z\re,\ y'\in\li x,y\re,\ y''\in\li y,z\re\quad\volgt\quad y\in\li y',y''\re.
\ee
To see this, we note that if the assumptions in (\ref{yyy}) hold but the conclusion does not, then by the triangle inequality
\be\bac
\dis d(x,z)=d(x,y)+d(y,z)&=&\dis d(x,y')+d(y',y)+d(y,y'')+d(y'',z)\\[5pt]
&>&\dis d(x,y')+d(y',y'')+d(y'',z),
\ec
which contradicts the triangle inequality. Now let $(\Xc,d)$ be a metric space with unique geodesics and let $\lli x,z\rre$ denote the unique geodesic with endpoints $x,z$. Clearly $\lli x,z\rre\sub\li x,z\re$. We claim that
\be\label{gconcat}
y\in\li x,z\re\quad\volgt\quad\lli x,y\rre\cup\lli y,z\rre=\lli x,z\rre.
\ee
To see this, let $r:=d(x,y)$, $r':=d(y,z)$, and let $\ga\cn[0,r]\to\Xc$ and $\ga''\cn[r,r+r']\to\Xc$ be the unique isometries with $\ga(0)=x$, $\ga(r)=\ga'(r)=y$, and $\ga'(r+r')=z$. We claim that $\ga''\cn[0,r+r']\to\Xc$ defined as $\ga''(t)=\ga(t)$ for $t\in[0,r]$ and $:=\ga'(t)$ for $t\in[r,r+r']$ is an isometry. So see this, let $0\leq t'<t''\leq r+r'$. We need to show that $d\big(\ga''(t'),\ga''(t'')\big)=t''-t'$. This is clear when $t''\leq r$ or $r\leq t''$, while in the remaining case $t'<r<t''$ the claim follows from (\ref{yyy}).

We now prove that if $(\Xc,d)$ is a metric space with unique geodesics, then conditions (i) and (ii) are satisfied. Condition~(i) is trivial. To prove (ii), let $x,z\in\Xc$, let $r,r'\geq 0$ satisfy $r+r':=d(x,z)$, and assume that $y_1,y_2\in\Xc$ satisfy $d(x,y_i)=r$, $d(y_i,z)=r'$ $(i=1,2)$. By (\ref{gconcat}), there exist isometries $\ga_i\cn[0,r+r']\to\Xc$ with $\ga_i(0)=x$, $\ga_i(r)=y_i$, and $\ga_i(r+r')=z$, so by the assumption that $(\Xc,d)$ has unique geodesics we conclude that $y_1=y_2$, proving~(ii).

Conversely, if $(\Xc,d)$ is a metric space for which (i) and (ii) hold, then for each $x,z\in\Xc$ with $r:=d(x,z)$, we can uniquely define $\ga_{x,z}\cn[0,r]\to\Xc$ by
\be\label{gaxz}
\ga_{x,z}(t):=y\quad\mbox{with}\quad d(x,y)=t,\ d(y,z)=r-t.
\ee
Clearly, if $\ga\cn[0,r]\to\Xc$ is an isometry with $\ga(0)=x$ and $\ga(r)=z$, then we must have $\ga=\ga_{x,z}$, so to prove that $(\Xc,d)$ has unique geodesics, it suffices to show that $\ga_{x,z}$ is an isometry. Let $0\leq t_1\leq t_2\leq r$ and let $y_i:=\ga_{x,z}(t_i)$ $(i=1,2)$. Set $y'_1:=\ga_{x,y_2}(t_1)$. Then $d(x,y'_1)=t_1$ and $d(y'_1,z)\leq d(y'_1,y_2)+d(y_2,z)=r-t_1$, which by the assumption (ii) implies $y'_1=y_1$. Since $d(y'_1,y_2)=t_2-t_1$, this proves that $\ga_{x,z}$ is an isometry. This completes the proof that a metric space $(\Xc,d)$ has unique geodesics if and only if (i) and (ii) hold.

Trivially, (ii)' implies (ii), so to complete the proof of the proposition, it suffices to prove that for proper metric spaces, (ii) implies (ii)'. Assume that (ii)' does not hold for some $x,z\in\Xc$. Let $0<\eps_n\leq 1$ satisfy $\eps_n\to 0$. Then for some $\de>0$, we can find $y^n_1,y^n_2\in\Xc$ with $d(y^n_1,y^n_2)\geq\de$, as well as $r_n,r'_n\geq 0$ with $r_n+r'_n\leq d(x,z)+\eps_n$, such that $d(x,y^n_i)\leq r_n$ and $d(y^n_i,z)\leq r'_n$ $(i=1,2)$. Since $d(x,y^n_i)\leq d(x,z)+1$ for all $n$, by the properness assumption, we can select a subsequence such that $y^n_i\to y_i$ for some $y_1,y_2\in\Xc$. Since $r_n+r'_n\leq d(x,z)+1$, by going to a further subsequence, we can assume that $r_n\to r$ and $r'_n\to r'$ for some $r,r'\geq 0$ with $r+r'=d(x,z)$. Then $d(x,y_i)\leq r$ and $d(y_i,z)\leq r'$ while $d(y_1,y_2)\geq\de$ which shows that $\eta_{x,z}(0)\geq\de$, violating (ii).
\epro

If a metric space $(\Xc,d)$ has unique geodesics, then by conditions (i) and (ii) of Proposition~\ref{P:geouni}, we can uniquely define a function $\vhi\cn\Xc^2\times[0,1]\to\Xc$ by
\be\label{phigeo}
\vhi(x,z,p):=y\quad\mbox{with}\quad d(x,y)=pd(x,z)\mbox{ and }d(y,z)=(1-p)d(x,z).
\ee
Lemma~\ref{L:interpol} is now implied by Proposition~\ref{P:geouni} and the following lemma (the statement in Lemma~\ref{L:interpol} about normed linear spaces being trivial).

\bl[Geodesic interpolation function]
Let $(\Xc,d)$ be a metric space with unique geodesics and let $\vhi$ be defined as in (\ref{phigeo}). Then for each $x,z\in\Xc$, the unique geodesic with endpoints $x,z$ is given by $\{\vhi(x,z,p):p\in[0,1]\}$. If condition (ii)' of Proposition~\ref{P:geouni} is satisfied, then $\vhi\cn\Xc^2\times[0,1]\to\Xc$ is continuous.
\el

\bpro
Let $x,z\in\Xc$ and $r:=d(x,z)$. We observe that $\{\vhi(x,z,p):p\in[0,1]\}=\big\{\ga_{x,z}(t):t\in[0,r]\}$ where $\ga_{x,z}$ is defined as in (\ref{gaxz}). It has already been shown in the proof of Proposition~\ref{P:geouni} that this is the unique geodesic with endpoints $x,z$. Therefore, to complete the proof, it suffices to show that condition (ii)' of Proposition~\ref{P:geouni} implies that $\vhi$ is continuous.

Assume that $x_n,x,z_n,z\in\Xc$ and $p_n,p\in[0,1]$ satisfy $x_n\to x$, $z_n\to z$, and $p_n\to p$. Set $y_n:=\vhi(x_n,z_n,p_n)$ and $y:=\vhi(x,z,p)$. We have to show that $y_n\to y$. We observe that
\bc
\dis d(x,y_n)+d(y_n,z)
&\leq&\dis d(x_n,y_n)+d(y_n,z_n)+d(x,x_n)+d(z,z_n)\\[5pt]
&=&\dis d(x_n,z_n)+d(x,x_n)+d(z,z_n)\asto{n}d(x,z).
\ec
Thus, for each $\eps>0$, we can find an $m$ such that $d(x,y_n)+d(y_n,z)\leq d(x,z)+\eps$ for all $n\geq m$. Since moreover $d(x,y)+d(y,z)=d(x,z)$, it follows that $d(y_n,y)\leq\eta_{x,z}(\eps)$ for all $n\geq m$. Since $\eps>0$ is arbitrary, by (ii)', this implies $d(y_n,y)\to 0$.
\epro

We conclude this subsection with the following lemma that will be needed in Subsection~\ref{S:cad} below.

\bl[Segments as ordered sets]
Let~\label{L:betord} $\Xc$ be a metrisable space that is equipped with a proper betweenness. Then for each $x,z\in\Xc$, the segment $\li x,z\re$ equipped with the total order $\leq_{x,z}$ is an element of $\Ki_{\rm tot}(\Xc)$, and the map $(x,z)\mapsto\li x,z\re$ is continuous with respect to the product topology on $\Xc^2$ and the topology on $\Ki_{\rm tot}(\Xc)$.
\el

\bpro
To show that $\li x,z\re$ is an element of $\Ki_{\rm tot}(\Xc)$, we must show that the total order $\leq_{x,z}$ is compatible with the induced topology on $\li x,z\re$. Assume that $y_n,y'_n,y,y'\in\li x,z\re$ satisfy $y_n\to y$, $y'_n\to y'$, and $y_n\leq_{x,z}y'_n$ for all $n$. Then $y_n\in\li x,y'_n\re$ for all $n$. Since the betweenness is compatible with the topology, $\li x,y'_n\re\to\li x,y'\re$ in the Hausdorff topology, which by Lemma~\ref{L:Hauconv} implies that $y\in\li x,y'\re$ and hence $y\leq_{x,z}y'$. This shows that the total order $\leq_{x,z}$ is compatible with the induced topology on $\li x,z\re$.

To show that the map $(x,z)\mapsto\li x,z\re$ is continuous with respect to the topology on $\Ki_{\rm tot}(\Xc)$, assume that $x_n\to x$, $z_n\to z$. We will show that
\be
d_{\rm part}\big(\li x_n,z_n\re,\li x,z\re\big)\asto{n}0,
\ee
which is equivalent to the statement that $\li x_n,z_n\re^{\li 2\re}$ converges to $\li x,z\re^{\li 2\re}$ in the Hausdorff topology on $\Ki_+(\Xc^2)$. We apply Lemma~\ref{L:Hauconv}. Since $\li x_n,z_n\re\to\li x,z\re$, there exists a compact $C\sub\Xc$ such that $\li x_n,z_n\re\sub C$ for all $n$ and hence $\li x_n,z_n\re^{\li 2\re}\sub C^2$ for all $n$. Thus, it suffices to check that (compare (\ref{liconv}))
\be\ba{l}\label{liconv2}
\dis\big\{\!(y,y')\!\in\!\Xc^2:\exists y_n,y'_n\in\li x_n,z_n\re\mbox{ with }y_n\leq_{x_n,z_n}y'_n\mbox{ s.t.\ $(y,y')$ is a cluster point of }(y_n,y'_n)_{n\in\N}\!\big\}\\[5pt]
\dis\quad\sub\li x,z\re^{\li 2\re}\sub\big\{(y,y')\in\Xc:\exists y_n,y'_n\in\li x_n,z_n\re\mbox{ with }y_n\leq_{x_n,z_n}y'_n\mbox{ s.t.\ }(y_n,y'_n)\to(y,y')\big\}.
\ec
If $(y,y')$ is an element of the set on the left-hand side of (\ref{liconv2}) and $(y_n,y'_n)$ fulfill the conditions of the definition of this set, then by going to a subsequence we may assume that $(y_n,y'_n)\to(y,y')$. Then $y,y'\in\li x,z\re$ since $\li x_n,z_n\re\to\li x,z\re$. Moreover $y_n\leq_{x_n,z_n}y'_n$ means $y_n\in\li x_n,y'_n\re$. Since $y_n\to y$ and $\li x_n,y'_n\re\to\li x,y'\re$, this implies $y\in\li x,y'\re$ and hence $y\leq_{x,z}y'$, proving that $(y,y')\in\li x,z\re^{\li 2\re}$.

To prove the second inclusion in (\ref{liconv2}), assume that $(y,y')\in\li x,z\re^{\li 2\re}$. Since $\li x_n,z_n\re\to\li x,z\re$, there exist $y_n,y'_n\in\li x_n,z_n\re$ such that $y_n\to y$ and $y'_n\to y'$. We now distinguish two cases: $y\neq y'$ and $y=y'$. If $y\neq y'$, then we claim that $y_n\leq_{x_n,z_n}y'_n$ for all $n$ large enough. Indeed, in the opposite case, since $\leq_{x_n,z_n}$ is a total order, by going to a subsequence, we can assume that $y'_n\leq_{x_n,z_n}y_n$ for all $n$, which by the arguments we have already seen implies $y'\leq_{x,z}y$, so that by the fact that $(y,y')\in\li x,z\re^{\li 2\re}$ we must have $y=y'$, contradicting our assumption. Since $y_n\leq_{x_n,z_n}y'_n$ for all $n$ large enough, changing the definitions of $y_n,y'_n$ for finitely many $n$, we see that there exist $y_n,y'_n\in\li x_n,z_n\re$ with $y_n\leq_{x_n,z_n}y'_n$ such that s.t. $(y_n,y'_n)\to(y,y')$. In the case $y=y'$ the argument is even simpler, since now $(y_n,y_n)\to(y,y')$ while obviously $y_n\leq_{x_n,z_n}y_n$ for all $n$.
\epro

\subsection{Cadlag functions}\label{S:cad}

In this subsection we prove Lemmas \ref{L:fpi} and \ref{L:cogra} as well as some related results that will be needed later. Let $\Xc$ be a metrisable space. In Subsection~\ref{S:paths}, we defined a cadlag function to be a triple $(I,f^-,f^+)$, where $I\sub\R$ is a closed set and $f^-,f^+\cn I\to\Xc$ are left and right continuous functions that satisfy (\ref{cadlag}). By Lemma~\ref{L:cadlag}, we can equivalently define a cadlag function to be a pair $(I,f)$ such that $I\sub\R$ is closed and $f\cn I_\sg\to\Xc$ is continuous, where $I_\sg$ is the subset of the split real line defined in (\ref{Isg}). We work with this equivalent definition from now on and moreover define an \emph{extended cadlag function} to be a pair $(\hat I,f)$ such that $\hat I$ is a nonempty closed subset of $\ov\R$ and setting $I:=\hat I\cap\R$, we have that $(I,f)$ is a cadlag function. We define the \emph{closed graph} $\Gi(\hat I,f)$ of an extended cadlag function $(\hat I,f)$ as
\be\label{Gi}
\Gi(\hat I,f):=\big\{(x,t):t\in I,\ x\in\{f(t-),f(t+)\}\big\}\cup\big\{(\ast,\pm\infty):\pm\infty\in\hat I\big\}.
\ee
If $\Xc$ is equipped with a betweenness (see Subsection~\ref{S:between}), then we define the \emph{filled graph}\footnote{Except for the points at infinity, this is what Whitt \cite{Whi02} calls the \emph{completed graph}.} of an extended cadlag function $(\hat I,f)$ as
\be\label{Gint}
\Gi_{\rm f}(\hat I,f):=\big\{(x,t):t\in I,\ x\in\li f(t-),f(t+)\re\big\}\cup\big\{(\ast,\pm\infty):\pm\infty\in\hat I\big\}.
\ee
Note that for the trivial betweenness, the filled and closed graphs coincide. The filled graph is naturally equipped with a total order, which is defined by setting $(x,s)\prec(y,t)$ if either $s<t$ and $x,y$ are arbitrary, or $s=t$ and $x<_{f(t-),f(t+)}y$, where $\leq_{f(t-),f(t+)}$ is the total order on the segment $\li f(t-),f(t+)\re$ defined in Lemma~\ref{L:partdef}.

\bl[Filled graphs]
Assume\label{L:graph} that $\Xc$ is a metrisable space that is equipped with a proper betweenness. Then for any extended cadlag function $(\hat I,f)$, the filled graph $\Gi_{\rm f}(\hat I,f)$ is a compact subset of the squeezed space $\Ri(\Xc)$, and the total order $\pre$ is compatible with the induced topology on $\Gi_{\rm f}(\hat I,f)$.
\el

\bpro
We will show that each sequence $(x_n,t_n)\in\Gi_{\rm f}(\hat I,f)$ has a subsequence that converges to a limit in $\Gi_{\rm f}(\hat I,f)$. Since $\hat I$ is compact, we can select a subsequence such that $t_n\to t$ for some $t\in\hat I$. If $t=\pm\infty$, then Lemma~\ref{L:sqz} tells us that $(x_n,t_n)\to(\ast,\pm\infty)\in\Gi_{\rm f}(\hat I,f)$ so we are done, so from now on we can assume that $t\in\R$. By going to a further subsequence, we can assume that we are in one of the following three cases: (i) $t_n<t$ for all $n$, (ii) $t_n>t$ for all $n$, and (iii) $t_n=t$ for all $n$. In case~(i), we use the continuity of $f\cn I_\sg\to\Xc$ and the fact that the betweenness is compatible with the topology to see, using Lemma~\ref{L:betwel}~(\ref{A5}), that
\be
x_n\in\li f(t_n-),f(t_n+)\re\asto{n}\li f(t-),f(t-)\re=\{f(t-)\},
\ee
where the convergence is in $\Ki_+(\Xc)$. We conclude from this that $(x_n,t_n)$ converges to $\big(f(t-),t\big)\in\Gi_{\rm f}(\hat I,f)$. In case~(ii), the same argument shows that $(x_n,t_n)$ converges to $\big(f(t+),t\big)\in\Gi_{\rm f}(\hat I,f)$. In case~(iii), finally, using the compactness of $\li f(t-),f(t+)\re$, we can select a further subsequence such that $(x_n,t)\to(x,t)$ for some $x\in\li f(t-),f(t+)\re$. Since also in this case the limit $(x,t)$ is an element of $\Gi_{\rm f}(\hat I,f)$, we are done.

To see that the total order $\pre$ on $\Gi_{\rm f}(\hat I,f)$ is compatible with the (induced) topology on $\Gi_{\rm f}(\hat I,f)$, it suffices to show that
\be\label{strct}
S:=\big\{\big((x,s),(y,t)\big)\in\Gi_{\rm f}(\hat I,f)^2:(x,s)\prec(y,t)\big\}
\ee
is an open subset of $\Gi_{\rm f}(\hat I,f)^2$. If $\big((x,s),(y,t)\big)$ is an element of $S$, then either: (i) $s<t$, or: (ii) $s=t\in\R$ and $x\neq y$. In case (i), we can choose $s<S<T<t$. Then
\be
O:=\big\{\big((x',s'),(y',t')\big)\in\Gi_{\rm f}(\hat I,f)^2:s'<S,\ T<t'\big\}
\ee
is an open subset of $\Gi_{\rm f}(\hat I,f)^2$ such that $\big((x,s),(y,t)\big)\in O\sub S$. In case (ii), we recall that by definition $(x,t)\pre(y,t)$ if $x\leq_{f(t-),f(t+)}z$, where by Lemma~\ref{L:betord} the total order $\leq_{f(t-),f(t+)}$ on $\li f(t-),f(t+)\re$ is compatible with the topology. It follows that we can choose $\eps>0$ small enough such that for $z\in\li f(t-),f(t+)\re$, if $d(z,x)<\eps$ then $(z,t)\prec(y,t)$, while if $d(z,y)<\eps$ then $(x,t)\prec(z,t)$. Next, we use the continuity of $f\cn I_\sg\to\Xc$ to choose $\de>0$ small enough such that $d\big(f(s\pm),x\big)>\eps$ for all $t<s<t+\de$ and $d\big(f(s\pm),y\big)>\eps$ for all $t-\de<s<t$. Then
\be
O:=\big\{\big((x',s),(y',u)\big)\in\Gi_{\rm f}(\hat I,f):|s-t|\vee|u-t|<\de,\ d(x',x)\vee d(y',y)<\eps\big\}
\ee
is an open subset of $\Gi_{\rm f}(\hat I,f)^2$ such that $\big((x,s),(y,t)\big)\in O\sub S$. Together, these observations prove that $S$ is an open subset of $\Gi_{\rm f}(\hat I,f)^2$.
\epro

Recall that for any metrisable space $\Xc$, in Subsection~\ref{S:orHa} we define $\Ki_{\rm tot}(\Xc)$ to be the set of pairs $(K,\pre)$ such that $K$ is a compact subset of $\Xc$ and $\pre$ is a total order on $K$ that is compatible with the induced topology on $K$.

\bl[Characterisation of filled graphs]
Let\label{L:intchar} $\Xc$ be a metrisable space that is equipped with a betweennesss that is compatible with the topology. Assume that $(G,\pre)\in\Ki_{\rm tot}(\Ri(\Xc))$. Then $(G,\pre)$ is the filled graph of an extended cadlag function $(\hat I,f)$ if and only if the following conditions are satisfied.
\begin{enumerate}
\item $(x,s)\pre(y,t)$ for all $(x,s),(y,t)\in G$ such that $s<t$.
\item For each $t\in\R$ and $(x_1,t),(x_2,t),(x_3,t)\in G$ with $(x_1,t)\pre(x_2,t)\pre(x_3,t)$, one has $x_2\in\li x_1,x_2\re$.
\end{enumerate}
Moreover, the extended cadlag function $(\hat I,f)$ is uniquely determined by its filled graph $(G,\pre)$.
\el

\bpro
By Lemma~\ref{L:graph}, the filled graph of an extended cadlag function $(\hat I,f)$ corresponds to an element of $\Ki_{\rm tot}(\Ri(\Xc))$. Properties (i) and (ii) now follow from the definition of the total order $\pre$ on $\Gi_{\rm f}(\hat I,f)$ and Lemma~\ref{L:subseg}.

Assume that conversely, $(G,\pre)\in\Ki_{\rm tot}(\Xc)$ satisfies (i) and (ii). We define
\be
I:=\hat I\cap\R\quad\mbox{with}\quad\hat I:=\big\{t\in\ov\R:\exists x\in\Xc\cup\{\ast\}\mbox{ s.t.\ }(x,t)\in G\big\}.
\ee
Since $\hat I$ is the image of $G$ under the continuous map $\Ri(\Xc)\ni(x,t)\mapsto t\in\ov\R$, we see that $\hat I$ is compact. We claim that for each $t\in I$, there exist unique $f(t-),f(t+)\in\Xc$ such that $f(t-)\pre f(t+)$ and
\be
S_t:=\big\{x\in\Xc:(x,t)\in G\big\}=\li f(t-),f(t+)\re.
\ee
Indeed, $S_t$ is a compact metrisable set, so we can choose a countable dense set $\{x_n:n\in\N\}\sub S_t$. Set $y_0:=x_0$ and define $y_n$ as the maximum of $x_n$ and $y_{n-1}$ in the total order $\pre$ $(n\geq 1)$. By the compactness of $S_t$, by going to a subsequence, we can assume that $y_n\to f(t+)$ for some $f(t+)\in S_t$. Then $y'\pre y$ for all $y'\in D$ and hence also for all $y'\in S_t$ since the order is compatible with the topology. In the same way, we see that $S_t$ has a (necessarily unique) minimal element $f(t-)$. By (ii), we conclude that $S_t=\li f(t-),f(t+)\re$.

We can now use the claim we have just proved to define $f\cn I_\sg\to\Xc$ by
\be
\li f(t-),f(t+)\re:=\{x\in\Xc:(x,t)\in G\}
\quad\mbox{with}\quad
f(t-)\pre f(t+)\qquad\big(t\in I\big)
\ee
which has the result that $(G,\pre)$ is the filled graph of $(\hat I,f)$. To show that $(\hat I,f)$ is an extended cadlag function, we need to show that $f\cn I_\sg\to\Xc$ is continuous. By symmetry, it suffices to show that if $\tau_n\in I_\sg$ and $t\in I$ satisfy $\un\tau_n>t$ for all $n$ and $\un\tau_n\to t$ as $n\to\infty$, then $f(\tau_n)\to f(t+)$. It suffices to show that $\{f(\tau_n):n\in\N\}$ is precompact and its only cluster point is $f(t+)$. Equivalently, we may show that each subsequence of $f(\tau_n)$ contains a further subsequence that converges to $f(t+)$. By the compactness of $G$, for any subsequence, we can select a further subsequence such that $f(\tau_n)\to x$ for some $x\in\Xc$ such that $(x,t)\in G$. By (i), we have $\big(f(t+),t\big)\pre\big(f(\tau_n),\un\tau_n\big)$ for all $n$, so using the fact that the total order is compatible with the topology, we see that $\big(f(t+),t\big)\pre(x,t)$, which using the fact that $f(t+)$ is the maximal element of $\li f(t-),f(t+)\re$ with respect to the order $\pre$ identifies $x$ as $f(t+)$.
\epro

\bl[Characterisation of paths]
Let\label{L:pathgra} $\Xc$ be a metrisable space, let $\pi$ be a subset of $\Ri(\Xc)$, and let $\pre$ be a total order on $\pi$. Then $(\pi,\pre)\in\Pi(\Xc)$ if and only if $(\pi,\pre)$ is the closed graph of an extended cadlag function $(\hat I,f)$.
\el

\bpro
In Section~\ref{S:paths}, we defined $\Pi(\Xc)$ to be the set of $(\pi,\pre)\in\Ki_{\rm tot}(\Ri(\Xc))$ such that:
\begin{itemize}
\item[{\rm(i)}] $(x,s)\pre(y,t)$ for all $(x,s),(y,t)\in\pi$ such that $s<t$,
\item[{\rm(ii)'}] for each $t\in\R$, the set $\{x\in\Xc:(x,t)\in\pi\}$ has at most two elements.
\end{itemize}
We observe that (ii)' is equivalent to condition (ii) of Lemma~\ref{L:intchar} for the trivial betweenness. In view of this, the claim follows from Lemma~\ref{L:intchar}.
\epro

\bl[Characterisation of filled paths]
Let\label{L:fillgra} $\Xc$ be a metrisable space that is equipped with a proper betweenness, let $\ov\pi$ be a subset of $\Ri(\Xc)$, and let $\pre$ be a total order on $\ov\pi$. Then $\ov\pi$ is the filled path associated with a path $\pi\in\Pi(\Xc)$ if and only if $(\ov\pi,\pre)$ is the filled graph of an extended cadlag function $(\hat I,f)$.
\el

\bpro
By Lemma~\ref{L:pathgra}, there is a one-to-one correspondence between paths $\pi\in\Pi(\Xc)$ and extended cadlag functions $(\hat I,f)$. The claim now follows by comparing the definition of the filled path $\ov\pi$ associated with a path $\pi$ in (\ref{fillpath}) with the definition of the filled graph of an extended cadlag function $(\hat I,f)$ in (\ref{Gint}).
\epro

\bpro[of Lemma~\ref{L:fpi}]
Immediate from Lemma~\ref{L:pathgra} and the definition of the closed graph of an extended cadlag function $(\hat I,f)$.
\epro

\bpro[of Lemma~\ref{L:cogra}]
Immediate from Lemmas \ref{L:intchar} and \ref{L:fillgra}.
\epro

By definition, an \emph{extended continuous function} is a pair $(\hat I,f)$ such that $\hat I$ is a nonempty closed subset of $\ov\R$ and $f\cn I\to\Xc$ is a continuous function, where $I:=\hat I\cap\R$. Equivalently, this is an extended cadlag function such that $f(t-)=f(t+)$ for all $t\in I$. The \emph{closed graph} of an extended continuous function is the set
\be
\Gi(\hat I,f):=\big\{(f(t),t):t\in I\big\}\cup\big\{(\ast,\pm\infty):\pm\infty\in\hat I\big\},
\ee
which coincides with the closed graph of $(\hat I,f)$, viewed as an extended cadlag function. We equip $\Gi(\hat I,f)$ with a total order $\pre$ by setting $(x,s)\prec(y,t)$ if $s<t$. Then Lemma~\ref{L:graph} tells us that $\Gi(\hat I,f)$ is compact and that the total order on $\Gi(\hat I,f)$ is compatible with the induced topology on $\Gi(\hat I,f)$. The following lemma is similar to Lemma~\ref{L:intchar}.

\bl[Characterisation of continuous graphs]
Let\label{L:conchar} $\Xc$ be a metrisable space. Assume that $G\in\Ki_+(\Ri(\Xc))$. Then $G$ is the closed graph of an extended continuous function $(\hat I,f)$ if and only if for each $t\in\R$, the set $\{x\in\Xc:(x,t)\in G\}$ has at most one element.
\el

\bpro
We have just observed that if $G=\Gi(\hat I,f)$ is the graph of an extended continuous function, then $G$ is compact. Clearly, the set $\{x\in\Xc:(x,t)\in G\}$ has at most one element for each $t\in\R$.

Conversely, if $G\in\Ki_+(\Ri(\Xc))$ has this property, then we define
\be
I:=\hat I\cap\R\quad\mbox{with}\quad\hat I:=\big\{t\in\ov\R:\exists x\in\Xc\cup\{\ast\}\mbox{ s.t.\ }(x,t)\in G\big\}.
\ee
Then $\hat I$ is compact since it is the continuous image of a compact set. We use the fact that $G$ contains at each time at most one point to define $f\cn I\to\Xc$ by
\be\label{piG}
\big\{f(t)\big\}:=\{x\in\Xc:(x,t)\in G\}\qquad(t\in I).
\ee
To show that $f$ is continuous, assume that $t_n,t\in I$ satisfy $t_n\to t$. By the same argument as in the proof of Lemma~\ref{L:graph}, using the compactness of $G$, it suffices to show that if $f(t_n)\to x$ along a subsequence then $x=f(t)$. This follows from the assumption that $\{x\in\Xc:(x,t)\in G\}$ has at most one element.
\epro

Continuous paths whose domain is an interval have the property that viewed as a subset of the squeezed space, they are connected. The following lemma says that this is an if and only if.

\bl[Connected graphs]
A\label{L:congraph} path $\pi\in\Pi(\Xc)$ is connected if and only if $\pi\in\Pi^|_{\rm c}(\Xc)$.
\el

\bpro
If $\pi\in\Pi_{\rm c}(\Xc)$, then $\pi$ is the image of the compact set $\hat I(\pi)$ under the continuous map from $\ov\R$ to $\Ri(\Xc)$ given by $t\mapsto\big(\pi(t),t)$ (with $\pm\infty\mapsto(\ast,\pm\infty)$). If $\pi\in\Pi^|_{\rm c}(\Xc)$, then $\hat I(\pi)$ is connected, so by the fact that the continuous image of a connected set is connected, we conclude that $\pi$ is connected.

Conversely, if $\pi$ is connected, then $\hat I(\pi)$ must be connected and hence $\pi\in\Pi^|(\Xc)$. To see that $\pi$ is moreover continuous, assume that conversely, $\pi(t-)\neq\pi(t+)$ for some $t\in I(\pi)$. Then we can define new paths $\pi',\pi''$ with extended domains $\hat I(\pi'):=[-\infty,t]\cap\hat I(\pi)$ and $\hat I(\pi''):=[t,\infty]\cap\hat I(\pi)$, by setting $\pi'(s\pm):=\pi(s\pm)$ and $\pi''(s\pm):=\pi(s\pm)$ for $s<t$ and $s>t$, respectively, and $\pi'(t\pm):=\pi(t-)$ and $\pi''(t\pm):=\pi(t+)$. By Lemma~\ref{L:graph}, $\pi'$ and $\pi''$ are compact sets. Since $\pi'\cap\pi''=\emptyset$ and $\pi'\cup\pi''=\pi$, this proves that $\pi$ is not connected.
\epro

%
%

\section{The ordered Hausdorff metric}\label{S:Hord}

\subsection{The ordered Hausdorff metric}\label{S:orH}

In this subsection, we study the metrics $d_{\rm part}$ and $d_{\rm tot}$ defined in (\ref{part}) and (\ref{tot}), preparing for the proofs of Theorems \ref{T:partot} and \ref{T:totcomp}. Let $(\Xc,d)$ be a metric space. Generalising the definition in (\ref{Kli}), for each $m\geq 1$ and $K\in\Ki_{\rm part}(\Xc)$, we set
\be\label{Km}
K^{\li m\re}:=\big\{(x_1,\ldots,x_m)\in K^m:x_1\pre\cdots\pre x_m\big\}.
\ee
It is straightforward to check that $K^{\li m\re}$ is a closed subset of $K^m$ and hence a compact subset of $\Xc^m$. Generalising the definition in (\ref{d2}), we equip $\Xc^m$ with the metric
\be\label{dm}
d^m\big((x_1,\ldots,x_m),(y_1,\ldots,y_m)\big):=\bigvee_{k=1}^md(x_k,y_k),
\ee
and we equip $\Ki_+(\Xc^m)$ with the associated Hausdorff metric $d^m_{\rm H}$. Generalising the definition in (\ref{part}), for each $m\geq 1$, we define a function $d^{\li m\re}$ on $\Ki_{\rm part}(\Xc)^2$ by
\be\label{dlmr}
d^{\li m\re}(K_1,K_2):=d^m_{\rm H}(K^{\li m\re}_1,K^{\li m\re}_2)\qquad\big(K_1,K_2\in\Ki_{\rm part}(\Xc)\big).
\ee
In particular, when $m\geq 2$, this is a metric on $\Ki_{\rm part}(\Xc)$ since $(K,\pre)$ is uniquely characterised by $K^{\li m\re}$ for $m\geq 2$. On the other hand, $d^{\li 1\re}(K_1,K_2)$ is simply the Hausdorff distance between $K_1$ and $K_2$ as sets, which gives no information about the partial order. The following lemma describes a simple property of the metric $d^{\li 2\re}$.

\bl[Ordered limit]
Let\label{L:orpres} $\Xc$ be a metrisable space. Assume that $K_n,K\in\Ki_{\rm part}(\Xc)$ satisfy $d^{\li 2\re}(K_n,K)\to 0$ and that $x_n,y_n\in K_n$ satisfy $x_n\to x$, $y_n\to y$, and $x_n\pre y_n$ for all $n$. Then $x,y\in K$ satisfy $x\pre y$.
\el

\bpro
Since $d^{\li 2\re}(K_n,K)\to 0$, we have $K_n^{\li 2\re}\to K^{\li 2\re}$ and hence by Lemma~\ref{L:Hauconv} $x,y\in K^{\li 2\re}$, which proves that $x\pre y$.
\epro

The folling lemma gives a one-sided bound between metrics of the form $d^{\li m\re}$ for different values of $m$.

\bl[One-sided bound]
One\label{L:monesid} has
\be
d^{\li m\re}(K_1,K_2)\leq d^{\li m+1\re}(K_1,K_2)\qquad\big(m\geq 1,\ K_1,K_2\in\Ki_{\rm part}(\Xc)\big).
\ee
\el

\bpro
By Lemma~\ref{L:Haus2}, there exists a correspondence $R$ between $K_1^{\li m+1\re}$ and $K_2^{\li m+1\re}$ such that $d^{m+1}(x,y)\leq d^{m+1}_{\rm H}(K^{\li m+1\re}_1,K^{\li m+1\re}_2)$ for all $(x,y)\in R$. Let $\psi\cn\Xc^{m+1}\to\Xc$ denote the projection $\psi(x_1,\ldots,x_{m+1}):=(x_1,\ldots,x_m)$. Then (\ref{dm}) implies that
\be
d^m\big(\psi(x),\psi(y)\big)\leq d^{m+1}(x,y)
\qquad(x,y\in\Xc^{m+1}).
\ee
Since $\psi(K_i^{\li m+1\re})=K_i^{\li m\re}$ $(i=1,2)$, it follows that $R':=\big\{\big(\psi(x),\psi(y)\big):(x,y)\in R\big\}$ is a correspondence between $K_1^{\li m\re}$ and $K_2^{\li m\re}$ such that $d^m(x',y')\leq d^{m+1}_{\rm H}(K^{\li m+1\re}_1,K^{\li m+1\re}_2)$ for all $(x',y')\in R'$. By Lemma~\ref{L:Haus2}, this proves that
\be
d^m_{\rm H}(K^{\li m\re}_1,K^{\li m\re}_2)\leq d^{m+1}_{\rm H}(K^{\li m+1\re}_1,K^{\li m+1\re}_2),
\ee
which in view of (\ref{dlmr}) implies the claim.
\epro

The following lemmas show that in general, the metrics $d^{\li m\re}$ for different values of $m$ are not comparable. More precisely, the one-sided bound in Lemma~\ref{L:monesid} is not matched by an opposite inequality of the form $d^{\li m+1\re}(K_1,K_2)\leq Cd^{\li m\re}(K_1,K_2)$ for any finite constant $C$, and convergence in $d^{\li m\re}$ does not imply convergence in $d^{\li m+1\re}$.

\bl[No opposite inequality]\hspace{-0.3pt}
Let\label{L:noop} $\Xc=[0,1]$, equipped with the usual distance. Then for each $m\geq 1$ and $0<\eps\leq 1/4$, there exist $K_1,K_2\in\Ki_{\rm tot}(\Xc)$ such that $d^{\li m\re}(K_1,K_2)\leq\eps$ and $d^{\li m+1\re}(K_1,K_2)\geq 1/2$.
\el

\bpro
We choose $K_1=\{x_1,\ldots,x_{m+1}\}$ with $x_k\in[0,\eps]$ when $k$ is even and $x_k\in[1-\eps,1]$ if $k$ is odd, and we choose $K_2=\{y_1,\ldots,y_{m+1}\}$ with $y_k\in[0,\eps]$ when $k$ is odd and $y_k\in[1-\eps,1]$ if $k$ is even. We equip $K_1$ and $K_2$ with total orders such that $x_1\prec\cdots\prec x_{m+1}$ and $y_1\prec\cdots\prec y_{m+1}$. It is easy to see that
\be
d\big((x_1,\ldots,x_{m+1}),K_2^{\li m+1\re}\big)\geq 1/2,
\ee
and hence $d^{\li m+1\re}(K_1,K_2)\geq 1/2$. On the other hand, it is easy to see that for each $(z_1,\ldots,z_m)\in K_1^{\li m\re}$, there exists a $(z'_1,\ldots,z'_m)\in K_2^{\li m\re}$ such that $|z_k-z'_k|\leq\eps$ for all $k$, and vice versa, so $d^{\li m\re}(K_1,K_2)\leq\eps$.
\epro

\bl[Different topologies]
Let\label{L:diftop} $\Xc=[0,1]$, equipped with the usual distance. Then for each $m\geq 1$, there exist $K_n\in\Ki_{\rm part}(\Xc)$ and $K\in\Ki_{\rm tot}(\Xc)$ such that $d^{\li m\re}(K_n,K)\to 0$ as $n\to\infty$ but $d^{\li m+1\re}(K_n,K)\geq 1/2$ for all $n$.
\el

\bpro
It will be convenient to use the notation $[m]:=\{1,\ldots,m\}$ $(m\geq 1)$. We choose $x_1,\ldots,x_{m+1}$, all different, with $x_k\in[0,1/4]$ when $k$ is even and $x_k\in[3/4,1]$ if $k$ is odd. We set $K=\{x_1,\ldots,x_{m+1}\}$ and we equip $K$ with a total order by setting $x_1\prec\cdots\prec x_{m+1}$. We choose points
\be
\big(x_k^l(n)\big)_{k\in[m+1]}^{l\in[m+1]}
\ee
in $[0,1]$, all different, such that $x_k^l(n)\to x_k$ as $n\to\infty$ for all $k,l\in[m+1]$, and we set
\be
K_n:=\big\{x_k^l:k,l\in[m+1],\ k\neq l\big\}.
\ee
We equip $K_n$ with a partial order such that
\be
x_k^l\pre x_{k'}^{l'}\quad\desd\quad k\pre k'\mbox{ and }l=l'.
\ee
Then it is easy to check that $d^{\li m\re}(K_n,K)\to 0$ as $n\to\infty$ but $d^{\li m+1\re}(K_n,K)\geq 1/2$ for all $n$.
\epro

We note that $d^{\li m\re}(K_1,K_2)\leq\sup_{(x_1,x_2)\in K_1\times K_2}d(x_1,x_2)$, which is finite by the continuity of $d$ and the compactness of $K_1\times K_2$. We use this and Lemma~\ref{L:monesid} to define $d^{\li\infty\re}$ as the increasing limit
\be
d^{\li\infty\re}(K_1,K_2):=\lim_{m\to\infty}d^{\li m\re}(K_1,K_2)\qquad\big(K_1,K_2\in\Ki_{\rm part}(\Xc)\big).
\ee
It is straightforward to check that $d^{\li\infty\re}$ is a metric on $\Ki_{\rm part}(\Xc)$: symmetry and the triangle inequality follow by taking the limit in the corresponding properties of the metrics $d^{\li m\re}$, and $d^{\li\infty\re}(K_1,K_2)=0$ clearly implies $d^{\li m\re}(K_1,K_2)=0$ for all $m\geq 1$ and hence equality of $K_1$ and $K_2$ as partially ordered spaces. In the special case that $K_1$ and $K_2$ are totally ordered, the following proposition identifies $d^{\li\infty\re}(K_1,K_2)$ as the metric $d_{\rm tot}$ defined in (\ref{tot}).

\bp[Monotone correspondences]
Let\label{P:moncor} $(\Xc,d)$ be a metric space. Then one has $d^{\li\infty\re}(K_1,K_2)=d_{\rm tot}(K_1,K_2)$ for all $K_1,K_2\in\Ki_{\rm tot}(\Xc)$.
\ep

The proof of Proposition~\ref{P:moncor} uses the following simple lemma.

\bl[Eventually ordered sequences]
Let\label{L:evor} $K$ be a compact metrisable set that is equipped with a total order that is compatible with the topology. Assume that $x\prec y$ and that $x_n,y_n\in K$ satisfy $x_n\to x$ and $y_n\to y$. Then $x_n\prec y_n$ for all $n$ sufficiently large.
\el

\bpro
Since $\pre$ is a total order, if the statement is not true, then $y_n\pre x_n$ for infinitely many $n$, so we can select a subsequence such that $y_n\pre x_n$ for all $n$. Taking the limit, using the fact that the total order that is compatible with the topology, we find that $y\pre x$, which contradicts $x\prec y$.
\epro

\bpro[of Proposition~\ref{P:moncor}]
We first prove the inequality $d^{\li\infty\re}(K_1,K_2)\leq d_{\rm tot}(K_1,K_2)$. Let $R$ be a monotone correspondence between $K_1$ and $K_2$. Let $x_1,\ldots,x_m\in K_1$ satisfy $x_1\pre\cdots\pre x_m$. Then we can choose $x'_1,\ldots,x'_m\in K_2$ such that $(x_k,x'_k)\in R$ for all $1\leq k\leq m$, and moreover $x'_k=x'_{k+1}$ whenever $x_k=x_{k+1}$ $(1\leq k<m)$. Since $R$ is monotone and $K_2$ is totally ordered, we must have $x'_1\pre\cdots\pre x'_m$. This shows that
\be
d^m\big((x_1,\ldots,x_m),K_2^{\li m\re}\big)\leq\sup_{(x,x')\in R}d(x,x')
\qquad\big((x_1,\ldots,x_m)\in K_1^{\li m\re}\big).
\ee
The same is true with the roles of $K_1$ and $K_2$ interchanged, so we conclude that
\be
d^{\li m\re}(K_1,K_2)=d^m_{\rm H}(K_1^{\li m\re},K_2^{\li m\re})\leq\sup_{(x,x')\in R}d(x,x').
\ee
Taking the infimum over all monotone correspondences between $K_1$ and $K_2$ and letting $m\to\infty$ we see that $d^{\li\infty\re}(K_1,K_2)\leq d_{\rm tot}(K_1,K_2)$.

To prove the opposite inequality, let $\eps_n$ be positive constants, tending to zero. Since $K_1$ is totally bounded, for each $n$, we can find an $m(n)\geq 1$ and $x^n_1,\ldots,x^n_{m(n)}\in K_1$ such that
\be\label{K1net}
d\big(x,\{x^n_1,\ldots,x^n_{m(n)}\}\big)\leq\eps_n\quad\forall x\in K_1.
\ee
Since $K_1$ is totally ordered, we can assume without loss of generality that $x^n_1\pre\cdots\pre x^n_{m(n)}$. Since
\be\label{dminf}
d^m_{\rm H}(K_1^{\li m\re},K_2^{\li m\re})=d^{\li m\re}(K_1,K_2)\leq d^{\li\infty\re}(K_1,K_2),
\ee
we can find $y^n_1,\ldots,y^n_{m(n)}\in K_2$ with $y^n_1\pre\cdots\pre y^n_{m(n)}$ such that
\be\label{und}
d(x^n_k,y^n_k)\leq d^{\li\infty\re}(K_1,K_2)\qquad(1\leq k\leq m(n)).
\ee
Using the fact that $K_2$ is totally bounded, adding points to $\{y^n_1,\ldots,y^n_{m(n)}\}$ and making $m(n)$ larger if necessary, we can arrange things so that also
\be\label{K2net}
d\big(y,\{y^n_1,\ldots,y^n_{m(n)}\}\big)\leq\eps\quad\forall y\in K_2.
\ee
Now using again (\ref{dminf}) we can add corresponding points in $K_1$ for the new points we have added to $K_2$ so that (\ref{und}) remains true. Adding points will not spoil (\ref{K1net}) so we can arrange things such that (\ref{K1net}), (\ref{dminf}), and (\ref{K2net}) are satisfied simultaneously.

Let $R_n\sub K_1\times K_2$ be the set
\be
R_n:=\big\{(x^n_k,y^n_k):1\leq k\leq m(n)\big\}.
\ee
We claim that $R_n$ is monotone in the sense that
\be\label{Rnmon}
\mbox{there are no $(x^n_k,y^n_k),(x^n_l,x^n_l)\in R_n$ such that $x^n_k\prec x^n_l$ and $y^n_l\prec y^n_k$.}
\ee
Indeed, $x^n_k\prec x^n_l$ implies $k<l$ and $y^n_l\prec y^n_k$ implies $l<k$, which is a contradiction.

Since $K_1\times K_2$ is compact, by Lemma~\ref{L:Haucomp}, we can select a subsequence such that $R_n\to R$ in the Hausdorff topology on $\Ki_+(K_1\times K_2)$, for some compact set $R\sub K_1\times K_2$. We claim that $R$ is a correspondence between $K_1$ and $K_2$. Indeed, by (\ref{K1net}), for each $x\in K_1$, we can choose $k(n)$ such that $x^n_{k(n)}\to x$. By the compactness of $K_1\times K_2$, the sequence $(x^n_{k(n)},y^n_{k(n)})$ has at least one cluster point $(x,y)$, and by Lemma~\ref{L:Hauconv} $(x,y)\in R$. Similarly, for each $y\in K_2$ there exists an $x\in K_1$ such that $(x,y)\in R$.

We next claim that $R$ is monotone. Assume that conversely, there exist $(x,y),(x',y')\in R$ such that $x\prec x'$ and $y'\prec y$. Then by Lemma~\ref{L:Hauconv}, there exist $k(n),k'(n)$ such that $(x^n_{k(n)},y^n_{k(n)})\to(x,y)$ and $(x^n_{k'(n)},y^n_{k'(n)})\to(x',y')$. By Lemma~\ref{L:evor}, $x^n_{k(n)}\prec y^n_{k(n)}$ and $y^n_{k'(n)}\prec x^n_{k'(n)}$ for all $n$ large enough, which contradicts (\ref{Rnmon}).

Taking the limit in (\ref{und}), using Lemma~\ref{L:Hauconv}, we see that
\be
d(x,y)\leq d^{\li\infty\re}(K_1,K_2)\quad\forall(x,y)\in R,
\ee
and hence by (\ref{tot}) $d_{\rm tot}(K_1,K_2)\leq d^{\li\infty\re}(K_1,K_2)$.
\epro

\subsection{The mismatch modulus}\label{S:mis}

In this subsection, we prove Theorem~\ref{T:partot}, except for the statement about Polishness. Generalising the definition in (\ref{mismatch}) for any $K_1,K_2\in\Ki_{\rm tot}(\Xc)$ and $\eps>0$, we define the \emph{mismatch modulus} $m_\eps(K_1,K_2)$ by
\be\ba{r@{\,}l}
\dis m_\eps(K_1,K_2):=\sup\big\{&\dis d(x_1,y_1)\vee d(x_2,y_2):
x_1,y_1\in K_1,\ x_2,y_2\in K_2,\\[5pt]
&\dis\quad d(x_1,x_2)\vee d(y_1,y_2)\leq\eps,\ x_1\pre y_1,\ y_2\pre x_2\big\}. 
\ec

\bl[Convergence of the mismatch modulus]
Let\label{L:misnul} $\Xc$ be a metrisable space. Assume that $K_n,K\in\Ki_{\rm part}(\Xc)$ satisfy $d^{\li 2\re}(K_n,K)\to 0$. Then
\be\label{misnul}
m_{\eps_n}(K_n,K)\asto{n}0\quad\mbox{with}\quad\eps_n:=d^{\li 1\re}(K_n,K).
\ee
\el

\bpro
If (\ref{misnul}) does not hold, then, by going to a subsequence, we can assume that there exists a $\de>0$ such that $m_{\eps_n}(K_n,K)\geq\de$ for all $n$. It follows that for each $n$, we can find $x_1(n),y_1(n)\in K_n$ and $x_2(n),y_2(n)\in K$ with $d(x_1(n),y_1(n))\vee d(x_2(n),y_2(n))\geq\de$ and $d(x_1(n),x_2(n))\vee d(y_1(n),y_2(n))\leq\eps_n$ such that $x_1(n)\pre y_1(n)$ and $y_2(n)\pre x_2(n)$. By Lemma~\ref{L:monesid}, our assumption $d^{\li 2\re}(K_n,K)\to 0$ implies $\eps_n=d^{\li 1\re}(K_n,K)\to 0$ and hence $K_n\to K$. Therefore, by Lemma~\ref{L:Hauconv}, there exists a compact set $C\sub\Xc$ such that $K_n\sub C$ for all $n$, so by going to a subsequence, we can assume that $x_1(n),x_2(n),y_1(n),y_2(n)$ converge to limits $x_1,x_2,y_1,y_2$ in $\Xc$. Since $d(x_1(n),x_2(n))\vee d(y_1(n),y_2(n))\leq\eps_n\to 0$, we have $x:=x_1=x_2$ and $y:=y_1=y_2$. On the other hand, our assumption that $d(x_1(n),y_1(n))\vee d(y_2(n),x_2(n))\geq\de$ implies that $d(x,y)\geq\de$. This leads to a contradiction, since by the assumption that $d^{\li 2\re}(K_n,K)\to 0$ and Lemma~\ref{L:orpres}, $x_1(n)\pre y_1(n)$ and $y_2(n)\pre x_2(n)$ imply $x\pre y$ and $y\pre x$ and hence $x=y$.
\epro

The following estimate essentially uses that the spaces are totally ordered.

\bl[Estimate in terms of mismatch modulus]
Let\label{L:disord} $\Xc$ be a metrisable space and let $K_1,K_2\in\Ki_{\rm tot}(\Xc)$ satisfy $d^{\li 1\re}(K_1,K_2)\leq\eps$. Then
\be\label{disord}
d^{\li m\re}(K_1,K_2)\leq m_\eps(K_1,K_2)+\eps\qquad(m\geq 1).
\ee
\el

\bpro
By symmetry, it suffices to show that for each $x_1=(x_1^1,\ldots,x_1^m)\in K_1^{\li m\re}$, there exists an $x_2=(x_2^1,\ldots,x_2^m)\in K_2^{\li m\re}$ such that $d^m(x_1,x_2)\leq m_\eps(K_1,K_2)+\eps$. In view of (\ref{dm}), the latter means that $d(x_1^k,x_2^k)\leq m_\eps(K_1,K_2)+\eps$ for all $1\leq k\leq m$. Since $d^{\li 1\re}(K_1,K_2)\leq\eps$, there exists a $z(1)\in K_2$ such that $d(x_1^1,z(1))\leq\eps$. We set $x_2^i=z(1)$ for all $1\leq i<I(1)$, where $I(1):=\inf\{i>1:d(x_1^i,z(1))>m_\eps(K_1,K_2)+\eps\}$. Using again that $d^{\li 1\re}(K_1,K_2)\leq\eps$, there exists a $z(2)\in K_2$ such that $d(x_1^{I(1)},z(2))\leq\eps$. Then
\be
d\big(z(1),z(2)\big)>d\big(x_1^{I(1)},z(2)\big)-\eps\geq m_\eps(K_1,K_2)
\ee
and hence by the definition of $m_\eps(K_1,K_2)$ and the fact that $x_1^1\pre x_1^{I(1)}$ we cannot have $z(2)\pre z(1)$. Since $K_2$ is totally ordered, we conclude that $z(1)\prec z(2)$. This allows us to set $x_2^i=z(2)$ for all $I(1)\leq i<I(2)$, where $I(2):=\inf\{i>I(1):d(x_1^i,z(2))>m_\eps(K_1,K_2)+\eps\}$. Continuing inductively, we obtain $(x_2^1,\ldots,x_2^m)\in K_2^{\li m\re}$ such that $d(x_1^k,x_2^k)\leq m_\eps(K_1,K_2)+\eps$ for all $1\leq k\leq m$.
\epro

As a consequence of Lemmas \ref{L:misnul} and \ref{L:disord}, we can prove that the metrics $d^{\li m\re}$ with $2\leq m\leq\infty$ all generate the same topology on $\Ki_{\rm tot}(\Xc)$. This may be a bit surprising in view of Lemmas \ref{L:noop} and \ref{L:diftop}. As the latter shows, the restriction to totally ordered sets is essential in the following lemma.

\bl[Convergence of totally ordered sets]
Let\label{L:totcoef} $\Xc$ be a metrisable space. Then $K_n,K\in\Ki_{\rm tot}(\Xc)$ satisfy $d^{\li 2\re}(K_n,K)\to 0$ if and only if $d^{\li\infty\re}(K_n,K)\to 0$.
\el

\bpro
Since $d^{\li m\re}(K_n,K)\leq d^{\li m+1\re}(K_n,K)$ by Lemma~\ref{L:monesid} and since $d^{\li\infty\re}(K_n,K)$ is defined as the increasing limit of $d^{\li m\re}(K_n,K)$ as $m\to\infty$, it is clear that $d^{\li\infty\re}(K_n,K)\to 0$ implies $d^{\li 2\re}(K_n,K)\to 0$. To prove the opposite implication, it suffices to show that $d^{\li 2\re}(K_n,K)\to 0$ implies\be
\sup_{m\geq 1}d^{\li m\re}(K_n,K)\asto{n}0.
\ee
By Lemma~\ref{L:monesid}, $d^{\li 2\re}(K_n,K)\to 0$ implies $\eps_n:=d^{\li 1\re}(K_n,K)\to 0$. Lemmas \ref{L:misnul} and \ref{L:disord} now imply that
\be
\sup_{m\geq 1}d^{\li m\re}(K_n,K)\leq m_{\eps_n}(K_n,K)+\eps_n\asto{n}0.
\ee
\epro

\bpro[of Theorem~\ref{T:partot} (in part)]
By the definition of $d_{\rm part}$ and Proposition~\ref{P:moncor} we have $d_{\rm part}=d^{\li 2\re}$ and $d_{\rm tot}=d^{\li\infty\re}$. By Lemma~\ref{L:totcoef} both metrics generate the same topology on $\Ki_{\rm tot}(\Xc)$. If $d$ and $d'$ generate the same topology on $\Xc$ and $d_{\rm part}$ and $d'_{\rm part}$ are defined in terms of $d$ and $d'$ as in (\ref{part}), then by Lemmas \ref{L:Hauprop} and \ref{L:Haucomp}, $d_{\rm part}$ and $d'_{\rm part}$ generate the same topology on $\Ki_{\rm tot}(\Xc)$. The inequalities (\ref{partot}) follow from the fact that $d^{\li 1\re}(K_1,K_2)\leq d^{\li 2\re}(K_1,K_2)\leq d^{\li\infty\re}(K_1,K_2)$ by Lemma~\ref{L:monesid}. On the other hand, Lemma~\ref{L:noop} shows that if $\Xc=[0,1]$, then for each $\eps>0$ we can find $K_1,K_2\in\Ki_{\rm tot}(\Xc)$ such that
\be
\dis d^{\li 2\re}(K_1,K_2)\leq\eps
\quad\mbox{while}\quad
\dis 1/2\leq d^{\li 3\re}(K_1,K_2)\leq d^{\li\infty\re}(K_1,K_2),
\ee
proving the final claim of the theorem. This proves all claims of the theorem except for the claim that $\Ki_{\rm tot}(\Xc)$ is Polish if $\Xc$ is Polish, which will be proved in the next subsection.
\epro

\subsection{Polishness}\label{S:Pol}

In this subsection, we complete the proof of Theorem~\ref{T:partot} by proving the following proposition.

\bp[Preservation of Polishness]
If\label{P:totPol} $\Xc$ is a Polish space, then so is $\Ki_{\rm tot}(\Xc)$, equipped with the ordered Hausdorff topology.
\ep

We start with the following lemma, announced in the introduction, that shows that even when $(\Xc,d)$ is complete, it is in general not true that the metrics $d_{\rm part}$ and $d_{\rm tot}$ are complete on $\Ki_{\rm tot}(\Xc)$.

\bl[Metric not complete]
Let\label{L:noncompl} $\Xc=[0,1]$, equipped with the usual distance. Then the metrics $d^{\li m\re}$ with $2\leq m\leq\infty$ are not complete on $\Ki_{\rm tot}(\Xc)$.
\el

\bpro
It suffices to construct a Cauchy sequence that does not converge. In view of Lemma~\ref{L:monesid}, it suffices to construct a Cauchy sequence in the metric $d^{\li\infty\re}$, which by Proposition~\ref{P:moncor} equals $d_{\rm tot}$. Let $\eps_n$ be positive constants converging to zero, and let $K_n:=\{0,1,\eps_n\}$ equipped with a total order such that $0\prec 1\prec\eps_n$. For each $n,m$, we define a monotone correspondence $R_{n,m}$ between $K_n$ and $K_m$ by $R_{n,m}:=\{(0,0),(1,1),(\eps_n,\eps_m)\}$. Then
\be
d_{\rm tot}(K_n,K_m)\leq\sup_{(x_1,x_2)\in R_{n,m}}|x_1-x_2|=|\eps_n-\eps_m|,
\ee
so $K_n$ is clearly a Cauchy sequence in $d_{\rm tot}$. However, the sequence $K_n$ does not converge in the ordered Hausdorff topology. If it had a limit $K$, then (in view of Lemma~\ref{L:orpres}) this would have to be the set $K=\{0,1\}$ equipped with a total order such that $0\pre 1$ and $1\pre 0$, but such a totally ordered set does not exist.
\epro

The proof of Proposition~\ref{P:totPol} needs some preparations. For each $L\in\Ki_+(\Xc^2)$ and $\eps\geq 0$, we set
\be
m^{\li 2\re}_\eps(L):=\sup\big\{d(x_1,y_1)\vee d(x_2,y_2):(x_1,y_1),(y_2,x_2)\in L,\ d(x_1,x_2)\vee d(y_1,y_2)\leq\eps\big\}.
\ee
In particular, this implies $m^{\li 2\re}_\eps(K^{\li 2\re})=m_\eps(K)$ $(K\in\Ki_{\rm tot}(\Xc))$.

\bl[Right continuity]
For\label{L:meps0} any metric space $\Xc$ and $L\in\Ki_+(\Xc^2)$, the function $\half\ni\eps\to m^{\li 2\re}_\eps(L)$ is right continuous.
\el

\bpro
The function $\eps\mapsto m^{\li 2\re}_\eps(L)$ is clearly nondecreasing, so it suffices to prove that
\be\label{meps0}
m^{\li 2\re}_\eps(L)\geq\lim_{\eta\down\eps}m^{\li 2\re}_\eta(L)\quad(\eps\geq 0)
\ee
where the limit exist by monotonicity. Fix $\eps_n>\eps$ such that $\eps_n\to\eps$. Then for each $\de>0$ and for each $n$, we can choose $(x^n_1,y^n_1),(y^n_2,x^n_2)\in L$ such that $d(x^n_1,x^n_2)\vee d(y^n_1,y^n_2)\leq\eps_n$ and $d(x^n_1,y^n_1)\vee d(x^n_2,y^n_2)\geq m^{\li 2\re}_{\eta_n}(L)-\de$. Since $L$ is compact, by going to a subsequence, we can assume that $(x^n_1,y^n_1)\to(x_1,y_1)$ and $(y^n_2,x^n_2)\to(y_2,x_2)$ for some $(x_1,y_1),(y_2,x_2)\in L$. Then $d(x_1,x_2)\vee d(y_1,y_2)\leq\eps$ and $d(x_1,y_1)\vee d(x_2,y_2)\geq\lim_{\eta\down\eps}m^{\li 2\re}_\eta(L)-\de$, which proves that $m^{\li 2\re}_\eps(L)\geq\lim_{\eta\down\eps}m^{\li 2\re}_\eta(L)-\de$. Since $\de>0$ is arbitrary, this implies (\ref{meps0}).
\epro

\bl[Upper semi-continuity]
Let\label{L:msemi} $\Xc$ be a metric space and let $L_n,L\in\Ki_+(\Xc^2)$ satisfy $L_n\to L$. Then
\be\label{msemi}
m^{\li 2\re}_\eps(L)\geq\limsup_{n\to\infty}m^{\li 2\re}_\eps(L_n)\quad(\eps\geq 0).
\ee
\el

\bpro
By the compactness of $[0,\infty]$ we can select a subsequence for which $\lim_{n\to\infty}m^{\li 2\re}_\eps(L_n)$ exists and is equal to the limit superior of the original sequence. Let $\de_n>0$ converge to zero and pick $(x^n_1,y^n_1),(y^n_2,x^n_2)\in L_n$ such that $d(x^n_1,x^n_2)\vee d(y^n_1,y^n_2)\leq\eps$ and $d(x^n_1,y^n_1)\vee d(x^n_2,y^n_2)\geq m^{\li 2\re}_\eps(L_n)-\de_n$. By Lemma~\ref{L:Hauconv}, there exists a compact $C\sub\Xc^2$ such that $L_n\sub C$ for all $n$, so by going to a further subsequence we can assume that $(x^n_1,y^n_1)\to(x_1,y_1)$ and $(y^n_2,x^n_2)\to(y_2,x_2)$ for some $(x_1,y_1),(y_2,x_2)\in\Xc^2$. Then $(x_1,y_1),(y_2,x_2)\in L$ by Lemma~\ref{L:Hauconv}, $d(x_1,x_2)\vee d(y_1,y_2)\leq\eps$, and hence
\be
m^{\li 2\re}_\eps(L)\geq d(x_1,y_1)\vee d(x_2,y_2)\geq\lim_{n\to\infty}\big(m^{\li 2\re}_\eps(L_n)-\de_n).
\ee
Since $\de_n\to 0$, this proves (\ref{msemi}).
\epro

Before we can prove Proposition~\ref{P:totPol} we need one more lemma. For any metric space $(\Xc,d)$, we define $\Li(\Xc)\sub\Ki_+(\Xc^2)$ by
\be\label{Lidef}
\Li(\Xc):=\big\{K^{\li 2\re}:K\in\Ki_{\rm tot}(\Xc)\big\},
\ee
and we let $\ov{\Li(\Xc)}$ denote the closure of $\Li(\Xc)$ in the metric space $\big(\Ki_+(\Xc^2),d^2_{\rm H}\big)$.

\bl[Totally ordered sets]
For\label{L:Lident} any metric space $\Xc$, one has
\be\label{Lident}
\Li(\Xc)=\big\{L\in\ov{\Li(\Xc)}:m^{\li 2\re}_0(L)=0\big\}.
\ee
\el

\bpro
To prove the inclusion $\sub$ in (\ref{Lident}), it suffices to observe that
\be
m^{\li 2\re}_0(K^{\li 2\re})=\sup\big\{d(x,y):(x,y),(y,x)\in K^{\li 2\re}\big\}=0
\ee
for all $K\in\Ki_{\rm tot}(\Xc)$, since $x\pre y$ and $y\pre x$ imply $x=y$.

We next prove the inclusion $\supset$ in (\ref{Lident}). Assume that $L\in\ov{\Li(\Xc)}$ satisfies $m^{\li 2\re}_0(L)=0$. Since $L\in\ov{\Li(\Xc)}$, there exist $K_n\in\Ki_{\rm tot}(\Xc)$ such that $K_n^{\li 2\re}\to L$ in the topology on $\Ki_+(\Xc^2)$. Let $\pi_i(x_1,x_2):=x_i$ $(i=1,2)$ denote the coordinate projections $\pi_i\cn\Xc^2\to\Xc$. Since $\pi_1(K^{\li 2\re}_n)=K_n=\pi_2(K^{\li 2\re}_n)$ for each $n$, using Lemma~\ref{L:contim}, we see that $K_n\to K$ in the Hausdorff topology on $\Ki_+(\Xc)$, where $K:=\pi_1(L)=\pi_2(L)$. We define a relation $\pre$ on $K$ by setting $x\pre y$ if and only if $(x,y)\in L$. To complete the proof, it suffices to show that $\pre$ is a total order on $K$, i.e.,
\begin{enumerate}
\item for each $x,y\in K$, either $x\pre y$ or $y\pre x$, or both,
\item $x\pre y$ and $y\pre x$ imply $x=y$,
\item $x\pre y\pre z$ imply $x\pre z$.
\end{enumerate}
To prove (i), let $x,y\in K$. Since $K_n\to K$ in the Hausdorff topology on $\Ki_+(\Xc)$, by Lemma~\ref{L:Hauconv}, there exist $x_n,y_n\in K_n$ such that $x_n\to x$ and $y_n\to y$. Since $K_n$ is totally ordered, either $x_n\pre y_n$ happens for infinitely many $n$, or $y_n\pre x_n$ happens for infinitely many $n$, or both. Since $K_n^{\li 2\re}\to L$ in the Hausdorff topology on $\Ki_+(\Xc^2)$, by Lemma~\ref{L:Hauconv}, it folows that either $x\pre y$ or $y\pre x$, or both. Property~(ii) follows immediately from the fact that $\sup\{d(x,y):(x,y),(y,x)\in L\}=m^{\li 2\re}_0(L)=0$. To prove (iii), assume that $x,y,z\in K$ satisfy $x\pre y\pre z$. If $x=y$ or $y=z$ then trivially $x\pre z$, so without loss of generality we may assume that $x\neq y\neq z$. Since $K_n\to K$ in the Hausdorff topology on $\Ki_+(\Xc)$, by Lemma~\ref{L:Hauconv}, there exist $x_n,y_n,z_n\in K_n$ such that $x_n\to x$, $y_n\to y$, and $z_n\to z$. Since $K_n$ is totally ordered, for each $n$ either $x_n\pre y_n$, or $y_n\pre x_n$, or both. But $y_n\pre x_n$ can happen only for finitely many $n$ since otherwise the fact that $K_n^{\li 2\re}\to L$ in the Hausdorff topology on $\Ki_+(\Xc^2)$ and Lemma~\ref{L:Hauconv} would imply that $y\pre x$, which together with our assumptions $x\pre y$ and $x\neq y$ contradicts (ii). We conclude that $x_n\pre y_n$ for all $n$ sufficiently large and by the same argument also $y_n\pre z_n$ for all $n$ sufficiently large. It follows that $x_n\pre z_n$ for all $n$ sufficiently large. Since $K_n^{\li 2\re}\to L$ in the Hausdorff topology on $\Ki_+(\Xc^2)$, Lemma~\ref{L:Hauconv} shows that $(x,z)\in L$ and hence $x\pre z$.
\epro

We are now ready to prove Proposition~\ref{P:totPol}. We need to recall one well-known fact. A subset $A\sub\Xc$ of a topological space $\Xc$ is called a \emph{$G_\de$-set} if $A$ is a countable intersection of open sets. Our proof of Proposition~\ref{P:totPol} makes use of the following fact, that we cite from \cite[\S 6 No.~1, Thm.~1]{Bou58} (see also \cite[Thms~12.1 and 12.3]{Oxt80}).

\bl[Subsets of Polish spaces]
A\label{L:subPol} subset $\Yi\sub\Xc$ of a Polish space $\Xc$ is Polish in the induced topology if and only if $\Yi$ is a $G_\de$-subset of $\Xc$.
\el

\bpro[of Proposition~\ref{P:totPol}]
We first observe that if $\Xc$ is Polish, then so is $\Xc^2$, equipped with the product topology, and hence, by Lemma~\ref{L:Hauprop}, also $\Ki_+(\Xc^2)$. Therefore, since $\big(\Ki_{\rm tot}(\Xc),d_{\rm part}\big)$ is isometric to $\big(\Li(\Xc),d^2_{\rm H}\big)$ defined in (\ref{Lidef}), in view of Lemma~\ref{L:subPol}, it suffices to show that $\Li(\Xc)$ is a $G_\de$-subset of $\Ki_+(\Xc^2)$.

It follows from Lemma~\ref{L:msemi} that for each $\eps,\de>0$, the set
\be
\Ai_{\de,\eps}:=\big\{L\in\Ki_+(\Xc^2):m^{\li 2\re}_\eps(L)\geq\de\big\}
\ee
is a closed subset of $\Ki_+(\Xc^2)$ and hence its complement $\Ai_{\eps,\de}^{\rm c}$ is open. As a consequence,
\be
\Gi:=\bigcap_{n=1}^\infty\bigcup_{m=1}^\infty\Ai_{1/n,1/m}^{\rm c}
\ee
is a $G_\de$-set. Since each closed set is a $G_\de$-set, and the intersection of two $G_\de$-sets is a $G_\de$-set, using Lemmas \ref{L:meps0} and \ref{L:Lident}, we conclude that
\be\ba{l}
\dis\ov{\Li(\Xc)}\cap\Gi
=\big\{L\in\ov{\Li(\Xc)}:\forall\de>0\ \exists\eps>0\mbox{ s.t. }m^{\li 2\re}_\eps(L)<\de\big\}\\[5pt]
\dis\quad=\big\{L\in\ov{\Li(\Xc)}:\lim_{\eps\to 0}m^{\li 2\re}_\eps(L)=0\big\}
=\big\{L\in\ov{\Li(\Xc)}:m^{\li 2\re}_0(L)=0\big\}
=\Li(\Xc)
\ec
is a $G_\de$-set.
\epro

\subsection{Compactness criterion}\label{S:totcomp}

In this subsection, we prove Theorem~\ref{T:totcomp}.\med

\bpro[of Theorem~\ref{T:totcomp}]
Since the map $K\mapsto K^{\li 2\re}$ is a homeomorphism from $\Ki_{\rm tot}(\Xc)$ to $\Li(\Xc)$, equipped with the induced topology from $\Ki_+(\Xc^2)$, a set $\Ai\sub K^{\li 2\re}$ is precompact if and only if $\Bi:=\{K^{\li 2\re}:K\in\Ai\}$ is a precompact subset of $\Ki_+(\Xc^2)$ and its closure $\ov\Bi$ is contained in $\Li(\Xc)$. We will show that $\Bi$ is a precompact subset of $\Ki_+(\Xc^2)$ if and only if (\ref{totcomp}~(i) holds. Moreover, if (\ref{totcomp}~(i) holds, then $\ov\Bi$ is contained in $\Li(\Xc)$ if and only if (\ref{totcomp}~(ii) holds.

If (\ref{totcomp}~(i) holds, then $C^2$ is a compact subset of $\Xc^2$ and $K^{\li 2\re}\sub C^2$ for all $K\in\Ai$, so $\Bi$ is a precompact subset of $\Ki_+(\Xc^2)$ by Lemma~\ref{L:Haucomp}. Conversely, if $\Bi$ is a precompact subset of $\Ki_+(\Xc^2)$, then by Lemma~\ref{L:Haucomp} there exists a compact subset $D\sub\Xc^2$ such that $K^{\li 2\re}\sub D$ for all $K\in\Ai$. Without loss of generality, we may assume that $D$ is of the form $D=C^2$ for some compact subset $C$ of $\Xc$; for example, we may take for $C$ the union of the two coordinate projections of $D$. Then $K\sub C$ for all $K\in\Ai$, proving that (\ref{totcomp}~(i) holds.

To complete the proof, assume that (\ref{totcomp}~(i) holds. We must show that $\ov\Bi$ is contained in $\Li(\Xc)$ if and only if (\ref{totcomp}~(ii) holds. Assume, first, that (\ref{totcomp}~(ii) does not hold. Then we can find a $\de>0$ and $\eps_n>0$ tending to zero, as well as $K_n\in\Ai$, such that $m^{\li 2\re}_{\eps_n}(K_n)\geq\de$ for each $n$. By (\ref{totcomp}~(i), going to a subsequence if necessary, we can assume that $K_n^{\li 2\re}\to L$ for some $L\in\Ki_+(\Xc^2)$. Now Lemma~\ref{L:msemi} implies that
\be
m^{\li 2\re}_\eps(L)\geq\limsup_{n\to\infty}m^{\li 2\re}_\eps(K_n^{\li 2\re})\geq\limsup_{n\to\infty}m^{\li 2\re}_{\eps_n}(K_n^{\li 2\re})\geq\de
\ee
for each $\eps>0$, so letting $\eps\down 0$, using Lemma~\ref{L:meps0}, we conclude that $m^{\li 2\re}_0(L)\geq\de$. By Lemma~\ref{L:Lident}, we conclude that $L\not\in\Li(X)$ and hence $\ov\Bi$ is not contained in $\Li(\Xc)$.

Assume now that (\ref{totcomp})~(ii) holds. We must show that $\ov\Bi$ is contained in $\Li(\Xc)$. Assume that $K_n^{\li 2\re}\to L$ for some $K_n\in\Ai$. Then clearly $L\in\ov{\Li(\Xc)}$ so by Lemma~\ref{L:Lident} it suffices to prove that $m^{\li 2\re}_0(L)=0$. Assume that, conversely, there exist $x,y\in\Xc$ with $x\neq y$ such that $(x,y),(y,x)\in L$. Then by Lemma~\ref{L:Hauconv}, there exist $x^n_1,x^n_2,y^n_1,y^n_2\in K_n$ with $x^n_1\pre y^n_1$ and $y^n_2\pre x^n_2$ such that $x^n_i\to x$ and $y^n_i\to y$ as $n\to\infty$ $(i=1,2)$. Then for each $\eps>0$, we can choose $n$ large enough such that $d(x^n_i,x)\leq\eps/2$ $(i=1,2)$. It follows that $d(x^n_1,y^n_1)\vee d(x^n_2,y^n_2)\geq d(x,y)-\eps$ and $d(x^n_1,x^n_2)\vee d(y^n_1,y^n_2)\geq\eps$ so that $m_\eps(K_n)\geq d(x,y)-\eps$. This clearly contradicts (\ref{totcomp})~(ii), so we conclude that $m^{\li 2\re}_0(L)=0$ as required.
\epro

\subsection{Cadlag curves}\label{S:curve}

In this subsection, we prove Proposition~\ref{P:curve}. Let $(\Xc,d)$ be a metric space. If $R$ is any subset of $\Xc^2$, then let us call
\be
{\rm dist}(R):=\sup_{(x_1,x_2)\in R}d_{\rm sqz}(x_1,x_2)
\ee
the \emph{distortion} of $R$. Then (\ref{Haus2}) and (\ref{tot}) say that
\be\label{totdis}
d_{\rm H}(K_1,K_2)=\inf_{R\in{\rm Corr}(K_1,K_2)}{\rm dist}(R)
\quand
d_{\rm tot}(K_1,K_2)=\inf_{R\in{\rm Corr}_+(K_1,K_2)}{\rm dist}(R),
\ee
where ${\rm Corr}(K_1,K_2)$ and ${\rm Corr}_+(K_1,K_2)$ denote the sets of all (monotone) correspondences between $K_1$ and $K_2$. Let $\ov R$ denote the closure of a set $R\sub\Xc^2$. Then
\be\label{closR}
{\rm dist}(R)={\rm dist}(\ov R)\qquad(R\sub\Xc^2).
\ee
Indeed, the inequality $\leq$ is trivial, while the opposite inequality follows from the fact that for each $(x_1,x_2)\in\ov R$, there exist $(x^n_1,x^n_2)\in R$ such that $(x^n_1,x^n_2)\to(x_1,x_2)$ and hence $d(x^n_1,x^n_2)\to d(x_1,x_2)$.

We need one preparatory lemma.

\bl[Fine partition]
Let\label{L:fine} $(\Xc,d)$ be a metric space. Then for each $\ga\in\Di_{[0,1]}(\Xc)$ and $\eps>0$, there exist $t_0<0<t_1<\cdots<t_{n-1}<1<t_n$ such that
\be\label{fine}
\sup\big\{d\big(\ga(s),\ga(s')\big):1\leq k\leq n,\ s,t\in[0,1]\cap[t_{k-1},t_k)\big\}<\eps.
\ee
\el

\bpro
By Lemma~\ref{L:cadlag}, writing $\ga(t+):=\ga(t)$, we can view $\ga$ as a continuous function on the split real interval $\lc 0-,1+\rc$, that moreover satisfies $\ga(0-)=\ga(0+)$ and $\ga(1-)=\ga(1+)$. Fix $\eps>0$ and let
\be
R:=\big\{(s,t)\in\R^2:s<t,\ s\neq 0,\ t\neq 1,\ d\big(\ga(\sig),\ga(\tau)\big)<\eps\ \forall\sig,\tau\in\lc s+,t-\rc\cap\lc 0-,1+\rc\big\}.
\ee
Using the properties of $\ga$, it is easy to see that
\be
\bigcup_{(s,t)\in R}\lc s+,t-\rc\supset\lc 0-,1+\rc.
\ee
Since the intervals $\lc s+,t-\rc$ are open in the topology of the split real line and since $\lc 0-,1+\rc$ is compact by Proposition~\ref{P:Rpmcompact}, there exists a finite subcover, i.e., there exists a finite set $S\sub R$ such that
\be
\bigcup_{(s,t)\in S}\lc s+,t-\rc\supset\lc 0-,1+\rc.
\ee
Let $T:=\{s:(s,t)\in S\}\cup\{t:(s,t)\in S\}$. Then, letting $t_0$ denote the largest element of $T\cap(-\infty,0)$, ordering the elements of $T\cap(0,1)$ as $t_1<\cdots<t_{n-1}$, and letting $t_n$ denote the smallest element of $T\cap(1,\infty)$, we obtain times $t_0<\cdots<t_n$ as in (\ref{fine}).
\epro

\bpro[of Proposition~\ref{P:curve}]
If $\ga_1,\ga_2$ are cadlag parametrisations of $K_1,K_2$, and $\la\in\La$, then let us set
\be
R_\la:=\big\{\big(\ga_1(t),\ga_2\big(\la(t)\big)\big):t\in[0,1]\big\}=\big\{\big(\ga_1\big(\la^{-1}(t)\big),\ga_2(t)\big):t\in[0,1]\big\},
\ee
and let $\ov R_\la$ denote its closure. We claim that $\ov R_\la$ is a correspondence between $K_1$ and $K_2$. Indeed, by the definition of a cadlag parametrisation, each element $x_1\in K_1$ is of the form $x_1=\ga_1(t)$ or $=\ga_1(t-)$ for some $t\in[0,1]$. If $x_1=\ga_1(t)$, then clearly there exists an $x_2\in K_2$ such that $(x_1,x_2)\in R_\la$, namely $x_2:=\ga_2(\la(t))$. If $x_1=\ga_1(t-)$, then we can choose $t_n\up t$ and set $x^n_1:=\ga_1(t_n)$. Then $x^n_1\to x_1$ by the left continuity of $t\mapsto\ga_1(t-)$. We have already seen that there exist $x^n_2\in K_2$ such that $(x^n_1,x^n_2)\in R_\la$. Since $K_2$ is compact, by going to a subsequence, we can assume that $x^n_2\to x_2$ for some $x_2\in K_2$. Then $(x_1,x_2)\in\ov R_\la$. In the same way, we see that for each $x_2\in K_2$, there exists an $x_1\in K_1$ such that $(x_1,x_2)\in\ov R_\la$. This completes the proof that $\ov R_\la$ is a correspondence. Using the fact that the total orders on $K_1$ and $K_2$ are compatible with the topology, it is easy to see that $\ov R_\la$ is monotone if $\la\in\La_+$.

Using these facts as well as (\ref{totdis}) and (\ref{closR}), we see that
\bc\label{curcorleq}
\dis d_{\rm H}(K_1,K_2)&\leq&\dis\inf_{\la\in\La}{\rm dist}(\ov R_\la)=\inf_{\la\in\La}{\rm dist}(R_\la)=\inf_{\la\in\La}\sup_{t\in[0,1]}d\big(\ga_1(t),\ga_2\big(\la(t)\big)\big),\\[5pt]
\dis d_{\rm tot}(K_1,K_2)&\leq&\dis\inf_{\la\in\La_+}\!{\rm dist}(\ov R_\la)=\inf_{\la\in\La_+}\!{\rm dist}(R_\la)=\inf_{\la\in\La_+}\sup_{t\in[0,1]}d\big(\ga_1(t),\ga_2\big(\la(t)\big)\big).
\ec
To complete the proof, we must show that:
\begin{enumerate}
\item For each $R\in{\rm Corr}(K_1,K_2)$ and $\eps>0$, there exists a $\la\in\La$ such that ${\rm dist}(R_\la)\leq{\rm dist}(R)+\eps$.
\item For each $R\in{\rm Corr}_+(K_1,K_2)$ and $\eps>0$, there exists a $\la\in\La_+$ such that ${\rm dist}(R_\la)\leq{\rm dist}(R)+\eps$.
\end{enumerate}
We first prove (i). Fix $R\in{\rm Corr}(K_1,K_2)$ and $\eps>0$. By Lemma~\ref{L:fine}, for $i=1,2$, there exist $t^i_0<0<t^i_1<\cdots<t^i_{n_i-1}<1<t^i_{n_i}$ such that
\be
\sup\big\{d\big(\ga_i(s),\ga_i(s')\big):1\leq k\leq n_i,\ s,t\in[0,1]\cap[t^i_{k-1},t^i_k)\big\}<\eps/2\quad(i=1,2).
\ee
For $1\leq k\leq n_i$, let us write $K^i_k:=\{\ga_i(t):t\in[0,1]\cap[t^i_{k-1},t^i_k)\}$ $(i=1,2)$. We can choose a correspondence $S$ between $\{1,\ldots,n_1\}$ and $\{1,\ldots,n_2\}$ such that for each $(k_1,k_2)\in S$, there exists an $(x_1,x_2)\in R$ with $x_i\in K^i_{k_i}$ $(i=1,2)$. Then
\be\label{Kpart}
\sup_{(k_1,k_2)\in S}\sup\subb{x_1\in K_1}{x_2\in K_2}d(x_1,x_2)\leq{\rm dist}(R)+\eps.
\ee
By refining the partitions $t^i_0,\ldots,t^i_{n_i}$, we can make sure that for each $k_1\in\{1,\ldots,n_1\}$, there is a unique $k_2\in\{1,\ldots,n_2\}$ such that $(k_1,k_2)\in S$, and vice versa. We can then construct a bijection $\la\cn[0,1]\to[0,1]$ such that for each $(k_1,k_2)\in S$, the restriction of $\la$ to $[0,1]\cap[t^1_{k_1-1},t^1_{k_1})$ is a bijection to $[0,1]\cap[t^2_{k_2-1},t^2_{k_2})$. Then (\ref{Kpart}) implies that ${\rm dist}(R_\la)\leq{\rm dist}(R)+\eps$, completing the proof of (i).

To also prove (ii), we observe that if $R$ is a monotone correspondence, then $S$ as we initially constructed it will be a monotone correspondence between $\{1,\ldots,n_1\}$ and $\{1,\ldots,n_2\}$, and monotonicity will be preserved after we refine the partitions so that they have the same size and $S$ is a bijection. Now $\la$ can be chosen monotone too, completing the proof of (ii).
\epro

\section{Proofs of the main results}\label{S:mainproof}

\subsection{Metrics on path space}\label{S:Pimet}

In this section we prove Theorems \ref{T:J1} and \ref{T:M1}. We fix a metrisable space $\Xc$ on which there is defined a betweenness that is compatible with the topology. Recall that $\ov\pi$ denotes the filled path associated with a path $\pi\in\Pi(\Xc)$. By Lemmas \ref{L:graph}, \ref{L:intchar}, and \ref{L:fillgra},
\bc\label{ovPi}
\dis\ov\Pi(\Xc)&:=&\dis\big\{\ov\pi:\pi\in\Pi(\Xc)\big\}\\[5pt]
&=&\dis\big\{K\in\Ki_{\rm tot}(\Ri(\Xc)):K\mbox{ satisfies conditions (i) and (ii) of Lemma~\ref{L:intchar}}\big\}.
\ec
By Lemma~\ref{L:intchar}, a path $\pi$ is uniquely characterised by its filled graph $\ov\pi$, so given a metric $d_{\rm sqz}$ on the squeezed space $\Ri(\Xc)$, we can define metrics $d_{\rm tot}$ and $d_{\rm part}$ on $\Pi(\Xc)$ by
\be
d_{\rm tot}(\pi_1,\pi_2):=d_{\rm tot}(\ov\pi_1,\ov\pi_2)
\quand
d_{\rm part}(\pi_1,\pi_2):=d_{\rm part}(\ov\pi_1,\ov\pi_2),
\ee
where the metrics on the right-hand side of each equation are those on $\Ki_{\rm tot}(\Ri(\Xc))$ defined in (\ref{part}) and (\ref{tot}), and those on the left-hand side of each equation coincide with those on $\Pi(\Xc)$ already defined in (\ref{dSkor}). Similarly, $d_{\rm H}(\pi_1,\pi_2)$ defined in (\ref{Hpi2}) is the Hausdorff distance between $\ov\pi_1$ and $\ov\pi_2$, which is in general only a pseudometric. Since the total order $\pre$ on a continuous path $\pi\in\Pi_{\rm c}(\Xc)$ is uniquely determined by the compact set $\pi\sub\Ri(\Xc)$, however, $d_{\rm H}$ is a metric on $\Pi_{\rm c}(\Xc)$. For the same reason, if $I\sub\R$ is a closed interval of positive length, then $d_{\rm H}$ is a metric on $\Di_I(\Xc)$, where we identify the latter with a subset of $\Pi(\Xc)$ as in (\ref{DiI}). The following proposition summarises what we already know.\med

\bp[Basic properties of the metrics]
Assume\label{P:M1} that $\Xc$ is a metrisable space that is equipped with a proper betweenness. Let $d_{\rm H}$, $d_{\rm tot}$, and $d_{\rm part}$ be defined as in (\ref{Hpi2}) and (\ref{dSkor}). Then
\be\label{dcomp3}
d_{\rm H}(\pi_1,\pi_2)\leq d_{\rm part}(\pi_1,\pi_2)\leq d_{\rm tot}(\pi_1,\pi_2)
\qquad\big(\pi_1,\pi_2\in\Pi(\Xc)\big).
\ee
The functions $d_{\rm tot}$ and $d_{\rm part}$ are metrics on $\Pi(\Xc)$ that generate the same topology. This topology does not depend on the choice of the metric $d_{\rm sqz}$ on the squeezed space $\Ri(\Xc)$.
\ep

\bpro
Immediate from Theorem~\ref{T:partot}.
\epro

The following proposition implies that on $\Pi_{\rm c}(\Xc)$, all three metrics $d_{\rm H}$, $d_{\rm tot}$, and $d_{\rm part}$ generate the same topology.

\bp[Space of continuous paths]
Assume\label{P:Pic} that $\Xc$ is a metrisable space that is equipped with a proper betweenness. Let $d_{\rm H}$, $d_{\rm tot}$, and $d_{\rm part}$ be defined as in (\ref{Hpi2}) and (\ref{dSkor}). Then for paths $\pi_n\in\Pi(\Xc)$ and $\pi\in\Pi_{\rm c}(\Xc)$, the following statements are equivalent:
\[
{\rm(i)}\ d_{\rm part}(\pi_n,\pi)\asto{n}0,\quad
{\rm(ii)}\ d_{\rm tot}(\pi_n,\pi)\asto{n}0,\quad
{\rm(iii)}\ d_{\rm H}(\pi_n,\pi)\asto{n}0.
\]
\ep

\bpro
By Theorem~\ref{T:partot}, $d_{\rm part}$ and $d_{\rm tot}$ generate the same topology on $\Pi(\Xc)$, and by (\ref{dcomp3}) convergence in any of these two metrics implies convergence in $d_{\rm H}$. Therefore, to show that the conditions (i)--(iii) are equivalent, it suffices to show that for $\pi_n\in\Pi(\Xc)$ and $\pi\in\Pi_{\rm c}(\Xc)$, condition (iii) implies (i). More precisely, we will show that for the filled paths defined in (\ref{fillpath}), $\ov\pi_n\to\ov\pi$ in the Hausdorff topology implies
\be
\ov\pi_n^{\li 2\re}\asto{n}\ov\pi^{\li 2\re}
\ee
in the Hausdorff topology. By Lemma~\ref{L:Hauconv}, convergence of $\ov\pi_n$ implies the existence of a compact set $C\sub\Ri(\Xc)$ such that $\ov\pi_n\sub C$ for all $n$, which implies $\ov\pi_n^{\li 2\re}\sub C^2$. To complete the proof, by Lemma~\ref{L:Hauconv}, we need to prove the following two statements.
\begin{enumerate}
\item For every $\big((x,s),(y,t)\big)\in\ov\pi^{\li 2\re}$, there exist $\big((x_n,s_n),(y_n,t_n)\big)\in\ov\pi_n^{\li 2\re}$ such that\\ $\big((x_n,s_n),(y_n,t_n)\big)\to\big((x,s),(y,t)\big)$.
\item If a sequence $\big((x_n,s_n),(y_n,t_n)\big)\in\ov\pi_n^{\li 2\re}$ has a cluster point $\big((x,s),(y,t)\big)\in\Ri(\Xc)^2$,\\ then $\big((x,s),(y,t)\big)\in\ov\pi^{\li 2\re}$.
\end{enumerate}
To prove (i), we use the fact that $\ov\pi_n\to\ov\pi$ to find $(x_n,s_n),(y_n,t_n)\in\ov\pi_n$ such that $(x_n,s_n)\to(x,s)$ and $(y_n,t_n)\to(y,t)$. If $s<t$, then $\big((x_n,s_n),(y_n,t_n)\big)\in\ov\pi_n^{\li 2\re}$ for $n$ large enough, so (i) follows. On the other hand, if $s=t$, then $\big((x_n,s_n),(x_n,s_n)\big)\in\ov\pi_n^{\li 2\re}$ so (i) also holds in this case.

To prove (ii), we use the fact that $\ov\pi_n\to\ov\pi$ to conclude that any cluster point $\big((x,s),(y,t)\big)$ satisfies $(x,s),(y,t)\in\ov\pi$ with $s\leq t$, and hence by the continuity of $\pi$ either $s<t$ or $(x,s)=(y,t)$, from which we conclude that $\big((x,s),(y,t)\big)\in\ov\pi^{\li 2\re}$.
\epro

\bl[Connected paths]
Let\label{L:connectp} $\Xc$ be a metrisable space that is equipped with the trivial betweenness. Then $\Pi^|_{\rm c}(\Xc)$ defined in (\ref{Pint}) is a closed subset of $\Pi(\Xc)$.
\el

\bpro
Immediate from Lemmas \ref{L:Hconnec} and \ref{L:congraph}.
\epro

\subsection{Polishness}\label{S:JPol}

In this subsection we prove that if $\Xc$ is Polish and equipped with a betweenness that is compatible with the topology, then so are $\Pi(\Xc)$, $\Pi_{\rm c}(\Xc)$, and $\Di_I(\Xc)$ equipped with the Skorokhod topology corresponding to the betweennness on $\Xc$. Recall the definition of the modulus of continuity $m_{T,\de}(\pi)$ in (\ref{modul}). By a slight abuse of notation, for any $K\in\Ki_+(\Ri(\Xc))$, we set
\be\label{moduli1}
m_{T,\de}(K):=\sup\big\{d(x_1,x_2):(x_i,t_i)\in K,\ -T\leq t_1\leq t_2\leq T,\ t_2-t_1\leq\de\big\}.
\ee

\bl[Upper semi-continuity]
Let\label{L:upper1} $\Xc$ be a metrisable space and assume that $K_n,K\in\Ki_+(\Ri(\Xc))$ satisfy $K_n\to K$. Then, for each $T<\infty$ and $\de>0$,
\be\label{upper1}
m_{T,\de}(K)\geq\limsup_{n\to\infty}m_{T,\de}(K_n).
\ee
\el

\bpro
By the compactness of $[0,\infty]$ we can select a subsequence for which $\lim_{n\to\infty}m_{T,\de}(K_n)$ exists and is equal to the limit superior of the original sequence. Let $\eps_n>0$ converge to zero and pick $(x^n_i,t^n_i)\in K_n$ with $-T\leq t^n_1\leq t^n_2\leq T$ and $t^n_2-t^n_1\leq\de$, such that $d(x^n_1,x^n_2)\geq m_{T,\de}(K_n)-\eps_n$. Since $K_n\to K$, by Lemma~\ref{L:Hauconv} there exists a compact $C\sub\Ri(\Xc)$ such that $K_n\sub C$ for all $n$. It follows that we can select a subsequence such that $(x^n_i,t^n_i)\to(x_i,t_i)$ for some $(x_i,t_i)\in K$ $(i=1,2)$. Then $-T\leq t_1\leq t_2\leq T$ and $t_2-t_1\leq\de$, so
\be
m_{T,\de}(K)\geq d(x_1,x_2)\geq\lim_{n\to\infty}\big(m_{T,\de}(K_n)-\eps_n\big)=\lim_{n\to\infty}m_{T,\de}(K_n).
\ee
Since we have chosen our subsequence such that the right-hand side is equal to the limit superior of the original sequence, this proves the claim.
\epro

\bp[The space of continuous paths is Polish]
If\label{P:PcPol} $\Xc$ is a Polish space, then so is $\Pi_{\rm c}(\Xc)$, equipped with the topology generated by $d_{\rm H}$.
\ep

\bpro
By Lemma~\ref{L:conchar} (using also Lemma~\ref{L:pathgra}), we may identify $\Pi_{\rm c}(\Xc)$ with the set of all $G\in\Ki_+(\Ri(\Xc))$ such that
\be\label{maxone}
\big|\big\{x\in\Xc:(x,t)\in G\big\}\big|\leq 1\quad\forall t\in\R.
\ee
We claim that this condition is equivalent to
\be\label{onelim}
\lim_{\de\to 0}m_{T,\de}(G)=0\ \forall T<\infty.
\ee
If for some $t\in[-T,T]$ there exist $x_1\neq x_2$ such that $(x_1,t),(x_2,t)\in G$, then $m_{T,\de}(G)\geq d(x_1,x_2)$ for all $\de>0$ so (\ref{onelim}) is violated. Conversely, if (\ref{onelim}) does not hold then for some $\eps>0$ we can choose $(x^n_i,t^n_i)\in G$ with $t^n_i\in[-T,T]$ $(i=1,2)$ such that $d(x^n_1,x^n_2)\geq\eps$ for all $n$. Taking the limit, using the compactness of $G$, we see that (\ref{maxone}) is violated.

By Lemma~\ref{L:upper1}, for each $T<\infty$ and $\eps,\de>0$, the set
\be
\Gi_{T,\eps,\de}:=\big\{G\in\Ki_+(\Ri(\Xc)):m_{T,\de}(G)\geq\eps\big\}
\ee
is a closed subset of $\Ki_+(\Xc^2)$ and hence its complement $\Gi_{T,\eps,\de}^{\rm c}$ is open. As a consequence, using (\ref{onelim}), we see that
\be
\Pi_{\rm c}(\Xc)=\bigcap_{N=1}^\infty\bigcap_{n=1}^\infty\bigcup_{m=1}^\infty\Gi_{N,1/n,1/m}^{\rm c}
\ee
is a $G_\de$-subset of $\Ki_+(\Ri(\Xc))$ and hence Polish in the induced topology by Lemma~\ref{L:subPol}.
\epro

Recall the definition of the Skorokhod modulus of continuity $m^{\rm S}_{T,\de}(\pi)$ in (\ref{mJM1}). By a slight abuse of notation, for any $K\in\Ki_{\rm tot}(\Ri(\Xc))$, we set
\be\ba{l@{\,}l}\label{moduli}
\dis m^{\rm S}_{T,\de}(K):=\sup\big\{d\big(x_2,\li x_1,x_3\re\big):&\dis(x_i,t_i)\in K\ \;\mbox{and}\!\ -T\leq t_i\leq T\ \forall i=1,2,3,\\[5pt]
&\dis(x_1,t_1)\pre(x_2,t_2)\pre(x_3,t_3),\ t_3-t_1\leq\de\big\}.
\ec
With this definitions, if $\ov\pi$ is the filled path associated with a path $\pi\in\Pi(\Xc)$, then $m^{\rm S}_{T,\de}(\ov\pi)$ as defined in (\ref{moduli}) corresponds to $m^{\rm S}_{T,\de}(\pi)$ as defined (\ref{mJM1}). The following lemma is similar to Lemmas \ref{L:msemi} and \ref{L:upper1}.

\bl[Upper semi-continuity]
Let\label{L:upper} $\Xc$ be a metrisable space that is equipped with a proper betweenness, and assume that $K_n,K\in\Ki_{\rm tot}(\Ri(\Xc))$ satisfy $K_n\to K$. Then, for each $T<\infty$ and $\de>0$,
\be\label{upper}
m^{\rm S}_{T,\de}(K)\geq\limsup_{n\to\infty}m^{\rm S}_{T,\de}(K_n).
\ee
\el

\bpro
The proof will be very similar to the proof of Lemma~\ref{L:msemi}. By the compactness of $[0,\infty]$ we can select a subsequence for which $\lim_{n\to\infty}m^{\rm S}_{T,\de}(K_n)$ exists and is equal to the limit superior of the original sequence. Let $\eps_n>0$ converge to zero and pick $(x^n_i,t^n_i)\in K_n$ with $-T\leq t^n_i\leq T$ $(i=1,2,3)$, $(x^n_1,t^n_1)\pre(x^n_2,t^n_2)\pre(x^n_3,t^n_3)$, and $t^n_3-t^n_1\leq\de$, such that $d\big(x^n_2,\li x^n_1,x^n_3\re\big)\geq m^{\rm S}_{T,\de}(K_n)-\eps_n$. Since $K_n\to K$ in the topology on $\Ki_{\rm tot}(\Ri(\Xc))$, by the first inequality in (\ref{partot}) they also converge in the topology on $\Ki_+(\Ri(\Xc))$, so by Lemma~\ref{L:Hauconv} there exists a compact $C\sub\Ri(\Xc)$ such that $K_n\sub C$ for all $n$. It follows that we can select a subsequence such that $(x^n_i,t^n_i)\to(x_i,t_i)$ for some $(x_i,t_i)\in K$ $(i=1,2,3)$. Recall from Proposition~\ref{P:moncor} that $d_{\rm tot}=d^{\li\infty\re}$ so by Lemma~\ref{L:monesid}, convergence in $\Ki_{\rm tot}(\Ri(\Xc))$ implies that $K^{\li m\re}_n\to K^{\li m\re}$ in the Hausdorff topology for any $1\leq m\leq\infty$. In particular, $K^{\li 3\re}_n\to K^{\li 3\re}$, which by Lemma~\ref{L:Hauconv} implies $(x_1,t_1)\pre(x_2,t_2)\pre(x_3,t_3)$. Moreover $-T\leq t_i\leq T$ $(i=1,2,3)$ and $t_3-t_1\leq\de$, so
\be
m^{\rm S}_{T,\de}(K)\geq d\big(x_2,\li x_1,x_3\re\big)\geq\lim_{n\to\infty}\big(m^{\rm S}_{T,\de}(K_n)-\eps_n\big)=\lim_{n\to\infty}m^{\rm S}_{T,\de}(K_n),
\ee
where we have used that $d\big(x^n_2,\li x^n_1,x^n_3\re\big)\to d\big(x_2,\li x_1,x_3\re\big)$ since the betweenness is compatible with the topology. Since we have chosen our subsequence such that the right-hand side is equal to the limit superior of the original sequence, this proves the claim.
\epro

\bp[Skorokhod topologies are Polish]
Let\label{P:J1M1Pol} $\Xc$ be a metrisable space equipped with a betweennesss that is compatible with the topology. Then if $\Xc$ is Polish, so is $\Pi(\Xc)$, equipped with the Skorokhod topology corresponding to the betweenness.
\ep

\bpro
We observe that if $\Xc$ is Polish, then, by Lemma~\ref{L:pointprop}, so is $\Ri(\Xc)$ and hence, by Proposition~\ref{P:totPol}, also $\Ki_{\rm tot}(\Ri(\Xc))$. By identifying a path with its filled path, we can as in (\ref{ovPi}) identity $\Pi(\Xc)$ with a subset of $\Ki_{\rm tot}(\Ri(\Xc))$. The Skorokhod topology on $\Pi(\Xc)$ is then the induced topology from $\Ki_{\rm tot}(\Ri(\Xc))$. Therefore, in view of Lemma~\ref{L:subPol}, it suffices to show that the set $\ov\Pi(\Xc)$ defined in (\ref{ovPi}) is a $G_\de$-subset of $\Ki_{\rm tot}(\Ri(\Xc))$.

We start by showing that condition~(ii) of Lemma~\ref{L:intchar} can be replaced by
\begin{itemize}
\item[{\rm(ii)'}] $\dis\lim_{\de\to 0}m^{\rm S}_{T,\de}(G)=0\ \forall T<\infty$.
\end{itemize}
To see this, we argue as follows. If (ii) does not hold, then for some $t\in\R$, there exist $(x_i,t_i)\in G$ $(i=1,2,3)$ with $(x_1,t_1)\pre(x_2,t_2)\pre(x_3,t_3)$ and $x_2\not\in\li x_1,x_3\re$, which implies that $m^{\rm S}_{T,\de}(G)\geq d\big(x_2,\li x_1,x_3\re\big)>0$ for all $\de>0$ and $T<\infty$ such that $-T\leq t\leq T$, so (ii)' clearly does not hold. Conversely, if (ii)' does not hold, then for some $T<\infty$ and $\eps>0$ we can choose $\de_n>0$ tending to zero and $(x^n_i,t^n_i)\in G$ $(i=1,2,3)$ with $(x^n_1,t^n_1)\pre(x^n_2,t^n_2)\pre(x^n_3,t^n_3)$ such that $d\big(x^n_2,\li x^n_1,x^n_3\re\big)\geq\eps$. By the compactness of $G$, we can select a subsequence such that $(x^n_i,t^n_i)\to(x_i,t_i)$ $(i=1,2,3)$. Then clearly $t_1=t_2=t_3=:t$ for some $-T\leq t\leq T$. Since the order is compatible with the topology moreover $(x_1,t_1)\pre(x_2,t_2)\pre(x_3,t_3)$. The fact that the betweenness is compatible with the topology allows us to conclude that $d\big(x_2,\li x_1,x_3\re\big)=\lim_{n\to\infty}d\big(x^n_2,\li x^n_1,x^n_3\re\big)\geq\eps$. This shows that $x_2\not\in\li x_1,x_3\re$ and hence (ii) does not hold.

Let $\Hi$ denote the set of all elements of $\Ki_{\rm tot}(\Ri(\Xc))$ that satisfy condition (i) of Lemma~\ref{L:intchar}. By what we have just proved,
\be\label{Piid}
\ov\Pi(\Xc)=\big\{G\in\Hi:\lim_{\de\to 0}m^{\rm S}_{T,\de}(G)=0\ \forall T<\infty\big\}.
\ee
It follows from Lemma~\ref{L:upper} that for each $T<\infty$ and $\eps,\de>0$, the set
\be
\Gi_{T,\eps,\de}:=\big\{G\in\Ki_{\rm tot}(\Ri(\Xc)):m^{\rm S}_{T,\de}(G)\geq\eps\big\}
\ee
is a closed subset of $\Ki_{\rm tot}(\Ri(\Xc))$ and hence its complement $\Gi_{T,\eps,\de}^{\rm c}$ is open. As a consequence,
\be
\Gi:=\bigcap_{N=1}^\infty\bigcap_{n=1}^\infty\bigcup_{m=1}^\infty\Gi_{N,1/n,1/m}^{\rm c}
\ee
is a $G_\de$-set. Formula (\ref{Piid}) says that $\ov\Pi(\Xc)=\Gi\cap\Hi$. It is easy to see that $\Hi$ is a closed subset of $\Ki_{\rm tot}(\Ri(\Xc))$, and hence a $G_\de$-set. Since the intersection of two $G_\de$-sets is a $G_\de$-set, this yields the statement we wanted to prove.
\epro

\bp[Polishness of the space of functions on a fixed interval]
Let\label{P:DIPol} $\Xc$ be a metrisable space, equipped with a betweenness that is compatible with the topology, let $I\sub\R$ be a closed real interval of positive length, and let $\Di_I(\Xc)$ be the space in (\ref{DiI}), equipped with the Skorokhod topology corresponding to the betweenness on $\Xc$. If $\Xc$ is Polish, then so is $\Di_I(\Xc)$.
\ep

\bpro
It follows from Lemma~\ref{L:contim} that the set
\be
\ov\Di_I(\Xc):=\big\{\pi\in\Pi(\Xc):\hat I(\pi)=I\cup\{-\infty,\infty\}\big\}
\ee
is a closed subset of $\Pi(\Xc)$, and hence Polish by Lemma~\ref{L:subPol} and Proposition~\ref{P:J1M1Pol}. We observe that for each $t\in I$, the set
\be\ba{l}
\dis\big\{\pi\in\ov\Di_I(\Xc):\pi(t-)=\pi(t+)\big\}\\[5pt]
\dis\quad=\bigcap_{n=1}^\infty\bigcup_{m=1}^\infty\big\{\pi\in\ov\Di_I(\Xc):d\big(\pi(s),\pi(t)\big)<\ffrac{1}{n}\ \forall s\in[t-\ffrac{1}{m},t+\ffrac{1}{m}]\cap I\big\}
\ec
is a $G_\de$-subset of $\ov\Di_I(\Xc)$, and hence by (\ref{DiI}), so is $\Di_I(\Xc)$. By Lemma~\ref{L:subPol}, it follows that $\Di_I(\Xc)$ is Polish.
\epro

\subsection{Compactness criteria}\label{S:cocrit}

In this subsection we prove Theorems \ref{T:ArzAsc}, \ref{T:ArzAsc2}, and \ref{T:ArzI}, as well as Lemma~\ref{L:J1vsM1}.\med

\bpro[of Theorem~\ref{T:ArzAsc}]
As in the proof of Proposition~\ref{P:PcPol} we identify $\Pi_{\rm c}(\Xc)$ with the subset of $\Ki_+(\Ri(\Xc))$ consisting of all $G$ that satisfy (\ref{maxone}) or equivalently (\ref{onelim}). Let $\Ai\sub\Pi_{\rm c}(\Xc)$. Then $\Ai$ is precompact if and only if its closure $\ov\Ai$ in $\Ki_+(\Ri(\Xc))$ is compact and $\ov\Ai\sub\Pi_{\rm c}(\Xc)$. By Lemmas  \ref{L:REcomp} and \ref{L:Haucomp}, $\ov\Ai$ is a compact subset of $\Ki_+(\Ri(\Xc))$ if and only if $\Ai$ satisfies the compact containment condition. Therefore, to complete the proof, it suffices to show that if $\Ai\sub\Pi_{\rm c}(\Xc)$ satisfies the compact containment condition, then $\ov\Ai\sub\Pi_{\rm c}(\Xc)$ if and only if $\Ai$ is equicontinuous.

Assume that $\Ai$ satisfies the compact containment condition and is equicontinuous, and that $G_n\in\Ai$ satisfy $G_n\to G$ for some $G\in\Ki_+(\Ri(\Xc))$. We need to show that $G\in\Pi_{\rm c}(\Xc)$. If this is not the case, then there exist $(x_1,t),(x_2,t)\in G$ with $t\in\R$ and $x_1\neq x_2$. Then by Lemma~\ref{L:Hauconv}, there exist $(x^n_1,t^n_1),(x^n_2,t^n_2)\in G_n$ such that $(x^n_i,t^n_i)\to(x_i,t)$ as $n\to\infty$ $(i=1,2)$. Choose $T<\infty$ such that $-T<t<T$. Then for each $\de>0$ we have for $n$ large enough that $-T<t^n_1,t^n_2<T$, $|t^n_1-t^n_2|\leq\de$, and $d(x^n_1,x^n_2)\geq d(x_1,x_2)/2$. This proves that
\be
\sup_{\pi\in\Ai}m_{T,\de}(\pi)\geq d(x_1,x_2)/2\qquad\forall\de>0,
\ee
contradicting the equicontinuity of $\Ai$.

Assume, on the other hand, that $\Ai$ satisfies the compact containment condition and is not equicontinuous. Let $\de_n$ be positive constants tending to zero. Since $\Ai$ is not equicontinuous, for some $T<\infty$ and $\eps>0$ we can find $G_n\in\Ai$ such that $m_{T,\de_n}(G_n)\geq\eps$ for all $n$. Since $\Ai$ satisfies the compact containment condition, $\ov\Ai$ is a compact subset of $\Ki_+(\Ri(\Xc))$, so by going to a subsequence we may assume that $G_n\to G$ for some $G\in\Ki_+(\Ri(\Xc))$. Since $m_{T,\de_n}(G_n)\geq\eps$ we can find $(x^n_1,t^n_1),(x^n_2,t^n_2)\in G_n$ such that $-T\leq t^n_1<t^n_2\leq T$, $t^n_2-t^n_1\leq\de_n$, and $d(x^n_1,x^n_2)\geq\eps$. By Lemma~\ref{L:Hauconv}, going to a further subsequence if necessary, we can assume that $(x^n_i,t^n_i)\to(x_i,t)$ as $n\to\infty$ $(i=1,2)$ for some $(x_1,t),(x_2,t)\in G$. Then $-T\leq t\leq T$ and $d(x_1,x_2)\geq\eps$, which shows that $\{x\in\Xc:(x,t)\in G\}$ has more than one element and hence $\ov\Ai$ is not contained in $\Pi_{\rm c}(\Xc)$.
\epro

\bpro[of Theorem~\ref{T:ArzAsc2}]
We identify $\Pi(\Xc)$ with $\ov\Pi(\Xc)$ which as in (\ref{ovPi}) we view as a subset of $\Ki_{\rm tot}(\Ri(\Xc))$. In this identification, the metrics $d_{\rm part}$ and $d_{\rm tot}$ on $\Ki_{\rm tot}(\Ri(\Xc))$ induce metrics $d_{\rm part}$ and $d_{\rm tot}$ on $\Pi(\Xc)$ that both generate the Skorokhod topology.

Let $\Ai\sub\Pi(\Xc)$. Then $\Ai$ is precompact in the Skorokhod topology if and only if its closure $\ov\Ai$ in $\Ki_{\rm tot}(\Ri(\Xc))$ is compact and $\ov\Ai\sub\Pi(\Xc)$. By Theorem~\ref{T:totcomp} and Lemma~\ref{L:REcomp}, $\ov\Ai$ is a compact subset of $\Ki_{\rm tot}(\Ri(\Xc))$ if and only if $\Ai$ satisfies the compact containment condition and
\be\label{mism}
\lim_{\eps\to 0}\sup_{G\in\Ai}m_\eps(G)=0,
\ee
where $m_\eps(G)$ denotes the mismatch modulus of $G$. To complete the proof we will prove the following three statements.
\begin{itemize}
\item[{\rm I}] Assume that $\Ai$ is satisfies the compact containment condition and is Skorokhod-equi\-con\-tinuous. Then $\Ai$ satisfies (\ref{mism}).
\item[{\rm II}] Assume that $\Ai$ is Skorokhod-equicontinuous. Then $\ov\Ai\sub\Pi(\Xc)$.
\item[{\rm III}] Assume that $\ov\Ai$ is a compact subset of $\Ki_{\rm tot}(\Ri(\Xc))$ and that $\ov\Ai\sub\Pi(\Xc)$. Then $\Ai$ is Skorokhod-equicontinuous.
\end{itemize}
Now if $\Ai$ satisfies the compact containment condition and is Skorokhod-equicontinuous, then by our earlier remarks I implies that $\ov\Ai$ is a compact subset of $\Ki_{\rm tot}(\Ri(\Xc))$ and II implies that $\ov\Ai\sub\Pi(\Xc)$, so $\Ai$ is precompact in the Skorokhod topology. Conversely, if $\Ai$ is precompact in the Skorokhod topology, then $\ov\Ai$ is a compact subset of $\Ki_{\rm tot}(\Ri(\Xc))$ and hence by our earlier remarks $\Ai$ satisfies the compact containment condition, and moreover $\ov\Ai\sub\Pi(\Xc)$ which by III implies that $\Ai$ is Skorokhod-equicontinuous.

We start by proving I. Since $\Ai$ satisfies the compact containment condition, by Lemma~\ref{L:REcomp}, we see that there exists a compact $C\sub\Ri(\Xc)$ such that $G\sub C$ for all $G\in\Ai$. Since $\sup_{G\in\Ai}m_\eps(G)$ is nondecreasing as a function of $\eps$, the limit in (\ref{mism}) always exists. Let $\eps_n$ be positive constants, tending to zero. If the limit in (\ref{mism}) is positive, then there exists a $\de>0$ such that for each $n$ we can find a $G\in\Ai$ and $(x^n_i,s^n_i),(y^n_i,t^n_i)\in G$ $(i=1,2)$ with $(x^n_1,s^n_1)\pre(y^n_1,t^n_1)$, $(y^n_2,t^n_2)\pre(x^n_2,s^n_2)$ such that
\[\ba{c}
\dis d_{\rm sqz}\big((x^n_1,s^n_1),(x^n_2,s^n_2)\big)\vee d_{\rm sqz}\big((y^n_1,t^n_1),(y^n_2,t^n_2)\big)\leq\eps_n,\\[5pt]
\dis d_{\rm sqz}\big((x^n_i,s^n_i),(y^n_i,t^n_i)\big)\geq\de\quad(i=1,2).
\ea\]
Since $G\sub C$ for all $G\in\Ai$, by going to a subsequence, we may assume that $(x^n_i,s^n_i)\to(x,s)$ and $(y^n_i,t^n_i)\to(y,t)$ $(i=1,2)$ for some $(x,s),(y,t)\in\Ri(\Xc)$. Then $d_{\rm sqz}\big((x,s),(y,t)\big)\geq\de$ and hence $(x,s)\neq(y,t)$. Since $(x^n_1,s^n_1)\pre(y^n_1,t^n_1)$ and $(y^n_2,t^n_2)\pre(x^n_2,s^n_2)$ we have $s^n_1\leq t^n_1$ and $t^n_2\leq s^n_2$ for all $n$ which implies $s=t$ and hence $x\neq y$, since $(x,t)\neq(y,t)$. By the structure of $\Ri(\Xc)$, this implies $t\in\R$. Let $(x^n_-,s^n_-)$ (resp.\ $(x^n_+,s^n_+)$) be the smallest (resp.\ largest) of the points $(x^n_i,s^n_i)$ $(i=1,2)$ with respect to the order $\pre$, and define $(y^n_\pm,t^n_\pm)$ similarly. Since $G$ is totally ordered, by going to a subsequence, we can assume that we are in one of the following two cases. 1.\ $(x^n_-,s^n_-)\pre(y^n_-,t^n_-)$ for all $n$, or 2.\ $(y^n_-,t^n_-)\pre(x^n_-,s^n_-)$ for all $n$. Let us assume that we are in case 1. Then $(x^n_-,s^n_-)\pre(y^n_-,t^n_-)\pre(x^n_+,s^n_+)$ for all $n$. Since the betweenness is compatible with the topology,
\be
d\big(y^n_-,\li x^n_-,x^n_+\re\big)\asto{n}d(y,x)>0,
\ee
which contradicts the Skorokhod-equicontinuity. Case~2 is completely the same, exchanging the roles of $x$ and~$y$.

We next prove II. Assume that $\pi_n\in\Ai$ converge in $\Ki_{\rm tot}(\Ri(\Xc))$ to a limit $G$. Recall from Subsection~\ref{S:orH} that $d_{\rm part}=d^{\li 2\re}\leq d^{\li m\re}\leq d^{\li\infty\re}=d_{\rm tot}$ for all $m\geq 2$. In particular, $\pi_n\to G$ in $\Ki_{\rm tot}(\Xc)$ implies that $\pi^{\li m\re}_n\to G^{\li m\re}$ in the Hausdorff topology for all $m\geq 2$. It suffices to check that $G$ satisfies conditions (i) and (ii) of Lemma~\ref{L:intchar}. Condition~(i) easily follows from the fact that $\pi^{\li 2\re}_n\to G^{\li 2\re}$ in the Hausdorff topology, using Lemma~\ref{L:Hauconv}. It remains to prove that $G$ satisfies condition~(ii) of Lemma~\ref{L:intchar}. Assume that conversely, for some $t\in\R$, there exist $(x_1,t),(x_2,t),(x_3,t)\in G$ with $(x_1,t)\pre(x_2,t)\pre(x_3,t)$ such that $x_2\not\in\li x_1,x_3\re$. Since $\pi^{\li 3\re}_n\to G^{\li 3\re}$ in the Hausdorff topology, by Lemma~\ref{L:Hauconv} there exist $\tau^n_i\in I_\sg(\pi_n)\cup\{\pm\binf\}$ $(i=1,2,3)$ with $\tau^n_1\leq\tau^n_2\leq\tau^n_3$ such that $\un\tau^n_i\to t$ and $\pi_n(\tau^n_i)\to x_i$ $(i=1,2,3)$. Using the fact that the betweenness is compatible with the topology, we see that
\be
d\big(\pi_n(\tau^n_2),\li\pi_n(\tau^n_1),\pi_n(\tau^n_3)\re\big)
\asto{n}d\big(x_2,\li x_1,x_3\re\big),
\ee
which is easily seen to contradict Skorokhod-equicontinuity, completing the proof of~II.

To prove III, finallly, we will show that if $\ov\Ai$ is a compact subset of $\Ki_{\rm tot}(\Ri(\Xc))$ and $\Ai$ is not Skorokhod-equicontinuous, then $\ov\Ai$ is not contained in $\Pi(\Xc)$. Let $\de_n$ be positive constants, tending to zero. Since $\Ai$ is not Skorokhod-equicontinuous, there exists a $\eps>0$, $T<\infty$, $\pi_n\in\Ai$, and $\tau^n_i\in I_\sg(\pi_n)$ $(i=1,2,3)$ such that $\tau_1\leq\tau_2\leq\tau_3$, $-T\leq\un\tau^n_1$, $\un\tau^n_3\leq T$, $\un\tau^n_3-\un\tau^n_1\leq\de_n$, and $d\big(\pi(\tau^n_2),\li\pi(\tau^n_1),\pi(\tau^n_3)\re\big)\geq\eps$. Since $\ov\Ai$ is a compact subset of $\Ki_{\rm tot}(\Ri(\Xc))$, we can select a subsequence such that $\pi_n\to G$ for some $G\in\Ki_{\rm tot}(\Ri(\Xc))$. Then $\pi^{\li 3\re}_n\to G^{\li 3\re}$ in the Hausdorff topology, so by Lemma~\ref{L:Hauconv} there exist $C\sub\Ri(\Xc)^3$ such that $\pi^{\li 3\re}_n\sub C$ for all $n$. It follows that by going to a further subsequence we can assume that $\big(\pi_n(\tau^n_i),\un\tau^n_i\big)\to(x_i,t_i)$ as $n\to\infty$ for some $(x_i,t)\in G$ $(i=1,2,3)$ with $(x_1,t_1)\pre(x_2,t_2)\pre(x_3,t_3)$ and $-T\leq t\leq T$. Using the fact that the betweenness is compatible with the topology, we see that $d(x_2,\li x_1,x_3\re)\geq\eps$, which shows that $G$ is not the filled graph of a path $\pi\in\Pi(\Xc)$ and hence $\ov\Ai$ is not contained in $\Pi(\Xc)$.
\epro

As fairly simple applications of the results proved in this subsection we can also prove Lemmas \ref{L:J1vsM1} and \ref{L:interpol1}.\med

\bpro[of Lemma~\ref{L:J1vsM1}]
We start by noting that for each $x,z\in\Xc$, the total order on $\li x,z\re$ coincides with the induced order from its embedding in $\li x,z\re'$. Indeed, $y,y'\in\li x,z\re$ satisfy $y\leq_{x,z}y'$ if and only if $y'\in\li y,z\re$, which implies $y'\in\li y,z\re'$ and hence $y\leq'_{x,z}y'$, where $\leq'_{x,z}$ denotes the total order on $\li x,z\re'$.

For any path $\pi$, let $\ov\pi$ and $\ov\pi'$ denote the filled paths defined in terms of the betweennesses $\li\,\cdot\,,\,\cdot\,\re$ and $\li\,\cdot\,,\,\cdot\,\re'$, respectively. The condition $\li x,z\re\subset\li x,z\re'$ $(x,z\in\Xc)$ implies that $\ov\pi\sub\ov\pi'$. By our earlier remarks, the total order on $\ov\pi$ coincides with the induced order from $\ov\pi'$. Using notation as in (\ref{Kli}), it follows that $\ov\pi^{\li 2\re}\sub{\ov\pi'}^{\li 2\re}$ for all $\pi\in\Pi(\Xc)$.

Let ${\rm S}$ and ${\rm S'}$ denote the Skorokhod topologies defined by the betweennesses $\li\,\cdot\,,\,\cdot\,\re$ and $\li\,\cdot\,,\,\cdot\,\re'$. By Theorem~\ref{T:partot}, one has $\pi_n\to\pi$ in the topology ${\rm S}$ if and only if $\ov\pi^{\li 2\re}_n\to\ov\pi^{\li 2\re}$ in the Hausdorff topology on $\Ki_+(\Ri(\Xc)^2)$, and an analogous statement holds for for convergence in the topology ${\rm S'}$. Let $m^{\rm S}_{T,\de}$ and $m^{\rm S'}_{T,\de}$ denote the Skorokhod moduli of continuity associated with both betweennesses. The condition $\li x,z\re\subset\li x,z\re'$ $(x,z\in\Xc)$ implies $d(y,\li x,z\re)\geq d(y,\li x,z\re')$ $(x,y,z\in\Xc)$ and hence $m^{\rm S}_{T,\de}(\pi)\geq m^{\rm S'}_{T,\de}(\pi)$ for all $\pi\in\Pi(\Xc)$, $T<\infty$, and $\de>0$.

Assume that $\pi_n,\pi\in\Pi(\Xc)$ satisfy $\pi_n\to\pi$ in the topology ${\rm S}$. Then by Theorem~\ref{T:ArzAsc2}, the set $\{\pi_n:n\in\N\}$ satisfies the compact containment condition and is Skorokhod equicontinuous with respect to $m^{\rm S}_{T,\de}$ and hence also with respect to $m^{\rm S'}_{T,\de}$, which by Theorem~\ref{T:ArzAsc2} implies that the set $\{\pi_n:n\in\N\}$ is precompact in the topology ${\rm S'}$. Therefore, to complete the proof, it suffices to show that if $\pi'$ is a cluster point of the sequence $\pi_n$ in the topology ${\rm S'}$, then $\pi'=\pi$.

By going to a subsequence, we may assume that $\pi_n\to\pi$ in the topology ${\rm S}$ and $\pi_n\to\pi'$ in the topology ${\rm S'}$. By our earlier remarks, $\ov\pi^{\li 2\re}_n\to\ov\pi^{\li 2\re}$ in the Hausdorff topology and similarly ${\ov\pi'_n}^{\!\li 2\re}\to{\ov\pi'}^{\li 2\re}$. It follows that $I(\pi')=I(\pi)$ and $\ov\pi^{\li 2\re}\sub{\ov\pi'}^{\li 2\re}$, which implies $\pi'=\pi$, as required.
\epro

\bpro[of Lemma~\ref{L:interpol1}]
If $\pi'_n\to\pi$, then by Theorem~\ref{T:ArzAsc} the $\pi'_n$ are equicontinuous and satisfy the compact containment condition. Since $\pi_n\sub\pi'_n$ the same is then true for the $\pi_n$ so by Theorem~\ref{T:ArzAsc}, going to a subsequence if necessary, we can assume that $\pi_n\to\ti\pi$ for some $\ti\pi\sub\pi$. Since for each $t\in I(\pi)$, there exist $t_n\in I(\pi_n)$ such that $t_n\to t$, each such subsequential limit must satisfy $\ti\pi=\pi$, which allows us to conclude that $\pi_n\to\pi$.

Assume, conversely, that $\pi_n\to\pi$. We claim that for each $z_n=(x_n,t_n)\in\pi'_n$, we can select a subsequence such that $z_n\to z\in\pi$. To see this, let $s_n:=\sup\{t\in I(\pi_n):t\leq t_n\}$ and $u_n:=\inf\{t\in I(\pi_n):t\geq t_n\}$. By going to a subsequence, we can assume that $s_n,t_n$, and $u_n$ converge to limits $s,t,u\in\ov\R$ with $s\leq t\leq u$. Since $\pi\in\Pi^|_{\rm c}(\Xc)$ and $\pi_n\to\pi$ we must have $s=u$. Since $\pi_n\to\pi$ both $(\pi_n(s_n),s_n)$ and $(\pi_n(u_n),u_n)$ converge to $(\pi(t),t)$. If $t=\pm\infty$, then $z_n\to(\ast,\pm\infty)=(\pi(t),t)$ and we are done. If $t\in\R$ then we use that $\pi'_n(t_n)\in\li\pi_n(s_n),\pi_n(u_n)\re$ and the betweenness is compatible with the topology to conclude $(\pi'_n(t_n),t_n)\to(\pi(t),t)$.

By the claim we have just proved, there exists a compact set $C\sub\Ri(\Xc)$ such that $\pi'_n\sub C$ for each $n$. By Lemma~\ref{L:Haucomp}, it follows that $\{\pi'_n:n\in\N\}$ is precompact in the topology on $\Ki_+(\Ri(\Xc))$. By the claim we have just proved, each subsequential limit $\pi^\ast$ satisfies $\pi^\ast\sub\pi$. Since $\pi_n\sub\pi'_n$ and $\pi_n\to\pi$ we conclude $\pi^\ast=\pi$. It follows that $\pi'_n\to\pi$.
\epro

\subsection{Identification of the J1, J2, M1, and M2 topologies}\label{S:classic}

In this subsection we complete the proofs of Theorems \ref{T:J1} and \ref{T:M1} by showing that for the trivial betweenness, the topology we have introduced on the path space $\Pi(\Xc)$ induces on the space $\Di_I(\Xc)$ from (\ref{DiI}) Skorokhod's J1 topology while the metric $d_{\rm H}$ generates the J2 topology, and similarly, if $\Xc=\R$ equipped with the linear betweenness, then we obtain the M1 and M2 topologies. In the process also establish Lemma~\ref{L:halfline}.

We view a function $f\in\Di_{[0,t]}(\Xc)$ as an extended cadlag function $(\hat I,f)$ with $\hat I=[0,t]\cup\{-\infty,\infty\}$, and similarly we view a function $f\in\Di_\half(\Xc)$ as an extended cadlag function $(\hat I,f)$ with $\hat I=[0,\infty]\cup\{-\infty\}$. We let $\Gi_{\rm f}(f):=\Gi_{\rm f}(\hat I,f)$ denote its filled graph, defined in (\ref{Gint}). For $f\in\Di_\half(\Xc)$ we let $f\big|_{[0,t]}$ denote the restriction of $f$ to the interval $[0,t]$.

\bl[Convergence of restricted functions]
Let\label{L:restr} $(\Xc,d)$ be a metric space that is equipped with a proper betweenness and let $d_{\rm H}$ and $d_{\rm tot}$ be defined as in (\ref{Hpi2}) and (\ref{dSkor}). Then for all $f_n,f\in\Di_\half(\Xc)$, one has
\be\ba{r@{\,}c@{\,}ll}\label{atco}
\dis d_{\rm H}(f_n,f)\asto{n}0\quad&\desd&\dis\quad d_{\rm H}\big(f_n\big|_{[0,t]},f\big|_{[0,t]}\big)\asto{n}0\quad&\dis\forall t>0\mbox{ s.t.\ }f(t-)=f(t),\\[5pt]
\dis d_{\rm tot}(f_n,f)\asto{n}0\quad&\desd&\dis\quad d_{\rm tot}\big(f_n\big|_{[0,t]},f\big|_{[0,t]}\big)\asto{n}0\quad&\dis\forall t>0\mbox{ s.t.\ }f(t-)=f(t).
\ec
\el

\bpro
We start by proving the statement for $d_{\rm H}$. Let $T:=\{t>0:f(t-)=f(t)\}$. It is well-known that the complement of $T$ is countable (see \cite[Lemma~3.5.1]{EK86}), so $T$ is dense in $\half$. Let $G:=\Gi_{\rm f}(f)$, $G_n:=\Gi_{\rm f}(f_n)$, $G^t:=\Gi_{\rm f}(f\big|_{[0,t]})$, and $G^t_n:=\Gi_{\rm f}(f_n\big|_{[0,t]})$.

We first prove the implication $\volgt$. Let $t\in T$. Since $G_n\to G$, by Lemma~\ref{L:Hauconv}, there exists a compact set $C$ such that $G_n\sub C$ for all $n$, and hence also $G^t_n\sub C$ for all $n$. By Lemma~\ref{L:Haucomp}, it follows that $\{G^t_n:n\in\N\}$ is compact, so to prove that $G^t_n\to G^t$, it suffices to show that $G^t$ is the only subsequential limit of the $G^t_n$. Let $G^\ast$ be such a subsequential limit. Since $G^t_n\sub G_n$ it is clear from Lemma~\ref{L:Hauconv} that $G^\ast\sub G$. Let $\psi$ denote the projection $\psi(x,s):=s$ and let $\psi(G^\ast)$ denote the image of $G^\ast$ under $\psi$. By Lemma~\ref{L:contim}, $\psi(G^\ast)=[0,t]$. It follows that $G^\ast\sub G^t$. To prove the opposite inclusion, assume that $(y,s)\in G^t$. If $s<t$, then we use that by Lemma~\ref{L:Hauconv}, there exist $(x_n,s_n)\in G_n$ such that $(x_n,s_n)\to(x,s)$. Since $s<t$, we have $(x_n,s_n)\in G^t_n$ for $n$ large enough and hence $(x,s)\in G^\ast$. If $s=t$, then we use that $\psi(G^\ast)=[0,t]$ to conclude that there must be at least one $y'\in\Xc$ such that $(y',t)\in G^\ast$. Since $G^\ast\sub G^t$ we must have $y,y'\in\li f(t-),f(t)\re=\{f(t)\}$, where we have used that $t\in T$, so we conclude that $y'=y$, concluding the proof that $G^\ast=G^t$.

We next prove the implication $\Leftarrow$. Since $G^t_n\to G^t$ for each $t\in T$, using Lemmas \ref{L:Hauconv} and \ref{L:REcomp}, we see that there exists a compact set $C$ such that $G_n\sub C$ for all $n$, so by Lemma~\ref{L:Haucomp} it suffices to show that if $G_\ast$ is a subsequential limit of the $G_n$, then $G_\ast=G$. Since $G^t_n\sub G_n$ for each $n$ it is clear from Lemma~\ref{L:Hauconv} that $G^t\sub G_\ast$ for each $t\in T$. We claim that conversely, for each $(x,s)\in G_\ast$ and $s<t\in T$, we have $(x,s)\in G^t$. Indeed, by Lemma~\ref{L:Hauconv} for some subsequence there exist $(x_n,s_n)\in G_n$ such that $(x_n,s_n)\to(x,s)$. Since $s_n<t$ for $n$ large enough, it follows that $(x,s)\in G^t$. These arguments show that $\{(x,t)\in G_\ast:t<\infty\}=\{(x,t)\in G:t<\infty\}$, which is enough to conclude $G_\ast=G$.

We next prove the statement for $d_{\rm tot}$. Let $\pi:=f$, $\pi_n:=f_n$, $\pi^t:=f\big|_{[0,t]}$, and $\pi^t_n:=f_n\big|_{[0,t]}$, which we view as elements of the path space $\Pi(\Xc)$. We first prove the implication $\volgt$. By Theorem~\ref{T:ArzAsc2}, $d_{\rm tot}(\pi_n,\pi)\to 0$ implies that $\{\pi_n:n\in\N\}$ is Skorokhod-equicontinuous and satisfies the compact containment condition, which implies the same is true for $\{\pi^t_n:n\in\N\}$ for any $t\in T$. Therefore, by Theorem~\ref{T:ArzAsc2}, it suffices to show that all subsequential limits of the $\pi^t_n$ are equal to $\pi^t$. Since by Theorem~\ref{T:partot}, convergence in $d_{\rm tot}$ implies convergence in $d_{\rm H}$, we can use what we have already proved for $d_{\rm H}$ to draw the desired conclusion. The implication $\Leftarrow$ follows in the same way, where now we use that if $\{\pi^t_n:n\in\N\}$ is Skorokhod-equicontinuous and satisfies the compact containment condition for any $t\in T$, then the same is true for $\{\pi_n:n\in\N\}$.
\epro

\bpro[of Theorems \ref{T:J1} and \ref{T:M1}]
Most statements of the theorems have already been proved in Propositions \ref{P:M1} and \ref{P:Pic} as well as Lemma~\ref{L:connectp}, while the statements about Polishness follow from Propositions \ref{P:PcPol}, \ref{P:J1M1Pol}, and \ref{P:DIPol}. Let $I$ be a closed real interval of positive length. The fact that $d_{\rm H}$ is a metric on $\Di_I(\Xc)$ has already been observed at the beginning of Subsection~\ref{S:Pimet}. If $I$ is a compact interval, then as discussed below Lemma~\ref{L:halfline} it follows from Proposition~\ref{P:curve} that for the trivial betweenness $d_{\rm tot}$ generates Skorokhod's J1 topology on $\Di_I(\Xc)$ as defined in \cite[Def~2.2.2]{Sko56} while $d_{\rm H}$ generates the J2 topology as defined in \cite[Def~2.2.3]{Sko56}. If $\Xc=\R$ equipped with the linear betweenness, then similarly one obtains the M1 and M2 topologies from \cite[Defs~2.2.4 and 2.2.6]{Sko56}. On half infinite intervals $I$ one moreover needs Lemma~\ref{L:restr} to see that our approach leads to the same topology as the classical definitions. If $I=\R$, then the argument is basically the same.
\epro

\bpro[of Lemma~\ref{L:halfline}]
This is of course well known, but alternatively it also follows from our Lemma~\ref{L:restr}.
\epro

\subsubsection*{Acknowledgments}

We thank Jan Seidler for useful discussions and for his help in understanding \cite{AU29}, and Petr Lachout for drawing our attention to \cite{BW71}. We thank the referee for bringing the references \cite{Fre03,Jak96,PM15} as well as example~95 in \cite{SS78} to our attention.


\begin{thebibliography}{DGW04}

\bibitem[Ald78]{Ald78}
D.~Aldous.
Stopping times and tightness.
\emph{Ann.\ Probab.}~6 (1978), 335--340. 

\bibitem[AU29]{AU29}
P.~Alexandroff and P.~Urysohn.
M\'emoire sur les espaces topologiques compacts, d\'edi\'e \`a Monsieur
D.~Egoroff.
\emph{Verhandelingen Amsterdam}~14(1) (1929).

\bibitem[BGS15]{BGS15}
N.~Berestycki, C.~Garban, and A.~Sen.
Coalescing Brownian flows: A new approach.
\emph{The Annals of Probability} 43(6) (2015), 3177–3215

\bibitem[Bil99]{Bil99}
P.~Billingsley.
\emph{Convergence of Probability Measures. 2nd ed.}
John Wiley \& Sons, 1999.

\bibitem[Bou58]{Bou58}
N.~Bourbaki.
\emph{\'El\'ements de Math\'ematique. VIII. Part. 1: Les Structures
  Fondamentales de l'Analyse. Livre III: Topologie G\'en\'erale. Chap. 9:
  Utilisation des Nombres R\'eels en Topologie G\'en\'erale. 2i\'eme \'ed.}
Actualit\'es Scientifiques et Industrielles~1045. Hermann \& Cie,
Paris, 1958.

\bibitem[BW71]{BW71}
P.J.~Bickel and M.J.~Wichura.
Convergence criteria for multiparameter stochastic processes and some
applications.
\emph{Ann. Math. Statist.}~42 (1971), 1656--1670.

\bibitem[DT96]{DT96}
A.W.M.~Dress and W.F.~Terhalle.
The real tree.
\emph{Adv.\ Math.}~120 (1996), 283--301.

\bibitem[EFS17]{EFS17}
A.~Etheridge, N.~Freeman, and D.~Straulino,
The Brownian net and selection in the spatial $\Lambda$-Fleming-Viot process.
\emph{Electronic Journal of Probability} 22 (2017).

\bibitem[EK86]{EK86}
S.N.~Ethier and T.G.~Kurtz.
\emph{Markov processes characterization and convergence}.
John Wiley \& Sons, 1986.

\bibitem[FINR04]{FINR04}
L.R.G.~Fontes, M.~Isopi, C.M.~Newman, and K.~Ravishankar.
The Brownian web: characterization and convergence.
\emph{Ann.\ Probab.}~32(4), 2857--2883, 2004.

\bibitem[Fre03]{Fre03}
D.H.~Fremlin.
\emph{Measure Theory. Volume~4.}
Torres Fremlin, Colchester, 2003.

\bibitem[FS23]{FS23}
N.~Freeman and J.~M.~Swart.
Weaves, webs and flows.
\emph{Electron.\ J.\ Probab.}~29 (2024), 1--82.

\bibitem[FS25]{FS25}
N.~Freeman and J.~M.~Swart.
Tightness criteria for random compact sets of cadlag paths.
Preprint (2025), arXiv:2501.06930.

\bibitem[Jak96]{Jak96}
A.~Jakubowski.
Convergence in various topologies for stochastic integrals driven by semimartingales.
\emph{Ann.\ Probab.}~24(4) (1996), 2141--2153. 

\bibitem[Kol56]{Kol56}
A.N.~Kolmogorov.
On Skorohod convergence.
\emph{Theor.\ Probability Appl.}~1 (1956), 213--222.

\bibitem[MRV19]{MRV19}
T.~Mountford, K.~Ravishankar, and G.~Valle
\emph{Lat. Am. J. Probab. Math. Stat.} 16 (2019), 787–-807.

\bibitem[NRS05]{NRS05}
C.~Newman, K.~Ravishankar, and R.~Sun.
Convergence of Coalescing Nonsimple Random Walks to The Brownian Web.
\emph{Electronic Journal of Probability} 10 (2005), 21--60.

\bibitem[Neu71]{Neu71}
G.~Neuhaus.
On weak convergence of stochastic processes with multidimensional time
parameter.
\emph{Ann.\ Math.\ Statist.}~42 (1971), 1285--1295.

\bibitem[Oxt80]{Oxt80}
J.C.~Oxtoby.
\emph{Measure and Category. Second Edition.}
Springer, New York, 1980.

\bibitem[Pom76]{Pom76}
J-M.L.~Pomarede.
A unified approach to Skorohod's topologies on the function space D.
\emph{Ph.D.~thesis}, Yale University, (1976).

\bibitem[PM15]{PM15}
I.~Pavlyukevich and M.~Riedle.
Non-standard Skorokhod convergence of Lévy-driven convolution integrals in Hilbert spaces.
\emph{Stochastic Anal.\ Appl.}~33(2) (2015), 271--305. 

\bibitem[Rog89]{Rog89}
L.C.G.~Rogers.
A guided tour through excursions.
\emph{Bull.\ Lond.\ Math.\ Soc.}~21(4) (1989), 305--341. 

\bibitem[Sko56]{Sko56}
A.V.~Skorohod.
Limit theorems for stochastic processes.
\emph{Theor.\ Probability Appl.}~1 (1956), 261--290.


\bibitem[SS78]{SS78}
L.A.~Steen and J.A.~Seebach Jr.
\emph{Counterexamples in Topology.} 2nd edition.
Springer, New York, 1978.


\bibitem[SSS14]{SSS14}
E.~Schertzer, R.~Sun, and J.M.~Swart.
Stochastic flows in the Brownian web and net.
\emph{Mem.\ Am.\ Math.\ Soc.} Vol.~227 (2014), Nr.~1065.

\bibitem[SSS17]{SSS17}
E.~Schertzer, R.~Sun, and J.M.~Swart.
The Brownian web, the Brownian net, and their universality.
Pages 270--368 in: P.~Contucci and C.~Giardin\`a (Eds.)
\emph{Advances in Disordered Systems, Random Processes and Some Applications.}
Cambridge University Press, 2017.

\bibitem[Whi02]{Whi02}
W.~Whitt.
\emph{Stochastic-Process Limits.
An introduction to stochastic-process limits and their application to queues.}
Springer, New York, 2002.

\end{thebibliography}
\end{document}